\documentclass[11pt]{article}
\usepackage{amsfonts}
\usepackage{amsfonts,latexsym,amsmath,amscd,geometry}
\usepackage{amssymb}
\usepackage{latexsym}
\renewcommand{\theequation}{\thesection.\arabic{equation}}
\newcommand \nc{\newcommand}
\newtheorem{theorem}{Theorem}[section]
\newtheorem{lemma}[theorem]{Lemma}
\newtheorem{proposition}[theorem]{Proposition}

\newtheorem{remark}[theorem]{Remark}
\renewcommand{\thetheorem}{\thesubsection.\arabic{theorem}}
\nc{\ba}{\begin{array}}\nc{\ea}{\end{array}}
\nc{\be}{\begin{eqnarray}}\nc{\ee}{\end{eqnarray}}
\nc{\beq}{\begin{equation}}\nc{\eeq}{\end{equation}}
\nc{\bex}{\begin{eqnarray*}}\nc{\eex}{\end{eqnarray*}}
\nc{\btm}{\begin{theorem}} \nc{\etm}{\end{theorem}}
\nc{\blm}{\begin{lemma}} \nc{\elm}{\end{lemma}}
\nc{\R}{\mathbb{R}} \nc{\va}{\varepsilon} \nc{\ls}{\limits} \nc{\wi}{\widetilde} \nc{\ga}{\gamma}

\def\pf{\noindent{\bf Proof.\quad}}\def\endpf{\hfill$\Box$}
\def\u{\dot{u}}\def\di{\mbox{div\,}}
\def\curl{\mbox{curl\,}}

\hfuzz=\maxdimen
\allowdisplaybreaks
\allowdisplaybreaks
\begin{document}
\title{{\bf Global classical solution to 3D isentropic compressible Navier-Stokes equations with large initial data and vacuum}}
\author{Xiaofeng Hou\thanks{School
of Mathematics and Statistics, Central China Normal University,
Wuhan 430079, China. Email: xiaofengh0513@163.com, }, \quad Hongyun
Peng\thanks{School of Mathematics and Statistics, Central China
Normal University, Wuhan 430079, China. Email: penghy010@163.com, },
\quad Changjiang Zhu\thanks{Corresponding author. School of
Mathematics, South China University of Technology, Guangzhou, {
510641}, China. Email: cjzhu@mail.ccnu.edu.cn. }  }
\date{}

\maketitle

\begin{abstract}
In this paper, we investigate the existence of a global classical
solution to 3D Cauchy problem of the isentropic compressible
Navier-Stokes equations with large initial data and vacuum.
Precisely, when the far-field density is vacuum ($\wi{\rho}=0$), we
get the global classical solution under the assumption that
$(\gamma-1)^\frac{1}{3}E_0\mu^{-1}$ is suitably small. In the case
that the far-field density is away from vacuum ($\wi{\rho}>0$), the
global classical solution is also obtained when
$\left((\gamma-1)^\frac{1}{36}+\wi{\rho}^\frac{1}{6}\right)E_0^{\frac{1}{4}}\mu^{-\frac{1}{3}}$
is suitably small. The above results show that the initial energy
$E_0$ could be large if $\gamma-1$ and $\wi{\rho}$ are small or the
viscosity coefficient $\mu$ is taken to be large. These results
improve the one obtained by Huang-Li-Xin in \cite{Huang-Li-Xin},
where the existence of the classical solution is proved with small
initial energy. It should be noted that in the theorems obtained in
this paper, no smallness restriction is put upon the initial data.
It can be viewed the first result on the existence of the global
classical solution to three-dimensional Navier-Stokes equations with
large initial energy and vacuum when $\gamma$ is near $1$.
\end{abstract}


\tableofcontents

\vspace{4mm}
\section {Introduction}
 \setcounter{equation}{0}\setcounter{theorem}{0}
\renewcommand{\theequation}{\thesection.\arabic{equation}}
\renewcommand{\thetheorem}{\thesection.\arabic{theorem}}

 In this paper, we consider the following isentropic compressible Navier-Stokes system in three-dimensional space
\be\label{3d-full N-S}
\begin{cases}
\rho_t+\mathrm{div}(\rho u)=0, \\
(\rho u)_t+\mathrm{div}(\rho u\otimes u)+\nabla
P=\mu\Delta u+(\mu+\lambda)\nabla\mathrm{div}u, \ \ x\in \mathbb{R}^3,\ \ t>0,\\
\end{cases}
\ee
 with the initial
condition \be\label{3d-initial-0} (\rho, u)|_{t=0}=(\rho_0,
u_0)(x),\ x\in\mathbb{R}^3, \ee
and the far-field behavior \be\label{3d-boundary-0}
\rho(x,t)\rightarrow \wi{\rho}\geq0,\ u(x,t)\rightarrow0,\
\ \mathrm{as}\ |x|\rightarrow\infty,\
\mathrm{for}\ t\ge0. \ee
Here $\rho=\rho(x,t)$ and $u=u(x, t)=(u_1, u_2, u_3)(x,t)$
represent the density and velocity of the fluid respectively; the pressure $P$ is given by
$$P(\rho)=A\rho^{\gamma},$$ where $\gamma>1$ is the adiabatic exponent, $A>0$ is a constant.
The constant viscosity coefficients $\mu$ and $\lambda$ satisfy the following physical restrictions
\be\label{viscosity assumption} \mu>0,\ \lambda+\frac{2\mu}{3}\ge0.
\ee

A great number of works have been devoted to the well-posedness of
solutions to Navier-Stokes equations. The one-dimensional problem
has been studied in many papers, for example \cite{Zhang-Fang 2006,
Kazhikhov-Shelukhi, Liutai, Mellet, Tani, Zhu2003jde, yang-zhao, Zhu
CMP, Zhang-Fang 2010} and so on. For the multi-dimensional case, in
the absence of vacuum, the local existence and uniqueness of
classical solutions are known in \cite{Nash, serrin}. Matsumura and
Nishida in \cite{Matsumura-Nishida: CMP} first proved the global
existence and uniqueness
 of classical solutions for the initial data close to a non-vacuum equilibrium in some Sobolev space $H^s$.
Later, the global existence of  weak solutions was proved by Hoff
\cite {Hoff-jde, Hoff-arma} for the discontinuous initial data with
small energy. In the presence of vacuum, due to the degeneration of
momentum equation, the issue becomes much more challenging. The
global existence of weak solutions with large initial data in $\R^N$
was first obtained by Lions in \cite{Lions2} for
$\gamma\geq\frac{3N}{N+2}\ (N=2, 3)$, where the initial energy is
finite, so that the density vanishes at far fields or even has
compact support. Later, E. Feireisl, A. Novotny and H. Petzeltov¨¢
in \cite{Feireisl2} extended Lions's result to the case
$\gamma>\frac{3}{2}$ for $N=3$, which include the monoatomic gas
$\gamma =\frac{5}{3}$. Jiang and Zhang in \cite{Jiang-Zhang-cmp,
Jiang-Zhang-jmpa} proved the global existence of weak solutions with
vacuum for any $\gamma>1$ for spherical symmetry or axisymmetric
initial data. However, the regularity and uniqueness of weak
solutions are basically open in general. Recently, when the
far-field density is away from vacuum ($\wi{\rho}>0$) and the
viscosity coefficients $\mu$ and $\lambda$ satisfy the assumption
that $\mu>\max\{4\lambda, -\lambda\}$,  Hoff and associates in
\cite{Hoff-jmfm, Hoff-arma-2, Hoff-jde-2} obtained a new type of
global weak solutions with small energy, which have extra regularity
compared with those large weak ones constructed by Lions
(\cite{Lions2}) and Feireisl et al. (\cite{Feireisl2}).

In spite of the huge amount of work, it is still a major open
problem whether or nor global (strong) classical solutions exist in
three space dimensions for general initial data with vacuum. The
local existence and uniqueness of (strong) classical solutions with
vacuum are known in \cite{Cho-Choe-Kim,Cho-Kim, cho-Kim: perfect
gas}. It seems that one should not expect better regularities of the
global solutions in general duo to Xin's results (\cite{Xin}) and
Rozanova's results (\cite{Rozanova}). It was proved that there is no
global smooth solution in $C^1([0,\infty);H^m(\mathbb{R}^d))$
($m>[\frac{d}{2}]+2$) to the Cauchy problem of the full compressible
Navier-Stokes system, if the initial density is nontrivial compactly
supported (\cite{Xin}), or the solutions are highly decreasing at
infinity (\cite{Rozanova}). Xin and Yan in \cite{Xin-Yan} improved
the blow-up results in \cite{Xin} by removing the assumptions that
the initial density has compact support and the smooth solution has
finite energy.

More recently, Huang-Li-Xin in \cite{Huang-Li-Xin} established the
surprising global existence and uniqueness of classical solutions
with constant state as far field which could be either vacuum or
nonvacuum to 3D isentropic compressible Navier-Stokes equations with
small total energy but possibly large oscillations. Then a natural
question of importance and interest is how to get a global classical
solution for large initial energy. In this paper, we devote
ourselves to this problem and show that, such a large classical
solution to (\ref{3d-full N-S})-(\ref{3d-boundary-0}) indeed exists
globally, when $\gamma$ is near $1$ or $\mu$ is taken to be large.

Before stating our main results, we would like to give some
notations which will be used throughout this paper.

\vspace{3mm}

\noindent{\bf Notations:}\\

(i)\ $\displaystyle\int_{\mathbb{R}^3} f =\int_{\mathbb{R}^3} f \,dx$, \ \  $\displaystyle\int_0^T g =\int_0^T g\,dt$.\\

(ii)\ For $1\le l\le \infty$, denote the $L^l$ spaces and the
standard Sobolev spaces as follows:
\be
&&L^l=L^l(\mathbb{R}^3),  \ D^{k,l}=\left\{ u\in L^1_{\rm{loc}}(\mathbb{R}^3): \|\nabla^k u \|_{L^l}<\infty\right\},\ D^k=D^{k,2},\nonumber\\[2mm]
&&W^{k,l}=L^l\cap D^{k,l},  \ H^k=W^{k,2}, \ D_0^1=\Big\{u\in L^6: \ \|\nabla u\|_{L^2}<\infty\Big\},\nonumber\\[2mm]
&&\dot{H}^{\beta}=\Big\{u: \mathbb{R}^3\rightarrow \mathbb{R}, \ \|u\|_{\dot{H}^{\beta}}^2=\int|\xi|^{2\beta}|\widehat{u}(\xi)|^2d \xi<\infty\Big\}.\nonumber
\ee

(iii)\ $G=(2\mu+\lambda)\mathrm{div}u-P$ is the effective viscous flux.\\

(iv)\ $\omega=\nabla\times u$ is the vorticity.\\

(v)\ $\dot{h}=h_t+u\cdot\nabla h$ denotes the material
derivative.\\

(vi)\ $\sigma =\sigma(t)=\min\{1, t\}.$\\

(vii)\ $\displaystyle
E_0=\int_{\mathbb{R}^3}\Big(\frac{1}{2}\rho_0|u_0|^2+G(\rho_0)\Big)$
is the initial energy, where $G$ denotes the potential energy
density given by
$$
\displaystyle G(\rho)\triangleq \rho\int_{\widetilde{\rho}}^{\rho}\frac{P(s)-P(\widetilde{\rho})}{s^2}ds.
$$
It is clear that
$$
\begin{cases}
\displaystyle G(\rho)=\frac{1}{\gamma-1}P \ \ \ \ \ \ \ \ \ \ \ \ \ \ \ \ \ \ \ \ \ \ \ \ \ \ \ \ \ \ \ \ \ \ \ &\rm{if} \  \widetilde{\rho}=0, \\[2mm]
\displaystyle \frac{1}{c(\bar{\rho}, \widetilde{\rho})}(\rho-\widetilde{\rho})^2\leq G(\rho)\leq c(\bar{\rho}, \widetilde{\rho})(\rho-\widetilde{\rho})^2 \ \ \ \  &\rm{if} \ \widetilde{\rho}>0,\
0\leq\rho\leq\bar{\rho},
\end{cases}
$$
for some positive constant $c(\bar{\rho}, \widetilde{\rho})$.

Now we state our main results. One of our main results is the following global existence to (\ref{3d-full N-S})-(\ref{3d-boundary-0}) with vacuum
 at infinite ($\wi{\rho}=0$).

\begin{theorem}\label{3d-th:1.1} For any given $M>0$ (not necessarily
small) and $\bar{\rho}\geq1$, assume that the initial data
$(\rho_0,u_0)$ satisfy \be\label{initial
data}\frac{1}{2}\rho_0|u_0|^2+\frac{A}{\gamma-1}\rho_0^\gamma\in L^1,\ \ u_0\in D^1\cap D^3
,\ \  (\rho_0, P(\rho_0))\in H^3,
 \ee
\be\label{3d-initial assumption} 0\le\rho_0\le\bar{\rho},\ \
\|\nabla u_0\|^2_{L^2}\le M  \ee and the compatibility condition
\be\label{3d-compatibility} -\mu\Delta
u_0-(\mu+\lambda)\nabla\mathrm{div}u_0+\nabla P(\rho_0)=\rho_0 g,
\ee where $g\in D^1$ and $\rho^{\frac{1}{2}}g\in L^2$. Then there
exists a unique global classical solution $(\rho,u)$ in
$\mathbb{R}^3\times[0,\infty)$ satisfying, for any
$0<\tau<T<\infty$,
\be\label{3d-th1.1-2} 0\le\rho\le2\bar{\rho},\ \
\ \  x\in\mathbb{R}^3,\ \ t\geq0, \ee \be\label{3d-th1.1-3}
\begin{cases}
(\rho, P)\in C([0,T]; H^3), \ \ \ \sqrt{\rho}u_t\in  L^{\infty}(0,T; L^2),\\[3mm]
u\in  C([0,T]; D^1\cap D^3)\cap L^2(0, T; D^4)\cap L^{\infty}(\tau, T; D^4),\\[3mm]
u_t\in  L^{\infty}(0,T; D^1)\cap L^2(0, T; D^2)\cap L^{\infty}(\tau, T; D^2)\cap H^1(\tau, T; D^1),
\end{cases}
\ee provided that\be\label{3d-th1.1-smallcondition}
\frac{(\gamma-1)^\frac{1}{3}E_0}{\mu}\le\varepsilon\triangleq
\min\left\{\varepsilon^2_3,\ (2C(\bar{\rho},M))^{-\frac{32}{3}}\mu^{8},(4C(\bar{\rho}))^{-4}\right\},
\ee where \bex\begin{aligned}
&\varepsilon_3=\min\left\{(CE_7)^{-3}\Big|_{(1< \gamma\leq
\frac{3}{2})}, (CE_{11})^{-2}\Big|_{( \gamma>\frac{3}{2})},
 \varepsilon_2
\right\},\\[3mm]
&\varepsilon_2=\min \left\{
C(\bar{\rho})^{-2}(\gamma-1)^{-\frac{2}{3}}E_2^{-3}\mu^5\Big|_{(1 <
\gamma\leq \frac{3}{2})},\
C(\bar{\rho})^{-1}\mu^{\frac{9}{4}}E_2^{-\frac{3}{4}}\Big|_{(\gamma >
\frac{3}{2})}, \varepsilon_1
\right\},\\[3mm]
&\varepsilon_1=\min\left\{\Big(4C(\bar{\rho})\Big)^{-6}, 1 \right\}.
\end{aligned}\eex
Here, $C$ depending on $\bar{\rho}, M$ and some other known
constants but independent of $\mu,\lambda, \gamma-1$ and $t$ $(
see \  \eqref{a1}, \eqref{a2})$. $E_2$, $E_7$ and $E_{11}$ are defined by
\eqref{e2}, \eqref{3d-dt-E-8} and \eqref{3d-dt-E-13} respectively.
\end{theorem}

Now we briefly outline the main ideas of the proof of Theorem 1.1,
some of which are inspired by \cite{Hoff-jde} and
\cite{Huang-Li-Xin}. The key point of this paper is to get the
time-independent upper bound on the density $\rho$ (see\
(\ref{3d-upper bound of rho})), and once that is obtained, the proof
of Theorem 1.1 follows in the same way as in \cite{Huang-Li-Xin}.
It's worth noting that the methods in all previous works
\cite{Hoff-jde, Huang-Li-Xin, Zhang 2011} depend crucially on the
small initial energy. Thus, some new ideas are needed to recover all
the $a \ priori$ estimates under only the assumption (\ref
{3d-th1.1-smallcondition}) while the initial energy $E_0$ could be
large (see  \eqref{3d-lemma-u-4}-\eqref{3d-dt-A-1+A-2-dd}). In fact,
the small initial energy could naturally yield the smallness of
$\int_0^{T}\int_{\mathbb{R}^3} |\nabla u|^2$, which plays a crucial
role in the analysis to prove the time-independent upper bound of
$\rho$. But the smallness of $\int_0^{T}\int_{\mathbb{R}^3} |\nabla
u|^2$  is not valid here without the small initial energy. The
crucial ideas to overcome this difficulty are as follows:

$\bullet$\ One key observation is that $A_1(T)$ could be bounded by
$\int_0^{\sigma(T)}\|\nabla u\|_{L^2}^2$ and some other terms, i.e.,
$$
A_1(T) \leq \cdots\cdots +\displaystyle
\frac{(2\mu+\lambda)}{\mu}\int_0^{\sigma(T)}\|\nabla
u\|_{L^2}^2+\frac{C(\gamma-1)E_0}{\mu^2}. \ \ \ \ \ \ \ \  {\rm (see
\  \eqref{3d-A1})}
$$
Thus, we can close the $a \ priori$ assumption on $A_1(T)$ by the
smallness of $\int_0^{\sigma(T)}\|\nabla u\|_{L^2}^2$ instead of the
smallness of $\int_0^T\|\nabla u\|_{L^2}^2$. And we can prove that
$\int_0^{\sigma(T)}\|\nabla u\|_{L^2}^2$ is small, when
$(\gamma-1)^\frac{1}{3}E_0\mu^{-1}$ is suitably small (see
(\ref{3d-basic-3})).

$\bullet$\  Another new ingredient in the proof is that we can use
the smallness of $(\gamma-1)^\frac{1}{3}E_0\mu^{-1}$ and the boundedness  of
$\int_0^T\int_{\mathbb{R}^3} |\nabla u|^2$ to estimate the
higher-order terms of $\nabla u$ as follows: \be
\displaystyle\int_0^T\int_{\mathbb{R}^3} \sigma^2|\nabla u|^4
&\le& \displaystyle CN_4^4(2\mu+\lambda)^{-3}\int_0^{T}\sigma^2\|\rho\u\|_{L^2}^3\|\nabla u\|_{L^2}+\cdots\cdots\nonumber\\[2mm]
&\le&\displaystyle CN_4^4(2\mu+\lambda)^{-3}\mu^{\frac{1}{2}}A^{\frac{1}{2}}_2(T)\sup\limits_{0\le t\le T}(\|\nabla u\|_{L^2}^2)^{\frac{1}{2}}\int_0^{T}\sigma\|\rho\|_{L^3}^2\|\nabla\u\|_{L^2}^2\nonumber\\[2mm]
&& +\cdots\cdots\nonumber\\[2mm]
&\le&(2\mu+\lambda)^{-3}\mu^\frac{1}{4}(\gamma-1)^{\frac{3}{4}}E_0^{\frac{11}{12}}E_6
\nonumber\ \ \ \ \ \ \ \ \ \ \ \ \ \ \ \ \ \ \ {\rm (see\
 \eqref{3d-lemma-u-4}-\eqref{3d-lemma-u-4-f}}) \ee
and \be
\displaystyle
\int_0^T \sigma \|\nabla u\|_{L^3}^3
&\le& \displaystyle\left(\int_0^T\int_{\mathbb{R}^3}  |\nabla u|^2\right)^{\frac{1}{2}}\left(\int_0^T\int_{\mathbb{R}^3} \sigma^2 |\nabla u|^4\right)^{\frac{1}{2}}\nonumber\\[3mm]
&\le&\displaystyle
C\frac{(\gamma-1)^{\frac{3}{8}}E_0^{\frac{23}{24}}E_6^{\frac{1}{2}}}{\mu^{\frac{3}{8}}(2\mu+\lambda)^{\frac{3}{2}}}.\nonumber
\ \ \ \ \ \ \ \ \ \ \ \ \ \ \ \ \ \ \ \ \ \ \  {\rm (see\
\eqref{3d-dt-II-4})}\ee We refer the details to Lemma
\ref{3d-le:3.8}. \vspace{1ex}

 On the other hand, from the proof of
Proposition 3.1, we know that it is important to find a suitable
match for $\mu$, $(\gamma-1)$ and $E_0$. That means much more
complicated estimates than those in \cite{Hoff-jde, Huang-Li-Xin,
Zhang 2011} are needed. To do this, we derive some more
sophisticated inequalities about $\mu$ (see Lemma 2.2). For more
details, please see the proof of Proposition 3.1.

\vspace{2mm}
Concerning  the global classical solutions for (\ref{3d-full N-S})-(\ref{3d-boundary-0}) in the case that the far-field density is away from vacuum ($\wi{\rho}>0$), we have

\begin{theorem}\label{3d-th:1.2} For any given $M>0$ (not necessarily
small) and $\bar{\rho}\geq \wi{\rho}+1$, assume that the initial
data $(\rho_0,u_0)$ satisfy \be\label{initial
data2}\frac{1}{2}\rho_0|u_0|^2+G(\rho_0)\in L^1,\ \ u_0\in H^1\cap
D^3 ,\ \  (\rho_0-\wi{\rho}, P(\rho_0)-P(\wi{\rho}))\in H^3,
 \ee
\be\label{3d-initial assumption2} 0\le\rho_0\le\bar{\rho},\ \
\|u_0\|^2_{L^2}\le E_0 ,\ \ \|\nabla u_0\|^2_{L^2}\le M  \ee and the
compatibility condition \be\label{3d-compatibility2} -\mu\Delta
u_0-(\mu+\lambda)\nabla\mathrm{div}u_0+\nabla P(\rho_0)=\rho_0 g,
\ee where $g\in D^1$ and $\rho^{\frac{1}{2}}g\in L^2$. Then there
exists a unique global classical solution $(\rho,u)$ in
$\mathbb{R}^3\times[0,\infty)$ satisfying, for any
$0<\tau<T<\infty$, \be\label{4d-th1.1-2} 0\le\rho\le2\bar{\rho},\ \
\ \   x\in\mathbb{R}^3,\ \ t\geq0, \ee \be\label{4d-th1.1-3}
\begin{cases}
(\rho-\wi{\rho}, P-P(\wi{\rho}))\in C([0,T]; H^3), \ \ \ \sqrt{\rho}u_t\in  L^{\infty}(0,T; L^2),\\[3mm]
u\in  C([0,T]; D^1\cap D^3)\cap L^2(0, T; D^4)\cap L^{\infty}(\tau, T; D^4),\\[3mm]
u_t\in  L^{\infty}(0,T; D^1)\cap L^2(0, T; D^2)\cap L^{\infty}(\tau, T; D^2)\cap H^1(\tau, T; D^1),
\end{cases}
\ee provided that
$\displaystyle\frac{(\gamma-1)^\frac{1}{36}E_0^{\frac{1}{4}}}{\mu^{\frac{1}{3}}}\le
\frac{\tilde{\rho}}{2 C}$ and \be\label{4d-th1.1-smallcondition}
\displaystyle\frac{\left((\gamma-1)^\frac{1}{36}+\wi{\rho}^{\frac{1}{6}}\right)E_0^{\frac{1}{4}}}{\mu^{\frac{1}{3}}}\le\varepsilon\triangleq
\min\left\{\varepsilon_6,(2C(\bar{\rho},M))^{-\frac{16}{3}}\mu^{4},(4C(\bar{\rho}))^{-2}\right\},
\ee where \bex\begin{aligned}
&\varepsilon_6=\min\left\{\Big(C(E_{18}+E_{19}+E_{20})\Big)^{-17},
\Big(C(E_{18}+E_{19}+E_{21})\Big)^{-8}, \varepsilon_5
\right\},\\[3mm]
&\varepsilon_5=\min \left\{
\Big(C(E_{15}E_{17}+E_{16})\Big)^{-4},\ \varepsilon_4\right\},\\[3mm]
&\varepsilon_4=\min\left\{\Big(4C(\bar{\rho})\Big)^{-6}, 1 \right\}.
\end{aligned}\eex
Here $C$ denotes a generic positive constant depending on
$\bar{\rho}, M$ and some other known constants but independent of
$\mu, \lambda, \gamma-1, \wi{\rho},$ and $t$. $E_{15}-E_{21}$ are defined by
\eqref{e15-17}, \eqref{s 4E14}, \eqref{s 4E15}, \eqref{s 4E16} and
\eqref{s 4E17}.
\end{theorem}

There are two key technical points in the proof of Theorem
\ref{3d-th:1.2}:

\vspace{2mm}

(1) The smallness of $\|\rho-\widetilde{\rho}\|_{L^2}$ plays a key
role in establishing the time-independent upper bound for the
density $\rho$. In \cite{Hoff-jde, Huang-Li-Xin, Zhang 2011}, the
smallness of $\|\rho-\wi{\rho}\|_{L^2}$ is easy to get because of
the small initial energy.  And in the proof of Theorem
\ref{3d-th:1.1} where $\wi{\rho}=0$, $\|\rho\|_{L^2}$ can be small
when $\ga$ is near 1. But, when $\wi{\rho}>0$ and the initial energy
is large, it seems impossible to get the smallness of
$\|\rho-\widetilde{\rho}\|_{L^2}$, even as $\ga\rightarrow 1$. Based
on the elaborate analysis of the potential energy density $G(\rho)$,
we succeed in deriving a new estimate to $\rho-\wi{\rho}$, i.e., \be
\begin{aligned}
G(\rho)=&\frac{1}{\gamma-1}[\rho^{\ga}-\wi{\rho}^{\ga}-\ga\wi{\rho}^{\ga-1}(\rho-\wi{\rho})]\nonumber\\[2mm]
\geq&\begin{cases}(\ga-1)^{-\frac{1}{4}}|\rho-\wi{\rho}|^{\ga-1}, \ \ \ \ \ \ &|\rho-\wi{\rho}|>(\ga-1)^{\frac{1}{3}},\nonumber\\[2mm]
(\ga-1)^{-\frac{1}{4}}|\rho-\wi{\rho}|^3, \ \ \ \ \ \ \ \ \
&|\rho-\wi{\rho}|\le(\ga-1)^{\frac{1}{3}},
\end{cases}
\end{aligned}
\ee which implies
 the $L^3$-norm of $\rho-\wi{\rho}$ can be small when $\ga\rightarrow 1$ (see Lemma \ref{3d-le-rho3}).
In fact, the $L^3$-norm of $\rho-\wi{\rho}$ is weaker than the
$L^2$-norm in the whole space, when $\rho$ is upper-bounded. Thus,
applying the smallness of $L^3$-norm of $\rho-\wi{\rho}$ to
establish the time-independent upper bound for the density $\rho$
will bring much more difficulties than that of $L^2$-norm.

\  $\bullet$\ Unlike Theorem \ref{3d-th:1.1}, the method in
obtaining the smallness of $\int_0^{\sigma(T)}\int_{\mathbb{R}^3}
|\nabla u|^2$ is no longer applicable here. To overcome this
difficulty, we have to give the estimate of $\|u\|_{L^2}$ (see
\eqref{3d-p-p-p}-\eqref{3d-p-u}). This together with some other
estimates yields the smallness of
$\int_0^{\sigma(T)}\int_{\mathbb{R}^3} |\nabla u|^2$ (see Lemma
\ref{4.2}). \vspace{1ex}

\ $\bullet$\ All the estimates containing $P-P(\widetilde{\rho})$ in
the proof of Theorem \ref{3d-th:1.2} can not be estimated as the
previous way. Hence, some new ideas are needed to recover all the
{\it a priori} estimates under the condition that
$\|\rho-\wi{\rho}\|_{L^3}$ is small.


\vspace{2mm} 
(2) Due to losing the smallness of $\int_{0}^T\int_{\mathbb{R}^3}
|\nabla u|^2$ and $\|\rho-\widetilde{\rho}\|_{L^2}$, the techniques
in \cite{Hoff-jde, Huang-Li-Xin, Zhang 2011} which are used to close
the $a \ priori$ assumptions on $A_1(T)$ and $A_2(T)$ fail in this
paper. From \eqref{four A one}\ (see $N_{11}$) \be A_{1}(T) \le
\cdots
+\frac{CE^\frac{1}{2}_0}{\mu^2}A^\frac{3}{4}_{1}(T)A^\frac{1}{4}_{2}(T)+\cdots,\nonumber
\ee we find that, in order to close the $a \ priori$ assumption on
$A_1(T)$, the bound of $A_2(T)$ should be given a higher order than
that of $A_{1}(T)$ \ (Roughly speaking, we can choose $A_2(T)\sim
A_1(T)^{\frac{8}{3}} $ in the present paper to close the $a \
priori$ assumptions on $A_1(T)$ and $A_2(T)$). This means we could
not close them simultaneously. Fortunately, we observe that $A_2(T)$
could be bounded by the boundedness of $A_1(\sigma(T))$ and some
other terms, i.e., \be
 \displaystyle A_2(T)&\leq& C A_1(\sigma(T))
+\cdots\cdots +\frac{CP(\wi{\rho})^2}{\mu^2}E_0.\nonumber\ \ \ \ {\rm
(see \ \eqref{4d-A2-c})} \ee And we can give a better estimate of
$A_1(\sigma(T))$ compared with $A_1(T)$'s. This is based on that
$\int_0^{\sigma(T)} \sigma \|\nabla u\|_{L^3}^3$ has a stronger
control than $\int_{\sigma(T)}^T \sigma \|\nabla u\|_{L^3}^3$ (see
\eqref{4A1-11} and \eqref{4A1t-12}). Therefore, to close the $a \
priori$ assumptions on $A_1(T)$ and $A_2(T)$, the first step is to
estimate $A_1(\sigma(T))$ (see \eqref{4d-A1t-c1}). Next, we can
bound $A_2(T)$ in \eqref{4a2}. Finally, the estimate of $A_1(T)$
is obtained in \eqref{four A one}.
\begin{remark}
This should be the first result concerning the global existence and
uniqueness of classical solutions to the Cauchy problem for the
three dimensional isentropic compressible Navier-Stokes equations
with large initial energy and vacuum. Compared to these results in
\cite{Hoff-jde, Huang-Li-Xin}, the thrust of this paper is to remove
the condition of smallness on the initial energy.
\end{remark}

\begin{remark}
The results in this paper generalize the ones in
\cite{Huang-Li-Xin}. More accurately,
in the case of $\wi{\rho}=0$ and $\wi{\rho}>0$, the existence of
classical solutions to (\ref{3d-full N-S})-(\ref{3d-boundary-0}) are
obtained respectively; in particular, the initial energy is allowed
to be large when $\gamma-1$ and $\wi{\rho}$ are suitable small or
$\mu$ is taken to be large. On the other hand, Theorems 1.1 and 1.2
are still applicable to the case that initial energy $E_0$ is small
for any given $\gamma$, $\wi{\rho}$ and $\mu$.
\end{remark}

\begin{remark}
It is well known that Nishida and Smoller in \cite{Nishida-Smoller:
CPAM} proved the Cauchy problem for 1D isentropic Euler equations
has a global solution provided that $(\gamma-1)\rm{T.V.}\{ \it u_0,
\rho_0\}$ is sufficiently small. This means that when $\gamma$ is
near $1$, one can allow large data, and conversely, as $\gamma$
increases, one must take correspondingly smaller data. Recently,
Tan-Yang-Zhao-Zou in \cite{Tan 2013} obtained a global smooth
solution to the one-dimensional compressible Navier-Stokes-Poisson
equations with large initial data under the assumption that $\gamma$
is near $1$. Soon later, a similar result on 1D compressible
Navier-Stokes equations was obtained by Liu-Yang-Zhao-Zou in
\cite{Liu 2014}. These works inspire us to look for the large
solution to (\ref{3d-full N-S})-(\ref{3d-boundary-0}) when $\gamma$
is near $1$.
\end{remark}

\begin{remark}
From physical viewpoint, it is very nature to get a global large
solution to (\ref{3d-full N-S})-(\ref{3d-boundary-0}), when the
viscosity coefficient is sufficiently large. Note that the
coefficients of viscosity are only required to satisfy the physical
restriction (\ref{viscosity assumption}) in the present paper.
\end{remark}

\begin{remark}
Though the initial data can be large if the adiabatic exponent
$\gamma$ goes to $1$ or the viscosity coefficient $\mu$ is taken to
be large, it is still unknown whether the global classical solution
exists when the initial data is large for any fixed $\gamma$ and
$\mu$. It should be noted that the similar question of whether there
exists a global smooth solution of the three-dimensional
incompressible Navier-Stokes equations with smooth initial data is
one of the most outstanding mathematical open problems
(\cite{fefferman}). Motivated by this, some blow-up criterions of
strong\ (classical) solutions to (\ref {3d-full N-S}) have been
studied, please refer for instance to \cite{Huang-Xin:Sci,
Huang-Li-Xin:SIAM, Sun-Wang-Zhang, Wen-Zhu 4} and references
therein. In fact, for initial-boundary-value problems or periodic
problems of compressible Naiver-Stokes equations with vacuum in one
dimension, or in two dimensions for isentropic flow, or in higher
dimensions with symmetric initial data, the existence of global
large regular solutions has been obtained, please refer to
\cite{Ding-Wen-Zhu, Ding-Wen-Yao-Zhu, Jiu-Wang-Xin, Wen-Zhu, Wen-Zhu
3} and references therein.
\end{remark}

\begin{remark}
In addition to the conditions of Theorem \ref{3d-th:1.1} and  Theorem \ref{3d-th:1.2}, if we assume further that $u_0\in \dot{H}^{\beta}(\beta\in(\frac{1}{2}, 1])$ and replace $\|\nabla u_0\|^2_{L^2}\le M$ with $\|u_0\|_{\dot{H}^{\beta}}\leq \widetilde{M}$, the conclusions in Theorem \ref{3d-th:1.1} and Theorem \ref{3d-th:1.2} will still hold, and the $\varepsilon$ will also depend on $\widetilde{M}$ instead of $M$ correspondingly. This can be achieved by a similar way as in \cite{Huang-Li-Xin}.
\end{remark}


\vspace{1.5mm}

The rest of the paper is organized as follows. In Section 2, we collect some known inequalities and facts which will be frequently used later.
In Section 3, we obtain the prove of Theorem 1.1. Then the proof of Theorem 1.2 is completed in Section 4.

\section{Preliminaries}

 \setcounter{equation}{0}\setcounter{theorem}{0}
\renewcommand{\theequation}{\thesection.\arabic{equation}}
\renewcommand{\thetheorem}{\thesection.\arabic{theorem}}
If the solutions are regular enough (such as strong solutions
and classical solutions), (\ref{3d-full N-S}) is equivalent to the
following system
 \be\label{full N-S+1}
\begin{cases}
\rho_t+\nabla \cdot (\rho u)=0, \\
\rho  u_t+\rho u\cdot\nabla u+\nabla
P(\rho)=\mu\Delta u+(\mu+\lambda)\nabla\mathrm{div}u.
\end{cases}
\ee
System (\ref{full N-S+1}) is supplemented with initial
condition \be\label{3d-initial} (\rho, u)|_{t=0}=(\rho_0,
u_0)(x),\ x\in\mathbb{R}^3, \ee
and the far-field behavior \be\label{3d-boundary}
\rho(x,t)\rightarrow \wi{\rho}\geq 0,\ u(x,t)\rightarrow0,\
\ \mathrm{as}\ |x|\rightarrow\infty,\
\mathrm{for}\ t\ge0. \ee

Since the exact value of $A$ in the pressure $P$ doesn't play a role
in this paper, we henceforth assume $A=1$. Next, we will list
several facts which will be used in the proof of the main results.
The first one is the well-known Gagliardo-Nirenberg inequality (see
\cite {LO}).
\begin{lemma}\label{GN}
 For any $p \in[2,6], q\in (1,\infty)$ and $r\in(3,\infty)$, there exists some generic
 constant $C>0$ that may depend on $q$ and $r$ such that for $f\in H^1(\mathbb{R}^3)$ and $g\in L^q(\mathbb{R}^3)\cap D^{1,r}(\mathbb{R}^3)$, we have
 \begin{equation}\label{GN1}
 \|f\|^p_{L^p(\mathbb{R}^3)}\leq C\|f\|^{\frac{6-p}{2}}_{L^{2}(\mathbb{R}^3)}\|\nabla f\|^{\frac{3p-6}{2}}_{L^{2}(\mathbb{R}^3)},\\
 \end{equation}
\begin{equation}\label{GN2}
 \|g\|_{C(\mathbb{R}^3)}\leq C\|g\|^{\frac{q(r-3)}{3r+q(r-3)}}_{L^{q}(\mathbb{R}^3)}\|\nabla g\|^{\frac{3r}{3r+q(r-3)}}_{L^{r}(\mathbb{R}^3)}.
\end{equation}
\end{lemma}

We now state some elementary estimates that follow from \eqref{GN1} and the standard $L^p$-estimate for the following elliptic system derived from the
 momentum equation in \eqref{3d-full N-S}:
\be\label{elp1}
 \Delta G= \di(\rho \dot{u}), \ \ \ \ \ \mu\Delta\omega=\nabla\times(\rho \dot{u}).
 \ee

\begin{lemma}\label{2.2}
Let $(\rho,u)$ be a smooth solution of \eqref{full N-S+1}, \eqref{3d-boundary}. Then there exists a generic positive constant $C$ such that for any $p\in [2,6]$
 \be
&\|\nabla G\|_{L^p} \leq& C \|\rho \dot{u}\|_{L^p},\ \ \   \|\nabla \omega\|_{L^p}\le \frac{C}{\mu}\|\rho \dot{u}\|_{L^p}, \label{hp4}\\[1.5mm]
&\|G\|_{L^p}
\leq&C \Big((2\mu+\lambda)\|\nabla u\|_{L^2}+\|P-P(\wi{\rho})\|_{L^2}\Big)^{\frac{6-p}{2p}}\|\rho\dot{u}\|_{L^2}^{\frac{3p-6}{2p}},\label{hp1}\\[1.5mm]
&\|G\|_{L^p}
\leq&C \Big((2\mu+\lambda)\|\nabla u\|_{L^3}+\|P-P(\wi{\rho})\|_{L^3}\Big)^{\frac{6-p}{p}}\|\rho\dot{u}\|_{L^2}^{\frac{2p-6}{p}},\label{hp5}\\[1.5mm]
&\|w\|_{L^p}
\leq&C\Big(\frac{1}{\mu}\Big)^{\frac{3p-6}{2p}}\|\nabla u\|_{L^2}^{\frac{6-p}{2p}}\|\rho\dot{u}\|_{L^2}^{\frac{3p-6}{2p}},\label{hp2}\\[1.5mm]
&\|\nabla u\|_{L^p}
\leq&C\Big((2\mu+\lambda)\|\nabla u\|_{L^3}+\|P-P(\wi{\rho})\|_{L^3}\Big)^{\frac{6-p}{p}}\|\rho\dot{u}\|_{L^2}^{\frac{2p-6}{p}}\nonumber\\[1.5mm]
&&+C\|P-P(\wi{\rho})\|_{L^p}+C\|w\|_{L^p}\label{hp6}
\ee
and
\be\label{hp3}
\displaystyle\|\nabla u\|_{L^p}&\leq&\displaystyle
CN_p(2\mu+\lambda)^{\frac{6-3p}{2p}}\|\nabla u\|_{L^2}^{\frac{6-p}{2p}}\|\rho\dot{u}\|_{L^2}^{\frac{3p-6}{2p}}\nonumber\\[3mm]
&&\displaystyle+\frac{C}{2\mu+\lambda}\Big(\|\rho\dot{u}\|_{L^2}^{\frac{3p-6}{2p}}\|P-P(\wi{\rho})\|_{L^2}^{\frac{6-p}{2p}}+\|P-P(\wi{\rho})\|_{L^p}\Big),
\ee
where $\displaystyle N_p=1+\left(2+\frac{\lambda}{\mu}\right)^{\frac{3p-6}{2p}}$.
 \end{lemma}
\pf
The standard $L^p$-estimate for elliptic system \eqref{elp1} yields \eqref{hp4}.
Using \eqref{GN1}, we obtain
\be
\displaystyle\|G\|_{L^p}&\leq&\displaystyle C\|G\|_{L^2}^{\frac{6-p}{2p}}\|G\|_{L^6}^{\frac{3p-6}{2p}}
\leq C\|G\|_{L^2}^{\frac{6-p}{2p}}\|\nabla G\|_{L^2}^{\frac{3p-6}{2p}}\nonumber\\[1.5mm]
&\leq&\displaystyle C\Big((2\mu+\lambda)\|\nabla u\|_{L^2}+\|P-P(\wi{\rho})\|_{L^2}\Big)^{\frac{6-p}{2p}}\|\rho\dot{u}\|_{L^2}^{\frac{3p-6}{2p}},\\[1.5mm]
\displaystyle\|G\|_{L^p}&\leq&\displaystyle C\|G\|_{L^3}^{\frac{6-p}{p}}\|G\|_{L^6}^{\frac{2p-6}{p}}
\leq C\|G\|_{L^3}^{\frac{6-p}{p}}\|\nabla G\|_{L^2}^{\frac{2p-6}{p}}\nonumber\\[1.5mm]
&\leq&\displaystyle C\Big((2\mu+\lambda)\|\nabla u\|_{L^3}+\|P-P(\wi{\rho})\|_{L^3}\Big)^{\frac{6-p}{p}}\|\rho\dot{u}\|_{L^2}^{\frac{2p-6}{p}},\\[1.5mm]
\displaystyle\|w\|_{L^p}&\leq&C\displaystyle\|w\|_{L^2}^{\frac{6-p}{2p}}\|w\|_{L^6}^{\frac{3p-6}{2p}}
\leq C\|w\|_{L^2}^{\frac{6-p}{2p}}\|\nabla w\|_{L^2}^{\frac{3p-6}{2p}}\nonumber\\[1.5mm]
&\leq&\displaystyle C\Big(\frac{1}{\mu}\Big)^{\frac{3p-6}{2p}}\|\nabla u\|_{L^2}^{\frac{6-p}{2p}}\|\rho\dot{u}\|_{L^2}^{\frac{3p-6}{2p}}.
\ee
Note that $-\Delta u=-\nabla\di u+ \nabla\times\omega $, which implies that
$$
\nabla u=-\nabla(-\Delta)^{-1}\nabla \di u+ \nabla(-\Delta)^{-1}\nabla\times \omega.
$$
Thus the standard $L^p$-estimate shows that
$$
\displaystyle\|\nabla u\|_{L^p}\leq C\left(\|\di u\|_{L^p}+\|\omega\|_{L^p}\right) \ \  \mathrm{for}\   p \in [2,6],
$$
which together with \eqref{hp1}, \eqref{hp2} and the definition of $G$, give
\be
\displaystyle\|\nabla u\|_{L^p}&\leq& C\|\di u\|_{L^p}+C\|\curl u\|_{L^p}\nonumber\\[1.5mm]
&\leq&\displaystyle\frac{C}{2\mu+\lambda}\Big(\|G\|_{L^p}+\|P-P(\wi{\rho})\|_{L^p}\Big)+C\|w\|_{L^p}\nonumber\\[1.5mm]
&\leq&\displaystyle C\Big((2\mu+\lambda)^{\frac{6-3p}{2p}}+\mu^{\frac{6-3p}{2p}}\Big)\|\nabla u\|_{L^2}^{\frac{6-p}{2p}}\|\rho\dot{u}\|_{L^2}^{\frac{3p-6}{2p}}\nonumber\\[1.5mm]
&&+\frac{C}{2\mu+\lambda}\Big(\|\rho\dot{u}\|_{L^2}^{\frac{3p-6}{2p}}\|P-P(\wi{\rho})\|_{L^2}^{\frac{6-p}{2p}}+\|P-P(\wi{\rho})\|_{L^p}\Big)\nonumber\\[1.5mm]
&\leq&\displaystyle CN_p(2\mu+\lambda)^{\frac{6-3p}{2p}}\|\nabla u\|_{L^2}^{\frac{6-p}{2p}}\|\rho\dot{u}\|_{L^2}^{\frac{3p-6}{2p}}\nonumber\\[1.5mm]
&&+\frac{C}{2\mu+\lambda}\Big(\|\rho\dot{u}\|_{L^2}^{\frac{3p-6}{2p}}\|P-P(\wi{\rho})\|_{L^2}^{\frac{6-p}{2p}}+\|P-P(\wi{\rho})\|_{L^p}\Big).
\ee
\endpf
\begin{lemma}[\cite{Zlo}]\label{zlo}
Let the function $y$ satisfy
$$
y'(t)=g(y)+b'(t)\ \ \  \mathrm{on}\ [0,T], \ \  y(0)=y^0,
$$
with $g\in C(\mathbb{R})$
 and $ y,b \in W^{1,1}(0,T)$. If $g(\infty)=-\infty$ and
$$
 b(t_2)-b(t_1)\leq  N_0+N_1(t_2-t_1),
$$
for all $0 \leq t_1\leq t_2\leq T$ with some $N_0\geq 0$ and $N_1\geq 0$, then
$$
y(t)\leq\max\{y^0,\bar{\xi}\}+N_0<\infty \ \ \ \rm{on}\ \  [0,T],
$$
where $\bar{\xi}$ is a constant such that
$$
g(\xi)\leq -N_1,\ \ \mathrm{for}\ \ \xi \geq \bar{\xi}.
$$
\end{lemma}

\section{The proof of Theorem 1.1}

 \setcounter{equation}{0}\setcounter{theorem}{0}
\renewcommand{\theequation}{\thesection.\arabic{equation}}
\renewcommand{\thetheorem}{\thesection.\arabic{theorem}}

In this section, we will first establish the time-independent upper
bound of the density $\rho$. Assume that $(\rho, u)$ is a smooth solution
to (\ref{3d-full N-S})-(\ref{3d-boundary-0}) on $\mathbb{R}^3 \times
(0, T)$ for some positive time $T>0$. Set
$\sigma=\sigma(t)=\min\{1,t\}$ and denote \be\label{A2 1}
\begin{cases}
\displaystyle A_1(T)=\sup\limits_{0\le t\le T}\sigma\int_{\mathbb{R}^3}|\nabla
u|^2+\int_0^T\int_{\mathbb{R}^3}\sigma\frac{\rho|\u|^2}{\mu},\\[3mm]
\displaystyle A_2(T)=\sup\limits_{0\le t\le T}\sigma^2\int_{\mathbb{R}^3}\frac{\rho
|\u|^2}{\mu}+\int_0^T\int_{\mathbb{R}^3}\sigma^2|\nabla\u|^2,\\[3mm]
\displaystyle A_3(T)=\sup\limits_{0\le t\le T}\int_{\mathbb{R}^3}\frac{\rho|u|^3}{\mu^3}.
\end{cases}
\ee The following proposition plays a crucial role in this section.
\begin{proposition}\label{prop 3.1}
Assume that the initial data satisfies (\ref{initial data}),
(\ref{3d-initial assumption}) and (\ref{3d-compatibility}). If the
solution $(\rho, u)$ satisfies

\be\label{a priori assumption 1} A_1(T)+A_2(T)\le\frac{2(\gamma-1)^\frac{1}{6}E_0^{\frac{1}{2}}}{\mu^{\frac{1}{2}}},\ A_3(\sigma(T))\le
\frac{2(\gamma-1)^\frac{1}{12}E_0^{\frac{1}{4}}}{\mu^{\frac{1}{4}}},\ 0\le \rho\le 2\bar{\rho},\
\ee then for any $(x,t) \in {\mathbb{R}}^3 \times [0,T]$, the following estimates hold:

\be\label{3d-result from the a p}  A_1(T)+A_2(T)\le\frac{(\gamma-1)^\frac{1}{6}E_0^{\frac{1}{2}}}{\mu^{\frac{1}{2}}},\ A_3(\sigma(T))\le
\frac{(\gamma-1)^\frac{1}{12}E_0^{\frac{1}{4}}}{\mu^{\frac{1}{4}}},\ 0\le \rho\le \frac{7}{4}\bar{\rho},
\ee
 provided $\displaystyle\frac{(\gamma-1)^{\frac{1}{3}}E_0}{\mu}\le\varepsilon$. Here
\bex
\begin{aligned}
\varepsilon= &\min\left\{\varepsilon^2_3,\
(2C(\bar{\rho},M))^{-\frac{32}{3}}\mu^{8},\
(4C(\bar{\rho}))^{-4}\right\},\end{aligned} \eex  and
\bex\begin{aligned} &\varepsilon_3=\min\left\{(CE_7)^{-3}\Big|_{(1<
\gamma\leq \frac{3}{2})}, (CE_{11})^{-2}\Big|_{(
\gamma>\frac{3}{2})},  \varepsilon_2
\right\},\\[3mm]
&\varepsilon_2=\min \left\{
C(\bar{\rho})^{-2}(\gamma-1)^{-\frac{2}{3}}E_2^{-3}\mu^5\Big|_{(1 <
\gamma\leq \frac{3}{2})},\
C(\bar{\rho})^{-1}\mu^{\frac{9}{4}}E_2^{-\frac{3}{4}}\Big|_{(\gamma >
\frac{3}{2})}, \varepsilon_1
\right\},\\[3mm]
&\varepsilon_1=\min\left\{\Big(4C(\bar{\rho})\Big)^{-6}, 1 \right\}.
\end{aligned}\eex
Here, $C$ depending on $\bar{\rho}, M$ and some other known
constants but independent of $\mu,\lambda, \gamma-1$ and $t$ $(see
\ \eqref{a1}, \eqref{a2})$. $E_2$, $E_7$ and $E_{11}$ are defined by
\eqref{e2}, \eqref{3d-dt-E-8} and \eqref{3d-dt-E-13} respectively.
\end{proposition}

\pf Proposition \ref{prop 3.1} can be derived from Lemmas \ref{3d-le:3.1-1}-\ref{3d-le:rho} below.

\begin{lemma}\label{3d-le:3.1-1}
Under the conditions of Proposition \ref{prop 3.1}, it holds that
\be \displaystyle\sup\limits_{0\le t\le T}\int_{\mathbb{R}^3} P \le (\gamma-1)E_0,\label{3d-basic-1}\\[2mm]
\displaystyle\int_0^{T}\int_{\mathbb{R}^3} |\nabla u|^2\le \frac{E_0}{\mu}.\label{3d-basic-2}
\ee
\end{lemma}

\pf Multiplying $(\ref{full N-S+1})_1$ by $G'(\rho)$ and $(\ref{full N-S+1})_2$ by $u$ and integrating over $\mathbb{R}^3\times [0,T]$,
then using $(\ref{3d-boundary})$, one gets
\be\label{3d-basic-pr-1} \displaystyle\sup\limits_{0\le t\le T}\int_{\mathbb{R}^3}
\left(\frac{1}{2}\rho|u|^2+G(\rho)\right)
+\int_0^T\int_{\mathbb{R}^3}\left(\mu|\nabla u|^2+(\lambda+\mu)|\mathrm{div}u|^2\right)\le E_0,\ee
which gives (\ref {3d-basic-1}) and (\ref{3d-basic-2}).
\endpf

The following is the key {\it a priori} estimate which is essential to close the {\it a priori} assumptions \eqref{a priori assumption 1}.
\begin{lemma}\label{3d-le:hp}
Under the conditions of Proposition \ref{prop 3.1}, assume further
that $ \displaystyle 1<\gamma\leq\frac{3}{2}$, we have
\be\label{3d-basic-3}\begin{aligned}
 \displaystyle\int_0^{\sigma(T)}\int_{\mathbb{R}^3}|\nabla u|^2\le \frac{(\gamma-1)^{\frac{2}{3}}E_0^{\frac{2}{3}}}{\mu}E_1,
\end{aligned}\ee
where $\displaystyle E_1=C(\overline{\rho},M)\left(1+\frac{(\gamma-1)^{\frac{2}{9}}}{\mu^{\frac{2}{3}}}\right)$.
\end{lemma}
\pf First, assume that $\displaystyle\frac{(\gamma-1)^\frac{1}{3}E_0}{\mu}\le 1$.
 Multiplying $(\ref{full N-S+1})_2$ by $u$ and then integrating the resulting equality over
$\mathbb{R}^3$, and using integration
by parts, we have
\be\label{3d-dt rho u}\begin{aligned}
\frac{1}{2}\frac{d}{dt}\int_{\mathbb{R}^3}\rho|u|^2+\int_{\mathbb{R}^3}\left(\mu|\nabla u|^2+(\lambda+\mu)|\mathrm{div}u|^2\right)
= \int_{\mathbb{R}^3}P \mathrm{div}u.
\end{aligned}
\ee Integrating (\ref{3d-dt rho u}) over $[0,\sigma(T)]$, using (\ref{3d-basic-1}) and Cauchy inequality, we have
\be\label{3d-basic-pr-2}&&\displaystyle\sup\limits_{0\le t\le \sigma(T)}\frac{1}{2}\int_{\mathbb{R}^3}
\displaystyle\rho|u|^2+\int_0^{\sigma(T)}\int_{\mathbb{R}^3}\left(\frac{\mu}{2}|\nabla u|^2+(\lambda+\mu)|\mathrm{div}u|^2\right)\nonumber\\[2mm]
&\le &\displaystyle\frac{1}{2}\int_{\mathbb{R}^3}\rho_0|u_0|^2+\frac{C}{\mu}\int_0^{\sigma(T)}\int_{\mathbb{R}^3}|P|^2 \nonumber\\[2mm]
&\le &\displaystyle\frac{1}{2}\|\rho_0\|_{L^{\frac{3}{2}}}\|u_0\|^2_{L^6}+\frac{C(\overline{\rho}) (\gamma-1)E_0}{\mu}\nonumber\\[2mm]
&\le &\displaystyle C(\overline{\rho},M)(\gamma-1)^{\frac{2}{3}}E_0^{\frac{2}{3}}\left(1+\frac{(\gamma-1)^{\frac{1}{3}}E_0^{\frac{1}{3}}}{\mu}\right)\nonumber\\[2mm]
&\le& \displaystyle
C(\overline{\rho},M)(\gamma-1)^{\frac{2}{3}}E_0^{\frac{2}{3}}\left(1+\frac{(\gamma-1)^{\frac{2}{9}}}{\mu^{\frac{2}{3}}}\right),
\ee where $\displaystyle 1<\gamma\le\frac{3}{2}$ has been used.
This completes the proof of Lemma \ref{3d-le:hp}.

\endpf
\begin{lemma}\label{3d-le:3.2} Under the conditions of Proposition
\ref{prop 3.1}, it holds that
\be \displaystyle
A_1(T)&\leq &\displaystyle\frac{C(2\mu+\lambda)}{\mu}\int_0^T \sigma
\|\nabla
u\|_{L^3}^3+\frac{C\gamma}{\mu}\int_0^{T}\int_{\mathbb{R}^3}
 \sigma P |\nabla u|^2\nonumber\\
&&\displaystyle+\frac{C
(2\mu+\lambda)}{\mu}\int_0^{\sigma(T)}\|\nabla
u\|_{L^2}^2+\frac{C(\gamma-1)E_0}{\mu^2},\label{3d-A1}\\[2mm]
\displaystyle A_2(T)&\leq& C
A_1(T)+\frac{C\gamma^2}{\mu^2}\int_0^T\int_{\mathbb{R}^3}
\sigma^2|P\nabla
u|^2+C\left(\frac{2\mu+\lambda}{\mu}\right)^2\int_0^T\int_{\mathbb{R}^3}
\sigma^2 |\nabla u|^4. \label{3d-A2}\ee
\end{lemma}
\pf The proof of (\ref{3d-A1}) and (\ref{3d-A2}) is duo to Hoff
\cite{Hoff-jde} and Huang-Li-Xin \cite{Huang-Li-Xin}. For $m\geq0$,
multiplying (\ref{full N-S+1})$_2$ by $\sigma^m \dot{u}$,
integrating the resulting equality over $\mathbb{R}^3$, we have
\be\label{3d-A1-eq}
\displaystyle\int_{\mathbb{R}^3}\sigma^m\rho|\dot{u}|^2&\displaystyle=&\int_{\mathbb{R}^3}\left(-\sigma^m \dot{u}\cdot\nabla P+\mu \sigma^m \Delta u\cdot \dot{u}+(\lambda+\mu)\sigma^m\nabla \mathrm{div}u\cdot \dot{u} \right)\nonumber\\
&\displaystyle=&\sum\limits_{i=1}^3I_i. \ee Integrating by parts
gives \be\label{3d-A1-I-1}
I_1&=&\displaystyle-\int_{\mathbb{R}^3}\sigma^m \dot{u}\cdot\nabla P\nonumber\\[2mm]
&=&\displaystyle\int_{\mathbb{R}^3}\sigma^m  \mathrm{div}u_tP+\int_{\mathbb{R}^3}\sigma^m  \mathrm{div}(u\cdot\nabla u)P\nonumber\\[2mm]
&=&\displaystyle\frac{d}{dt}\int_{\mathbb{R}^3}\sigma^m\mathrm{div}uP-m\int_{\mathbb{R}^3}\sigma^{m-1}\sigma'\mathrm{div}uP\nonumber\\[2mm]
&&-\int_{\mathbb{R}^3}\sigma^m  \mathrm{div}u P'\rho_t+\int_{\mathbb{R}^3}\sigma^m  \mathrm{div}(u\cdot\nabla u)P\nonumber\\[2mm]
&\le&\displaystyle\left(\int_{\mathbb{R}^3} \sigma^m \mathrm{div}u P \right)_t-m\sigma^{m-1}\sigma'\int_{\mathbb{R}^3}\mathrm{div}u P\nonumber\\[2mm]
&&\displaystyle+(\gamma-1)\sigma^m \int_{\mathbb{R}^3}
P|\mathrm{div}u|^2+C\int_{\mathbb{R}^3} \sigma^m P|\nabla u|^2,\\[2mm]
I_2&=&\displaystyle\int_{\mathbb{R}^3}\mu \sigma^m \Delta u\cdot \dot{u}\nonumber\\[2mm]
&=&\displaystyle-\int_{\mathbb{R}^3}\mu \sigma^m \nabla u\cdot\nabla u_t+\int_{\mathbb{R}^3}\mu \sigma^m\Delta u(u\cdot\nabla u)\nonumber\\[2mm]
&\le&\displaystyle-\frac{\mu}{2}\left(\int_{\mathbb{R}^3} \sigma^m
|\nabla u|^2\right)_t+\frac{\mu
m}{2}\sigma^{m-1}\sigma'\int_{\mathbb{R}^3}|\nabla u|^2
 +C\mu \sigma^m \int_{\mathbb{R}^3} |\nabla u|^3\label{3d-A1-I-2}
\ee and \be\label{3d-A1-I-3}
I_3&=&\displaystyle\int_{\mathbb{R}^3}(\lambda+\mu)\sigma^m\nabla \mathrm{div}u\cdot \dot{u}\nonumber\\[2mm]
&\le&\displaystyle-\frac{\mu+\lambda}{2}\left(\sigma^m \int_{\mathbb{R}^3} |\mathrm{div}u|^2\right)_t+\frac{m(\mu+\lambda)}{2}\sigma^{m-1}\sigma'\int_{\mathbb{R}^3}|\mathrm{div}u|^2\nonumber\\[2mm]
&&\displaystyle+C(\mu+\lambda)\sigma^m \int_{\mathbb{R}^3} |\nabla
u|^3. \ee Substituting (\ref{3d-A1-I-1})-(\ref{3d-A1-I-3}) into
(\ref{3d-A1-eq}) shows that \be\label{3d-A1-I-4}
&\displaystyle&\frac{d}{dt}\left(\frac{\mu}{2}\sigma^m\|\nabla u\|_{L^2}^2+\frac{(\lambda+\mu)}{2}\sigma^m\|\mathrm{div}u\|_{L^2}^2-\sigma^m\int_{\mathbb{R}^3}\mathrm{div}u P\right)+\int_{\mathbb{R}^3}\sigma^m \rho|\dot{u}|^2\nonumber\\[2mm]
&\leq&
\displaystyle-m\sigma^{m-1}\sigma'\int_{\mathbb{R}^3}\mathrm{div}u
P+C\gamma\int_{\mathbb{R}^3} \sigma^m P|\nabla u|^2
+C(2\mu+\lambda)\sigma^m \int_{\mathbb{R}^3} |\nabla u|^3\nonumber\\[2mm]
&&\displaystyle+\frac{m(4\mu+3\lambda)}{2}\sigma^{m-1}\sigma'\int_{\mathbb{R}^3}|\nabla
u|^2. \ee
 Applying
\eqref{3d-basic-1}, integrating (\ref{3d-A1-I-4}) over $(0,T)$, we get \be\label{3d-A1-I-5}
&&\displaystyle\sup_{0\leq t\leq T}\left(\frac{\mu}{4}\sigma^m\|\nabla u\|_{L^2}^2+\frac{(\lambda+\mu)}{2}\sigma^m\|\mathrm{div}u\|_{L^2}^2\right)+\int_0^{T}\int_{\mathbb{R}^3}\sigma^m \rho|\dot{u}|^2\nonumber\\[2mm]
&\leq&\displaystyle \frac{C(\gamma-1)E_0}{\mu}+C\int_0^{\sigma(T)}\int_{\mathbb{R}^3}|P||\mathrm{div}u| +\frac{m(4\mu+3\lambda)}{2}\int_0^{\sigma(T)}\int_{\mathbb{R}^3}|\nabla u|^2\nonumber\\[2mm]
&&\displaystyle+C(2\mu+\lambda)\int_0^{T}\sigma^m \int_{\mathbb{R}^3} |\nabla u|^3+(3\gamma-2)\int_0^{T}\int_{\mathbb{R}^3} \sigma^m P|\nabla u|^2\nonumber\\[2mm]
&\leq&\displaystyle
\frac{C(\gamma-1)E_0}{\mu}+C(2\mu+\lambda)\int_0^{\sigma(T)}\|\nabla
u\|_{L^2}^2
+C(2\mu+\lambda)\int_0^{T}\sigma^m \int_{\mathbb{R}^3} |\nabla u|^3\nonumber\\[2mm]
&&\displaystyle+C\gamma\int_0^{T}\int_{\mathbb{R}^3} \sigma^m
P|\nabla u|^2. \ee Choosing $m=1$, one gets (\ref{3d-A1}).

\vspace{2mm}

Next, for $m\geq0$, operating $\sigma^m\dot{u}^j[\partial/\partial
t+\mathrm{div}(u\cdot)]$ on $(\ref{full N-S+1})_2^j$, summing over
$j$, and integrating the resulting equation over $\mathbb{R}^3\times
[0,T]$, we obtain after integration by parts \be\label{3d-A2-eq}
&&\displaystyle\sup_{0\leq t\leq T}\left(\frac{1}{2}\sigma^m\int_{\mathbb{R}^3}\rho|\dot{u}|^2\right)-\frac{m}{2}\int_0^T\sigma^{m-1}\sigma'\int_{\mathbb{R}^3}\rho|\dot{u}|^2\nonumber\\[2mm]
&=&\displaystyle-\int_0^T\int_{\mathbb{R}^3}\sigma^m\dot{u}^j\left[\partial_jP_t+\mathrm{div}(\partial_jP u)\right]+\mu\int_0^T\int_{\mathbb{R}^3}\sigma^m\dot{u}^j\left[\Delta u_t^j+\mathrm{div}(u \Delta u^j)\right]\nonumber\\[2mm]
&&\displaystyle+(\lambda+\mu)\int_0^T\int_{\mathbb{R}^3}\sigma^m\dot{u}^j\left[\partial_t\partial_j\mathrm{div}u+\mathrm{div}(u\partial_j\mathrm{div}u )\right]\nonumber\\[2mm]
&=&\displaystyle\sum\limits_{i=1}^3II_i. \ee Integrating by parts
leads to
\be
II_1&=&\displaystyle-\int_0^T\int_{\mathbb{R}^3}\sigma^m\dot{u}^j\left[\partial_jP_t+\mathrm{div}(\partial_jP u)\right]\nonumber\\[2mm]
&=&\displaystyle\int_0^T\int_{\mathbb{R}^3}\sigma^m\partial_{j}\dot{u}^jP_t+\int_0^T\int_{\mathbb{R}^3} \sigma^m\partial_{k}\dot{u}^j(\partial_{j}Pu^k)\nonumber\\[2mm]
&=&\displaystyle-\int_0^T\int_{\mathbb{R}^3}\sigma^m\mathrm{div}\dot{u}P'\left(\rho\mathrm{div}u+u\cdot \nabla \rho\right)+\int_0^T\int_{\mathbb{R}^3}\sigma^m\mathrm{div}\dot{u}\mathrm{div}uP\nonumber\\[2mm]
&&\displaystyle+\int_0^T\int_{\mathbb{R}^3}\sigma^m\mathrm{div}\dot{u}u\cdot\nabla
P
-\int_0^T\int_{\mathbb{R}^3}\sigma^m \partial_k\dot{u}^j \partial_j u^k P\nonumber\\[2mm]
&\leq&\displaystyle C\gamma\int_0^T\int_{\mathbb{R}^3} \sigma^m P|\nabla \dot{u}||\nabla u|\nonumber\\[2mm]
&\le&\displaystyle \frac{\mu}{4}\int_0^T\int_{\mathbb{R}^3} \sigma^m
|\nabla \dot{u}|^2+\frac{C\gamma^2}{\mu}\int_0^T\int_{\mathbb{R}^3}
\sigma^m |P\nabla u|^2,\label{3d-A2-I-1-1}\\[2mm]
II_2&=&\displaystyle\mu\int_0^T\int_{\mathbb{R}^3}\sigma^m\dot{u}^j\left[\Delta u_t^j+\mathrm{div}(u \Delta u^j)\right]\nonumber\\[2mm]
&=&\displaystyle-\mu\int_0^T\int_{\mathbb{R}^3}\sigma^m|\partial_{i}\dot{u}^j|^2
+\mu\int_0^T\int_{\mathbb{R}^3}\sigma^m\partial_{i}\dot{u}^j\partial_{i}(u^k\partial_ku^j)\nonumber\\[2mm]
&&\displaystyle-\mu\int_0^T\int_{\mathbb{R}^3}\sigma^m\partial_{k}\dot{u}^j\partial_{i}(u^k\partial_iu^j)
-\mu\int_0^T\int_{\mathbb{R}^3}\sigma^m\partial_{k}\dot{u}^j\partial_{i}u^k\partial_iu^j\nonumber\\[2mm]
&\le&\displaystyle-\mu\int_0^T\int_{\mathbb{R}^3}\sigma^m|\nabla\dot{u}|^2+C\mu\int_0^T\int_{\mathbb{R}^3}\sigma^m|\nabla\dot{u}||\mathrm{div}u|^2
+C\mu\int_0^T\int_{\mathbb{R}^3}\sigma^m|\nabla\dot{u}||\nabla u|^2\nonumber\\[2mm]
&\le&\displaystyle-\mu\int_0^T\int_{\mathbb{R}^3} \sigma^m |\nabla \dot{u}|^2+C\mu\int_0^T\int_{\mathbb{R}^3} \sigma^m |\nabla \dot{u}||\nabla u|^2\nonumber\\[2mm]
&\le&\displaystyle -\frac{3\mu}{4}\int_0^T\int_{\mathbb{R}^3}
\sigma^m |\nabla \dot{u}|^2+C\mu\int_0^T\int_{\mathbb{R}^3} \sigma^m
|\nabla u|^4.\label{3d-A2-I-2}
 \ee
Similarly, \be\label{3d-A2-I-3}
II_3&\le&\displaystyle -\frac{\mu+\lambda}{2}\int_0^T\int_{\mathbb{R}^3} \sigma^m |\mathrm{div}\dot{u}|^2\nonumber\\[2mm]
&&\displaystyle+C(\mu+\lambda)\left(1+\frac{\lambda}{\mu}\right)\int_0^T\int_{\mathbb{R}^3}
\sigma^m |\nabla u|^4+ \frac{\mu}{4}\int_0^T\int_{\mathbb{R}^3}
\sigma^m |\nabla \dot{u}|^2. \ee Substituting
(\ref{3d-A2-I-1-1})-(\ref{3d-A2-I-3}) into (\ref{3d-A2-eq}) shows
that \be\label{3d-A2-I-4}
&&\displaystyle\sup_{0\leq t\leq T}\sigma^m\int_{\mathbb{R}^3}\rho|\dot{u}|^2+\mu\int_0^T\int_{\mathbb{R}^3} \sigma^m |\nabla \dot{u}|^2+(\mu+\lambda)\int_0^T\int_{\mathbb{R}^3} \sigma^m |\mathrm{div}\dot{u}|^2\nonumber\\[2mm]
&\le &\displaystyle C\int_0^T\sigma^{m-1}\sigma'\int_{\mathbb{R}^3}\rho|\dot{u}|^2+\frac{C\gamma^2}{\mu}\int_0^T\int_{\mathbb{R}^3} \sigma^m |P\nabla u|^2\nonumber\\[2mm]
&&\displaystyle+\frac{C(2\mu+\lambda)^2}{\mu}\int_0^T\int_{\mathbb{R}^3}
\sigma^m |\nabla u|^4,
\ee
where (\ref{viscosity assumption}) has
been used. Taking $m=2$, we immediately obtain (\ref{3d-A2}). The
proof of Lemma \ref{3d-le:3.2} is completed.
\endpf

\begin{lemma}\label{3d-le:3.3}Under the conditions of Proposition \ref{prop 3.1}, it holds that
\be\sup\limits_{0\le t\le
\sigma(T)}\int_{\mathbb{R}^3}|\nabla
u|^2+\int_0^{\sigma(T)}\int_{\mathbb{R}^3}\frac{\rho|\u|^2}{\mu}\le
E_2\label{3d-A-1},\\[2mm]
\sup\limits_{0\le t\le
\sigma(T)}t\int_{\mathbb{R}^3}\frac{\rho
|\u|^2}{\mu}+\int_0^{\sigma(T)}\int_{\mathbb{R}^3}t|\nabla\u|^2\le
E_3,\label{3d-A-2}\ee provided
$\displaystyle\frac{(\gamma-1)^{\frac{1}{6}}E_0^\frac{1}{2}}{\mu^\frac{1}{2}}
\leq \min\left\{\Big(4C(\bar{\rho})\Big)^{-6},1 \right\}\triangleq
\varepsilon_1$.

\end{lemma}

\pf  First, we assume that
$\displaystyle\frac{(\gamma-1)^\frac{1}{3}E_0}{\mu}\le
1$. Multiplying (\ref{full N-S+1})$_2$ by $u_t$, integrating the
resulting equality over $\mathbb{R}^3$ and using \eqref{hp4}, we
have \be \label{3d-le:3.3-eq}
&&\displaystyle\frac{d}{dt}\left(\frac{\mu}{2}\|\nabla u\|_{L^2}^2+\frac{(\lambda+\mu)}{2}\|\mathrm{div}u\|_{L^2}^2-\int_{\mathbb{R}^3}\mathrm{div}u P\right)+\int_{\mathbb{R}^3}\rho|\dot{u}|^2\nonumber\\[2mm]
&=&\displaystyle\int_{\mathbb{R}^3}\rho\dot{u}(u\cdot\nabla u)-\int_{\mathbb{R}^3}\mathrm{div}u P_t\nonumber\\[2mm]
&\le&\displaystyle C(\overline{\rho})\left(\int_{\mathbb{R}^3}\rho|\dot{u}|^2\right)^{\frac{1}{2}}\left(\int_{\mathbb{R}^3}\rho|u|^3\right)^{\frac{1}{3}}\|\nabla u\|_{L^6}+\int_{\mathbb{R}^3}\mathrm{div}u\mathrm{div}(Pu)+(\gamma-1)\int_{\mathbb{R}^3}P|\mathrm{div}u|^2\nonumber\\[2mm]
&\le&\displaystyle{\frac{C(\overline{\rho})}{2\mu+\lambda}}\left(\int_{\mathbb{R}^3}\rho|\dot{u}|^2\right)^{\frac{1}{2}}\left(\int_{\mathbb{R}^3}\rho|u|^3\right)^{\frac{1}{3}}
\Big((CN_6+1)\|\rho\u\|_{L^2}+\|P\|_{L^6}\Big)\nonumber\\[2mm]
&&\displaystyle-\int_{\mathbb{R}^3}Pu\cdot\nabla\mathrm{div}u
+C(\overline{\rho})(\gamma-1)\int_{\mathbb{R}^3}|\nabla u|^2\nonumber\\[2mm]
&\le&\displaystyle\frac{C(\overline{\rho})\mu(CN_6+1)}{2\mu+\lambda}\left(\int_{\mathbb{R}^3}\rho|\dot{u}|^2\right)A^{\frac{1}{3}}_3(\sigma(T))
+\frac{C(\overline{\rho})\mu }{2\mu+\lambda}\left(\int_{\mathbb{R}^3}\rho|\dot{u}|^2\right)^{\frac{1}{2}}A^{\frac{1}{3}}_3(\sigma(T))\|P\|_{L^6}\nonumber\\[2mm]
&&\displaystyle-\frac{1}{2\mu+\lambda}\int_{\mathbb{R}^3}Pu\cdot\nabla G+\frac{1}{2(2\mu+\lambda)}\int_{\mathbb{R}^3}\mathrm{div}uP^2+C(\bar{\rho})(\gamma-1)\|\nabla u \|_{L^2}^2\nonumber\\[2mm]
&\le&\displaystyle{C(\overline{\rho})}\left(\int_{\mathbb{R}^3}\rho|\dot{u}|^2\right)A^{\frac{1}{3}}_3(\sigma(T))
+{C(\bar{\rho})}A_3^{\frac{1}{3}}(\sigma(T))\|P\|_{L^6}^2
+C(\bar{\rho})(\gamma-1)\|\nabla u \|_{L^2}^2\nonumber\\[2mm]
&&\displaystyle+\frac{C}{2\mu+\lambda}\|P\|_{L^3}\|\nabla
u\|_{L^2}\|\rho
\dot{u}\|_{L^2}+\frac{C}{2\mu+\lambda}\left(\|\nabla u\|_{L^2}^2+\|P\|_{L^4}^4\right)\nonumber\\[2mm]
&\le&\displaystyle{C(\overline{\rho})}\left(\int_{\mathbb{R}^3}\rho|\dot{u}|^2\right)A^{\frac{1}{3}}_3(\sigma(T))
+{C(\bar{\rho})}A_3^{\frac{1}{3}}(\sigma(T))(\gamma-1)^\frac{1}{3}E_0^\frac{1}{3}
+C(\bar{\rho})(\gamma-1)\|\nabla u \|_{L^2}^2\nonumber\\[2mm]
&&\displaystyle+\frac{1}{4}\|\sqrt{\rho}
\dot{u}\|_{L^2}^2+\frac{C(\overline{\rho})}{(2\mu+\lambda)^2}(\gamma-1)^{\frac{2}{3}}E^\frac{2}{3}_0\|\nabla
u\|_{L^2}^2
+\frac{C}{2\mu+\lambda}\|\nabla u\|_{L^2}^2+\frac{C(\overline{\rho})}{2\mu+\lambda}(\gamma-1)E_0\nonumber\\[2mm]
&\le&\displaystyle{C(\overline{\rho})}\left(\int_{\mathbb{R}^3}\rho|\dot{u}|^2\right)A^{\frac{1}{3}}_3(\sigma(T))
+{C(\bar{\rho})}A_3^{\frac{1}{3}}(\sigma(T))(\gamma-1)^\frac{1}{3}E_0^\frac{1}{3}
+C(\bar{\rho})(\gamma-1)\|\nabla u \|_{L^2}^2\nonumber\\[2mm]
&&\displaystyle+\frac{1}{4}\|\sqrt{\rho}
\dot{u}\|_{L^2}^2+\frac{C(\overline{\rho})}{(2\mu+\lambda)^2}(\gamma-1)^{\frac{2}{3}}E^\frac{2}{3}_0\|\nabla
u\|_{L^2}^2
+\frac{C}{2\mu+\lambda}\|\nabla u\|_{L^2}^2\nonumber\\[2mm]
&&+\frac{C(\overline{\rho})}{2\mu+\lambda}(\gamma-1)E_0, \ee where
we have used \eqref{hp3}. Integrating (\ref{3d-le:3.3-eq}) over $(0,
\sigma(T))$, we obtain that \be \label{3d-le:3.3-eq-1}
&&\displaystyle \sup_{0 \leq t
\leq \sigma(T)}\left\{\frac{\mu}{2}\|\nabla u\|_{L^2}^2+\frac{(\lambda+\mu)}{2}\|\mathrm{div}u\|_{L^2}^2-\int_{\mathbb{R}^3}\mathrm{div}u P\right\}+\frac{1}{2}\int_0^{\sigma(T)}\int_{\mathbb{R}^3}\rho|\dot{u}|^2\nonumber\\[2mm]
&\le&\displaystyle{C(\bar{\rho})}A_3^{\frac{1}{3}}(\sigma(T))(\gamma-1)^\frac{1}{3}E_0^\frac{1}{3}
+C(\bar{\rho})(\gamma-1)\int_0^{\sigma(T)}\|\nabla u \|_{L^2}^2\nonumber\\[2mm]
&&\displaystyle+\frac{C(\overline{\rho})}{(2\mu+\lambda)^2}(\gamma-1)^{\frac{2}{3}}E^\frac{2}{3}_0\int_0^{\sigma(T)}\|\nabla
u\|_{L^2}^2
+\frac{C}{2\mu+\lambda}\int_0^{\sigma(T)}\|\nabla u\|_{L^2}^2\nonumber\\[2mm]
&&+\frac{C(\overline{\rho})}{2\mu+\lambda}(\gamma-1)E_0+C\mu M, \ee
provided
$\displaystyle\frac{(\gamma-1)^{\frac{1}{6}}E_0^\frac{1}{2}}{\mu^\frac{1}{2}}
\leq \Big(4C(\bar{\rho})\Big)^{-6} $.

\noindent\underline{\textbf{Case 1:}} \ $\displaystyle1 < \gamma \leq \frac{3}{2}$.

Using\eqref{3d-basic-3}, we get \be\label{1}\begin{aligned} &\sup_{0 \leq t
\leq \sigma(T)}\left\{\|\nabla
u\|_{L^2}^2+(\lambda+\mu)\|\mathrm{div}u\|_{L^2}^2\right\}
+\frac{1}{\mu}\int_0^{\sigma(T)}\int_{\mathbb{R}^3}\rho|\dot{u}|^2
\le E^1_2,
\end{aligned}
\ee where
\begin{equation*}
\displaystyle
E^1_{2}=\frac{C(\gamma-1)^\frac{2}{9}}{\mu^\frac{2}{3}}+\frac{C(\gamma-1)^\frac{13}{9}E_1}{\mu^\frac{4}{3}}
+\frac{C(\gamma-1)^\frac{8}{9}E_1}{\mu^\frac{8}{3}}
+\frac{C(\gamma-1)^\frac{4}{9}E_1}{\mu^\frac{7}{3}}
+\frac{C(\gamma-1)^\frac{2}{3}}{\mu} +CM,
\end{equation*}
and  we have also used the facts that
$\displaystyle\frac{(\gamma-1)^\frac{1}{3}E_0}{\mu}\le
1$ and $\mu+\lambda>0$.

\noindent\underline{\textbf{Case 2:}} \ $\displaystyle \gamma > \frac{3}{2}$.

Using\eqref{3d-basic-2}, we have \be\label{2}\begin{aligned} &\sup_{0 \leq
t \leq \sigma(T)}\left\{\|\nabla
u\|_{L^2}^2+(\lambda+\mu)\|\mathrm{div}u\|_{L^2}^2\right\}
+\frac{1}{\mu}\int_0^{\sigma(T)}\int_{\mathbb{R}^3}\rho|\dot{u}|^2
\le E^2_2,
\end{aligned}
\ee where
\begin{equation*}
\displaystyle
E^2_{2}=\frac{C(\gamma-1)^\frac{2}{9}}{\mu^\frac{2}{3}}+\frac{C(\gamma-1)^\frac{2}{3}}{\mu}
+\frac{C(\gamma-1)^\frac{1}{9}}{\mu^\frac{7}{3}} +\frac{C}{\mu^2}
+\frac{C(\gamma-1)^\frac{2}{3}}{\mu} +CM.
\end{equation*}
Combining \eqref{1}-\eqref{2}, the result leads to (\ref{3d-A-1}),
 where
\be\label{e2} E_2=\max\Big\{E^1_{2},\ E^2_{2}\Big\}. \ee

Taking $m=1$ and $T=\sigma(T)$ in (\ref{3d-A2-I-4}), we obtain that
\be\label{3d-A2-II-4}
&&\displaystyle\sigma\int_{\mathbb{R}^3}\rho|\dot{u}|^2+\mu\int_0^{\sigma(T)}\int_{\mathbb{R}^3}
\sigma|\nabla \dot{u}|^2
+(\mu+\lambda)\int_0^{\sigma(T)}\int_{\mathbb{R}^3} \sigma|\mathrm{div}\dot{u}|^2\nonumber\\[2mm]
&\leq&\displaystyle \mu
E_2+\frac{C\gamma^2}{\mu}\int_0^{\sigma(T)}\int_{\mathbb{R}^3}\sigma|P\nabla
u|^2
+\frac{C(2\mu+\lambda)^2}{\mu}\int_0^{\sigma(T)}\int_{\mathbb{R}^3}\sigma|\nabla u|^4\nonumber\\[2mm]
&\leq &\displaystyle\mu
E_2+\frac{C\gamma^2}{\mu}\left(\int_0^{\sigma(T)}\int_{\mathbb{R}^3}P^4\right)^{\frac{1}{2}}
\left(\int_0^{\sigma(T)}\int_{\mathbb{R}^3}\sigma^2|\nabla u|^4\right)^{\frac{1}{2}}\nonumber\\[2mm]
&&\displaystyle+\frac{C(2\mu+\lambda)^2}{\mu}\int_0^{\sigma(T)}\int_{\mathbb{R}^3}\sigma|\nabla u|^4\nonumber\\[2mm]
&\leq &\displaystyle\mu
E_2+\frac{C\gamma^2(\gamma-1)^\frac{1}{3}E_2^\frac{3}{4}}{\mu^\frac{5}{4}}+\frac{C\gamma^2(\gamma-1)^\frac{1}{2}E_2^\frac{1}{4}}{\mu^\frac{7}{4}}
+C\left(2+\frac{\lambda}{\mu}\right)^2\left(\frac{E_2^\frac{3}{2}}{\mu^\frac{1}{2}}
+\frac{(\gamma-1)^\frac{1}{3}E_2^\frac{1}{2}}{\mu^\frac{3}{2}}\right)\nonumber\\[2mm]
&\leq&\displaystyle\mu E_3, \ee where
$$ E_3=E_2+\frac{C\gamma^2(\gamma-1)^\frac{1}{3}E_2^\frac{3}{4}}{\mu^\frac{9}{4}}+\frac{C\gamma^2(\gamma-1)^\frac{1}{2}E_2^\frac{1}{4}}{\mu^\frac{11}{4}}
+C\left(2+\frac{\lambda}{\mu}\right)^2\left(\frac{E_2^\frac{3}{2}}{\mu^\frac{3}{2}}
+\frac{(\gamma-1)^\frac{1}{3}E_2^\frac{1}{2}}{\mu^\frac{5}{2}}\right).
$$
To get (\ref{3d-A2-II-4}), we have used the following estimate \be
&&\displaystyle\int_0^{\sigma(T)}\int_{\mathbb{R}^3}\sigma|\nabla u|^4\nonumber\\[2mm]
&\le&\displaystyle\sup_{0\leq t \leq \sigma(T)}\|\nabla
u\|_{L^2}\int_0^{\sigma(T)}\sigma\|\nabla
u\|^3_{L^6}\nonumber\\[2mm]
&\le&\displaystyle
E_2^\frac{1}{2}\int_0^{\sigma(T)}\sigma\left(\frac{(4\mu+\lambda)^3}{(2\mu+\lambda)^{3}\mu^3}\|\rho
\dot{u}\|_{L^2}^3+\frac{\|P\|^3_{L^6}}{(2\mu+\lambda)^{3}}\right)\nonumber\\[2mm]
&\le&\displaystyle\frac{CE_2^\frac{1}{2}}{\mu^3}\sup_{0\leq t \leq
\sigma(T)}\sigma\|\rho \dot{u}\|_{L^2}\int_0^{\sigma(T)}\|\rho
\dot{u}\|_{L^2}^2+\frac{CE_2^\frac{1}{2}(\gamma-1)^\frac{1}{2}E_0^\frac{1}{2}}{\mu^3}\nonumber\\[2mm]
&\le&\displaystyle\frac{CE_2^\frac{1}{2}}{\mu^3}\left(\mu
A_2(T)\right)^\frac{1}{2}\mu
E_2+\frac{CE_2^\frac{1}{2}(\gamma-1)^\frac{1}{3}}{\mu^\frac{5}{2}}\nonumber\\[2mm]
&\le&\displaystyle
\frac{CE_2^\frac{3}{2}}{\mu^\frac{3}{2}}+\frac{CE_2^\frac{1}{2}(\gamma-1)^\frac{1}{3}}{\mu^\frac{5}{2}},
\ee due to H\"older inequality, \eqref{hp3}, \eqref{a priori
assumption 1}, the facts that
$\displaystyle\frac{(\gamma-1)^\frac{1}{3}E_0}{\mu}\le
1$ and $\mu+\lambda>0$. The proof of Lemma \ref{3d-le:3.3} is
completed.
\endpf

\begin{lemma}\label{bound}
Under the conditions of Proposition \ref{prop 3.1}, we have
\be\sup\limits_{0\le t\le T}\int_{\mathbb{R}^3}|\nabla
u|^2+\int_0^{T}\int_{\mathbb{R}^3}\frac{\rho|\u|^2}{\mu}\le
C(E_2+1),\\[2mm]
\sup\limits_{0\le t\le
T}\sigma\int_{\mathbb{R}^3}\frac{\rho
|\u|^2}{\mu}+\int_0^{T}\int_{\mathbb{R}^3}\sigma|\nabla\u|^2\le
C(E_3+1),\ee provided
$\displaystyle\frac{(\gamma-1)^{\frac{1}{6}}E^\frac{1}{2}_0}{\mu^\frac{1}{2}}
\leq \varepsilon_1$.

\end{lemma}
\pf By \eqref{a priori assumption 1} and Lemma \ref{3d-le:3.3}, we
immediately get Lemma \ref{bound}.
\endpf\\

Next, we will close the $a \ priori$ assumption on $A_3(T)$.
\begin{lemma}\label{3d-A3}
Under the conditions of Proposition \ref{prop 3.1}, it holds that
\be A_3(\sigma(T))\leq
\left\{\frac{(\gamma-1)^{\frac{1}{6}}E_0^\frac{1}{2}}{\mu^\frac{1}{2}}\right\}^\frac{1}{2},
\ee provided \be
\frac{(\gamma-1)^{\frac{1}{6}}E_0^\frac{1}{2}}{\mu^\frac{1}{2}} &\le
&\min \left\{
C(\bar{\rho})^{-2}(\gamma-1)^{-\frac{2}{3}}E_2^{-3}\mu^5\Big|_{(1 <
\gamma\leq \frac{3}{2})},\
C(\bar{\rho})^{-1}\mu^{\frac{9}{4}}E_2^{-\frac{3}{4}}\Big|_{(\gamma >
\frac{3}{2})}, \varepsilon_1
\right\}\nonumber\\[2mm]
&\triangleq & \varepsilon_2. \ee

\end{lemma}
\pf
\underline{\textbf{Case 1:}} \ $\displaystyle1 < \gamma \leq \frac{3}{2}$.  \be
\displaystyle A_3(\sigma(T))=\displaystyle\sup\limits_{0\le t\le \sigma(T)}\displaystyle\int_{\mathbb{R}^3}\frac{\rho|u|^3}{\mu^3}&\le&\displaystyle\frac{1}{\mu^3}\sup\limits_{0\le t\le \sigma(T)}\|\rho\|_{L^2}\|u\|_{L^6}^3\nonumber\\[2mm]
&\leq& \displaystyle\frac{C(\bar{\rho})}{\mu^3}(\gamma-1)^{\frac{1}{2}}E_0^\frac{1}{2}E_2^\frac{3}{2}\nonumber\\[2mm]
&\leq&
\displaystyle\left(\frac{(\gamma-1)^{\frac{1}{6}}E_0^\frac{1}{2}}{\mu^\frac{1}{2}}\right)^{\frac{1}{2}},
\ee provided $\displaystyle
\frac{(\gamma-1)^{\frac{1}{6}}E_0^\frac{1}{2}}{\mu^\frac{1}{2}}\le
C(\bar{\rho})(\gamma-1)^{-\frac{2}{3}}E_2^{-3}\mu^5$, where we have
used  H\"older inequality, Sobolev inequality and Lemma
\ref{3d-le:3.3}.

\vspace{2ex}
\noindent\underline{\textbf{Case 2:}} \ $\displaystyle \gamma> \frac{3}{2}$.

Using  H\"older inequality, (\ref{3d-basic-pr-1}), (\ref{3d-A-1}) and $\displaystyle
2(\gamma-1)\geq 1$, we obtain \be \displaystyle A_3(\sigma(T))
=\displaystyle\sup\limits_{0\le t\le
\sigma(T)}\displaystyle\int_{\mathbb{R}^3}\frac{\rho|u|^3}{\mu^3}
&\le&\displaystyle\frac{C(\bar{\rho})}{\mu^3}\sup\limits_{0\le t\le\sigma(T)}\|\rho^{\frac{1}{2}}u\|_{L^2}^{\frac{3}{2}}\|u\|_{L^6}^{\frac{3}{2}}\nonumber\\[2mm]
&\leq & \displaystyle \frac{C(\bar{\rho})(\gamma-1)^{\frac{1}{4}}E_0^\frac{3}{4}E_2^\frac{3}{4}}{\mu^3}\nonumber\\[2mm]
&\leq&\displaystyle\left(\frac{(\gamma-1)^{\frac{1}{6}}E_0^\frac{1}{2}}{\mu^\frac{1}{2}}\right)^{\frac{1}{2}},
\ee provided
$\displaystyle\frac{(\gamma-1)^{\frac{1}{6}}E_0^\frac{1}{2}}{\mu^\frac{1}{2}}
\le C(\bar{\rho})\mu^{\frac{9}{4}}E_2^{-\frac{3}{4}}$.
\endpf

\begin{lemma}\label{3d-le:3.8}Under the conditions of Proposition \ref{prop 3.1}, it holds that
\be\label{3d-A-1-A-2}A_1(T)+A_2(T)\le
\frac{(\gamma-1)^\frac{1}{6}E_0^{\frac{1}{2}}}{\mu^{\frac{1}{2}}},\ee
provided
$$ \frac{(\gamma-1)^\frac{1}{6}E_0^{\frac{1}{2}}}{\mu^{\frac{1}{2}}}\le
\min\left\{(CE_7)^{-3}\Big|_{(1< \gamma\leq \frac{3}{2})},
(CE_{11})^{-2}\Big|_{( \gamma>\frac{3}{2})},
\varepsilon_2 \right\}\triangleq \varepsilon_3,
$$
where $E_7$ and $E_{11}$ will be determined in \eqref{3d-dt-E-8} and \eqref{3d-dt-E-13}.
\end{lemma}
\pf First we assume that
$\displaystyle\frac{(\gamma-1)^\frac{1}{3}E_0}{\mu}\le
1$. It follows from Lemma \ref{3d-le:3.2} that
\be\label{3d-dt-A-1+A-2}
&&\displaystyle A_1(T)+A_2(T)\nonumber\\[3mm]
&\le&\displaystyle\frac{C(\gamma-1)E_0}{\mu^2}+\frac{C(2\mu+\lambda)^2}{\mu^2}\int_0^T\int_{\mathbb{R}^3}
\sigma^2 |\nabla u|^4
+\frac{C(2\mu+\lambda)}{\mu}\int_0^T \sigma \|\nabla u\|_{L^3}^3\nonumber\\[3mm]
&&\displaystyle+\frac{C\gamma^2}{\mu^2}\int_0^T\int_{\mathbb{R}^3}
\sigma^2 |P\nabla
u|^2+\frac{C\gamma}{\mu}\int_0^{T}\int_{\mathbb{R}^3}
 \sigma P |\nabla u|^2+C\left(\frac{2\mu+\lambda}{\mu}\right)\int_0^{\sigma(T)}\|\nabla u\|_{L^2}^2\nonumber\\[3mm]
&\le&\displaystyle
\frac{C(\gamma-1)E_0}{\mu^2}+\sum\limits_{i=4}^8II_i. \ee Now we give the
estimates of $II_4-II_8$ in the following two cases.
The subsequent estimates proceed with different techniques for
$\displaystyle1<\gamma\le\frac{3}{2}$ and $\displaystyle\gamma>\frac{3}{2}$.

\vspace{2ex}

\noindent\underline{ \textbf{Case 1:}}\ $\displaystyle1 < \gamma \leq \frac{3}{2}$.

For $II_4$, due to \eqref{hp3}, we just estimate $\displaystyle
\int_0^T\int_{\mathbb{R}^3} \sigma^2|\nabla u|^4$ as follows:
\be\label{3d-lemma-u-4} \displaystyle\int_0^T\int_{\mathbb{R}^3}
\sigma^2|\nabla u|^4
&\le& \displaystyle CN_4^4(2\mu+\lambda)^{-3}\int_0^{T}\sigma^2\|\rho\u\|_{L^2}^3\|\nabla u\|_{L^2}\nonumber\\[2mm]
&&\displaystyle+C(2\mu+\lambda)^{-4}\int_0^{T}\sigma^2\|\rho\u\|_{L^2}^3\|P\|_{L^2}\nonumber\\[2mm]
&&\displaystyle+C(2\mu+\lambda)^{-4}\int_0^{T}\sigma^2\|P\|_{L^4}^4\nonumber\\[2mm]
&=&\sum\limits_{i=1}^3III_i. \ee Using H\"older inequality, \eqref{a
priori assumption 1}, (\ref{3d-basic-1}) and Lemma \ref{bound}, we
have \be\label{3d-III 1}
III_1&\le&\displaystyle CN_4^4(2\mu+\lambda)^{-3}\sup\limits_{0\le t\le T}(\sigma^2\|\rho\u\|_{L^2}^2)^{\frac{1}{2}}\int_0^{T}\sigma\|\rho\u\|_{L^2}^2\|\nabla u\|_{L^2}\nonumber\\[2mm]
&\le&\displaystyle CN_4^4(2\mu+\lambda)^{-3}\mu^{\frac{1}{2}}A^{\frac{1}{2}}_2(T)\sup\limits_{0\le t\le T}(\|\nabla u\|_{L^2}^2)^{\frac{1}{2}}\int_0^{T}\sigma\|\rho\|_{L^3}^2\|\nabla\u\|_{L^2}^2\nonumber\\[2mm]
&\le&\displaystyle CN_4^4(2\mu+\lambda)^{-3}\mu^{\frac{1}{2}}A^{\frac{1}{2}}_2(T)(\gamma-1)^\frac{2}{3}E_0^\frac{2}{3}\Big(E_2+1\Big)^\frac{1}{2}\Big(E_3+1\Big)\nonumber\\[2mm]
&\le&\displaystyle
CN_4^4(2\mu+\lambda)^{-3}\mu^{\frac{1}{4}}(\gamma-1)^\frac{3}{4}E_0^{\frac{11}{12}}\Big(E_2+1\Big)^\frac{1}{2}\Big(E_3+1\Big).
\ee Next, it follows from H\"older inequality, \eqref{a priori
assumption 1}, (\ref{3d-basic-1}) and Lemma 3.6 that
\be\label{3d-III_2}
III_2&\le&\displaystyle C(2\mu+\lambda)^{-4}\sup\limits_{0\le t\le T}(\sigma^2\|\rho\u\|_{L^2}^2)^{\frac{1}{2}}\int_0^{T}\sigma\|\rho\u\|_{L^2}^2\|P\|_{L^2}\nonumber\\[2mm]
&\le& \displaystyle C(2\mu+\lambda)^{-4}\mu^{\frac{1}{2}}A^{\frac{1}{2}}_2(T)\sup\limits_{0\le t\le T}\|P\|_{L^2}\int_0^{T}\sigma\|\rho\|_{L^3}^2\|\nabla\u\|_{L^2}^2\nonumber\\[2mm]
&\le&\displaystyle C(2\mu+\lambda)^{-4}\mu^{\frac{1}{4}}(\gamma-1)^{\frac{5}{4}}E_0^{\frac{17}{12}}\Big(E_3+1\Big).
\ee Thus, one gets from \eqref{hp1}, (\ref{3d-III 1}) and
(\ref{3d-III_2}) that
 \be\label{3d-lemma-u-4-e}
\displaystyle\int_0^T\sigma^2\|G\|_{L^4}^4&\le& C\Big((2\mu+\lambda)\|\nabla u\|_{L^2}+\|P\|_{L^2}\Big)\|\rho\dot{u}\|_{L^2}^3\nonumber\\[2mm]
&\le&\displaystyle C(2\mu+\lambda)\mu^{\frac{1}{4}}(\gamma-1)^\frac{3}{4}E_0^{\frac{11}{12}}\Big(E_2+1\Big)^\frac{1}{2}\Big(E_3+1\Big)\nonumber\\[2mm]
&&+C\mu^{\frac{1}{4}}(\gamma-1)^{\frac{5}{4}}E_0^{\frac{17}{12}}\Big(E_3+1\Big)\nonumber\\[2mm]
&\le&\displaystyle
(2\mu+\lambda)\mu^{\frac{1}{4}}(\gamma-1)^{\frac{3}{4}}E_0^{\frac{11}{12}}E_4,
\ee where $\displaystyle
E_4=C\Big(E_2+1\Big)^\frac{1}{2}\Big(E_3+1\Big)+C(\gamma-1)^{\frac{1}{3}}\mu^{-\frac{1}{2}}\Big(E_3+1\Big)$.
Here we have used the facts that
$\displaystyle\frac{(\gamma-1)^\frac{1}{3}E_0}{\mu}\le
1$ and $\mu+\lambda>0$.

To estimate $III_3$, one deduces from $(\ref{full N-S+1})_1$ that
$P$ satisfies \be\label{3d-eq-P}
\begin{aligned}
P_t+u\cdot\nabla P+\gamma P\mathrm{div}u=0.
\end{aligned}
\ee Multiplying (\ref{3d-eq-P}) by $3\sigma^2 P^2$ and integrating
the resulting equality  over $\mathbb{R}^3\times[0,T]$, one gets
that \be\label{3d-es-P-4}
&&\displaystyle\frac{3\gamma-1}{2\mu+\lambda}\int_0^T\sigma^2\|P\|_{L^4}^4\nonumber\\[2mm]
&=&\displaystyle-\sigma^2\|P\|_{L^3}^3+2\sigma\sigma'\int_0^{T}\|P\|_{L^3}^3-\frac{3\gamma-1}{2\mu+\lambda}\int_0^T\sigma^2\int_{\mathbb{R}^3} P^3G\nonumber\\[2mm]
&\le&\displaystyle
C\|P\|_{L^3}^3+\frac{3\gamma-1}{4\mu+2\lambda}\int_0^T\sigma^2\|P\|_{L^4}^4+\frac{C(3\gamma-1)}{2\mu+\lambda}\int_0^T\sigma^2\|G\|_{L^4}^4.
\ee The combination of (\ref{3d-lemma-u-4-e}) with (\ref{3d-es-P-4})
implies \be\label{3d-es-P-5}
\displaystyle\int_0^T\sigma^2\|P\|_{L^4}^4&\le&
C(2\mu+\lambda)(\gamma-1)E_0
+C(2\mu+\lambda)\mu^{\frac{1}{4}}(\gamma-1)^{\frac{3}{4}}E_0^{\frac{11}{12}}E_4\nonumber\\[2mm]
&\le&\displaystyle
(2\mu+\lambda)(\gamma-1)^{\frac{3}{4}}E_0^{\frac{11}{12}}\mu^{\frac{1}{4}}E_5,
\ee where $\displaystyle
E_5=C(\gamma-1)^\frac{2}{9}\mu^{-\frac{1}{6}}+CE_4$.

 Substituting
(\ref{3d-III 1}), (\ref{3d-III_2}) and (\ref{3d-es-P-5})  into
(\ref{3d-lemma-u-4}) shows that \be\label{3d-lemma-u-4-f}
\begin{aligned}
\int_0^T\int_{\mathbb{R}^3} \sigma^2|\nabla u|^4\le
(2\mu+\lambda)^{-3}\mu^\frac{1}{4}(\gamma-1)^{\frac{3}{4}}E_0^{\frac{11}{12}}E_6,
\end{aligned}
\ee where $\displaystyle E_6=CN_4^4E_4+CE_5$. Thus,
\be\label{3d-dt-II-4}
II_4=\displaystyle C\left(\frac{2\mu+\lambda}{\mu}\right)^2\int_0^T\int_{\mathbb{R}^3} \sigma^2 |\nabla u|^4
\le\displaystyle
\frac{(\gamma-1)^\frac{3}{4}E_0^\frac{11}{12}E_6}{(2\mu+\lambda)\mu^\frac{7}{4}}.
\ee To estimate $II_5$, using H\"older inequality,
(\ref{3d-basic-2}) and (\ref {3d-lemma-u-4-f}), we have
\be\label{3d-dt-II-5}
II_5&=&\displaystyle\frac{C(2\mu+\lambda)}{\mu}\int_0^T \sigma \|\nabla u\|_{L^3}^3\nonumber\\[3mm]
&\le& \displaystyle\frac{C(2\mu+\lambda)}{\mu} \left(\int_0^T\int_{\mathbb{R}^3}  |\nabla u|^2\right)^{\frac{1}{2}}\left(\int_0^T\int_{\mathbb{R}^3} \sigma^2 |\nabla u|^4\right)^{\frac{1}{2}}\nonumber\\[3mm]
&\le&\displaystyle
C\frac{(\gamma-1)^{\frac{3}{8}}E_0^{\frac{23}{24}}E_6^{\frac{1}{2}}}{\mu^{\frac{11}{8}}(2\mu+\lambda)^{\frac{1}{2}}}.
\ee For $II_6$, using H\"older inequality, (\ref{3d-es-P-5}) and
(\ref {3d-lemma-u-4-f}), one gets \be\label{3d-dt-II-6}
II_6&=&\displaystyle\frac{C\gamma^2}{\mu^2}\int_0^T\int_{\mathbb{R}^3} \sigma^2 |P\nabla u|^2\nonumber\\[3mm]
&\le& \displaystyle\frac{C\gamma^2}{\mu^2}\left(\int_0^{T}\int_{\mathbb{R}^3}\sigma^2 |P|^4\right)^{\frac{1}{2}} \left(\int_0^{T}\int_{\mathbb{R}^3}\sigma^2|\nabla u|^4\right)^{\frac{1}{2}}\nonumber\\[3mm]
&\le&\displaystyle
\frac{C\gamma^2E_5^{\frac{1}{2}}E_6^{\frac{1}{2}}(\gamma-1)^{\frac{3}{4}}E_0^{\frac{11}{12}}}{\mu^{\frac{7}{4}}(2\mu+\lambda)}.
\ee It holds from H\"older inequality, (\ref{3d-basic-2}) and
(\ref{3d-es-P-5})-(\ref {3d-lemma-u-4-f}) that \be\label{3d-dt-II-7}
II_7&=&\displaystyle\frac{C\gamma}{\mu}\int_0^{T}\int_{\mathbb{R}^3}
 \sigma P |\nabla u|^2\nonumber\\[3mm]
&\le&\displaystyle \frac{C\gamma}{\mu}\left(\int_0^{T}\int_{\mathbb{R}^3}\sigma^2 |P|^4\right)^{\frac{1}{4}} \left(\int_0^{T}\int_{\mathbb{R}^3}\sigma^2 |\nabla u|^4\right)^{\frac{1}{4}}\left(\int_0^{T}\int_{\mathbb{R}^3} |\nabla u|^2\right)^{\frac{1}{2}}\nonumber\\[3mm]
&\le&\displaystyle \frac{C\gamma
E_5^{\frac{1}{4}}E_6^{\frac{1}{4}}(\gamma-1)^{\frac{3}{8}}E_0^{\frac{23}{24}}}{\mu^{\frac{11}{8}}(2\mu+\lambda)^{\frac{1}{2}}}.
\ee Finally, relations \eqref{3d-basic-3},
(\ref{3d-dt-A-1+A-2}) and (\ref{3d-dt-II-4})-(\ref{3d-dt-II-7}) give rise to
\be\label{3d-dt-A-1+A-2-f} A_1(T)+A_2(T)
&\le&\displaystyle\frac{C(\gamma-1)E_0}{\mu^2}
+\frac{C(\gamma-1)^\frac{3}{4}E_0^\frac{11}{12}E_6}{(2\mu+\lambda)\mu^\frac{7}{4}}
+\frac{C(\gamma-1)^{\frac{3}{8}}E_0^{\frac{23}{24}}E_6^{\frac{1}{2}}}{\mu^{\frac{11}{8}}(2\mu+\lambda)^{\frac{1}{2}}}\nonumber\\[3mm]
&&\displaystyle+\frac{C\gamma^2E_5^{\frac{1}{2}}E_6^{\frac{1}{2}}(\gamma-1)^{\frac{3}{4}}E_0^{\frac{11}{12}}}{\mu^{\frac{7}{4}}(2\mu+\lambda)}
+\frac{C \gamma E_5^{\frac{1}{4}}E_6^{\frac{1}{4}}(\gamma-1)^{\frac{3}{8}}E_0^{\frac{23}{24}}}{\mu^{\frac{11}{8}}(2\mu+\lambda)^{\frac{1}{2}}}\nonumber\\[3mm]
&&\displaystyle+\frac{C(2\mu+\lambda)(\gamma-1)^\frac{2}{3}E_0^\frac{2}{3}E_1}{\mu^2}\nonumber\\[3mm]
&\le&\displaystyle C
\left(\frac{(\gamma-1)^{\frac{1}{6}}E_0^{\frac{1}{2}}}{\mu^{\frac{1}{2}}}\right)^{\frac{4}{3}}E_7,
\ee where \be\label{3d-dt-E-8} E_7&=&\displaystyle
\frac{(\gamma-1)^\frac{2}{3}}{\mu}+\frac{(\gamma-1)^\frac{4}{9}E_6}{\mu^\frac{11}{6}}
+\frac{(\gamma-1)^\frac{1}{18}E_6^\frac{1}{2}}{\mu^\frac{11}{12}}
+\frac{(\gamma-1)^\frac{4}{9}E_5^\frac{1}{2}E_6^\frac{1}{2}}{\mu^\frac{11}{6}}\nonumber\\[3mm]
&&\displaystyle+\frac{\gamma(\gamma-1)^\frac{1}{18}E_5^\frac{1}{4}E_6^\frac{1}{4}}{\mu^{\frac{11}{12}}}
+\left(2+\frac{\lambda}{\mu}\right)\frac{(\gamma-1)^\frac{4}{9}E_1}{\mu^\frac{1}{3}},
\ee and we have also used the facts that
$\displaystyle\frac{(\gamma-1)^\frac{1}{3}E_0}{\mu}\le
1$ and $\mu+\lambda>0$. It thus follows from (\ref{3d-dt-A-1+A-2-f})
that \be\label{3d-dt-A-1+A-2-dd}
\begin{aligned}
A_1(T)+A_2(T) \le
C\left(\frac{(\gamma-1)^{\frac{1}{6}}E_0^{\frac{1}{2}}}{\mu^{\frac{1}{2}}}\right)^{\frac{4}{3}}E_7\le
\frac{(\gamma-1)^{\frac{1}{6}}E_0^{\frac{1}{2}}}{\mu^{\frac{1}{2}}},
\end{aligned}
\ee provided $\displaystyle
\frac{(\gamma-1)^{\frac{1}{6}}E_0^{\frac{1}{2}}}{\mu^{\frac{1}{2}}}\le
(CE_7)^{-3}$.

\bigbreak

\noindent\underline{\textbf{Case 2:}} \ $\displaystyle\gamma > \frac{3}{2}$.

In view of \eqref{3d-lemma-u-4}, one gets
\be\label{3d-lemma-u-41}
\int_0^T\int_{\mathbb{R}^3} \sigma^2|\nabla
u|^4\leq\sum\limits_{i=1}^3III_i. \ee
 For $III_1$, using \eqref{a priori assumption 1} and Lemma \ref{bound}, we have
\be\label{3d-IIII 1}
III_1&\le&\displaystyle CN_4^4(2\mu+\lambda)^{-3}\sup\limits_{0\le t\le T}(\sigma^2\|\rho\u\|_{L^2}^2)^{\frac{1}{2}}\int_0^{T}\sigma\|\rho\u\|_{L^2}^2\|\nabla u\|_{L^2}\nonumber\\[2mm]
&\le&\displaystyle CN_4^4(2\mu+\lambda)^{-3}\mu^{\frac{1}{2}}A^{\frac{1}{2}}_2(T)\sup\limits_{0\le t\le T}(\|\nabla u\|_{L^2}^2)^{\frac{1}{2}}\int_0^{T}\sigma\|\rho\u\|_{L^2}^2\nonumber\\[2mm]
&\le&\displaystyle CN_4^4(2\mu+\lambda)^{-3}\mu^{\frac{3}{2}}A^{\frac{1}{2}}_2(T)A_1(T)(E_2+1)^\frac{1}{2}\nonumber\\[2mm]
&\le&\displaystyle
CN_4^4(2\mu+\lambda)^{-3}\mu^{\frac{3}{4}}(\gamma-1)^\frac{1}{4}E_0^\frac{3}{4}(E_2+1)^\frac{1}{2}.
\ee It follows from \eqref{a priori assumption 1},
(\ref{3d-basic-1}) and Lemma \ref{bound} that \be\label{3d-IIII_2}
III_2&\le &\displaystyle C(2\mu+\lambda)^{-4}\sup\limits_{0\le t\le T}(\sigma^2\|\rho\u\|_{L^2}^2)^{\frac{1}{2}}\int_0^{T}\sigma\|\rho\u\|_{L^2}^2\|P\|_{L^2}\nonumber\\[2mm]
&\le&\displaystyle
C(2\mu+\lambda)^{-4}\mu^{\frac{3}{4}}(\gamma-1)^{\frac{3}{4}}E_0^\frac{5}{4}(E_2+1)^\frac{1}{2}.
\ee Thus, one gets from \eqref{hp1}, (\ref{3d-IIII 1}) and
(\ref{3d-IIII_2}) that
 \be\label{3d-lemma-u-4-e1}
\begin{aligned}
\int_0^T\sigma^2\|G\|_{L^4}^4\le(2\mu+\lambda)\mu^\frac{3}{4}(\gamma-1)^\frac{1}{4}E_0^\frac{3}{4}E_8,
\end{aligned}
\ee where \be\label{e8} \displaystyle
E_8=C(E_2+1)^\frac{1}{2}\left(1+\frac{(\gamma-1)^{\frac{1}{3}}}{\mu^{\frac{1}{2}}}\right).
\ee
Using \eqref{3d-es-P-4}, $\displaystyle \frac{1}{\gamma-1} < 2$ and
$\displaystyle\frac{(\gamma-1)^\frac{1}{3}E_0}{\mu}\le1$,
we obtain \be\label{3d-es-P-41}
\displaystyle\int_0^T\sigma^2\|P\|_{L^4}^4&\le& C(2\mu+\lambda)(\gamma-1)E_0+C(3\gamma-1)(2\mu+\lambda)\mu^\frac{3}{4}(\gamma-1)^\frac{1}{4}E_0^\frac{3}{4}E_8\nonumber\\[2mm]
&\le&\displaystyle C(2\mu+\lambda)\mu^\frac{3}{4}(\gamma-1)E_0^\frac{3}{4}\Big[E_0^\frac{1}{4}\mu^{-\frac{3}{4}}+(3\gamma-1)(\gamma-1)^{-\frac{3}{4}}E_8\Big]\nonumber\\[2mm]
&\le&\displaystyle C(2\mu+\lambda)\mu^\frac{3}{4}(\gamma-1)E_0^\frac{3}{4}\Big[E_0^\frac{1}{4}\mu^{-\frac{3}{4}}+2^{\frac{3}{4}}(3\gamma-1)E_8\Big]\nonumber\\[2mm]
&\le&\displaystyle
(2\mu+\lambda)\mu^\frac{3}{4}(\gamma-1)E_0^\frac{3}{4}E_{9}, \ee
where $\displaystyle E_{9}=\displaystyle
C\mu^{-\frac{1}{2}}+C(3\gamma-1)E_8$.

 Substituting \eqref{3d-IIII 1}, \eqref{3d-IIII_2} and
\eqref{3d-es-P-41} into \eqref{3d-lemma-u-41}, we get
\be\label{3d-lemma-u-4-f1}
\begin{aligned}
\int_0^T\int_{\mathbb{R}^3} \sigma^2|\nabla u|^4\le
(2\mu+\lambda)^{-3}\mu^\frac{3}{4}(\gamma-1)^\frac{1}{4}E_0^\frac{3}{4}E_{10},
\end{aligned}
\ee where $\displaystyle
E_{10}=CN_4^4(E_2+1)^\frac{1}{2}+C\mu^{\frac{1}{2}}(\gamma-1)^\frac{1}{3}(E_2+1)^\frac{1}{2}+C(\gamma-1)^\frac{3}{4}E_{9}$.
Thus, \be\label{3d-dt-II-41} II_4\le
C\Big(\frac{2\mu+\lambda}{\mu}\Big)^{2}(2\mu+\lambda)^{-3}\mu^\frac{3}{4}(\gamma-1)^\frac{1}{4}
E_0^\frac{3}{4}E_{10}\le\frac{C(\gamma-1)^\frac{1}{4}E_0^\frac{3}{4}E_{10}}{\mu^\frac{5}{4}(2\mu+\lambda)}.
\ee Using the same spirit of deriving
\eqref{3d-dt-II-5}-\eqref{3d-dt-II-7}, we obtain
\begin{equation}\label{3d-dt-II-8}
\begin{cases}
II_5&\displaystyle\le \frac{C (2\mu+\lambda)}{\mu}\Big(\frac{E_0}{\mu}\Big)^\frac{1}{2}\Big((2\mu+\lambda)^{-3}\mu^\frac{3}{4}(\gamma-1)^\frac{1}{4}E_0^\frac{3}{4}E_{10}\Big)^\frac{1}{2},\\[3mm]
&\displaystyle\le\frac{C(\gamma-1)^\frac{1}{8}E_0^\frac{7}{8}E_{10}^\frac{1}{2}}{\mu^\frac{9}{8}(2\mu+\lambda)^\frac{1}{2}},\\[4.5mm]
\displaystyle II_6&\displaystyle\le\frac{C\gamma^2(\gamma-1)^\frac{5}{8}E_0^\frac{3}{4}E_{9}^\frac{1}{2}E_{10}^\frac{1}{2}}{\mu^{\frac{5}{4}}(2\mu+\lambda)},\\[4.5mm]
II_7&\displaystyle\le \frac{C\gamma(\gamma-1)^\frac{5}{16}E_0^\frac{7}{8}E_{9}^\frac{1}{4}E_{10}^\frac{1}{4}}{\mu^{\frac{9}{8}}(2\mu+\lambda)^\frac{1}{2}},\\[4.5mm]
II_8&\displaystyle\le\frac{C(2\mu+\lambda)E_0}{\mu^2}\le \frac{C
(2\mu+\lambda)(\gamma-1)E_0}{\mu^2}.
\end{cases}
\end{equation}

Consequently, after a bit tedious but straightforward calculation, relations
\eqref{3d-basic-3}, \eqref{3d-dt-A-1+A-2} and \eqref{3d-dt-II-8} give rise to
\be\label{3d-dt-A-1+A-2-f1}
A_1(T)+A_2(T)&\le&\displaystyle\frac{C(\gamma-1)E_0}{\mu^2}
+\frac{C(\gamma-1)^\frac{1}{4}E_0^\frac{5}{4}E_{10}}{\mu^\frac{3}{4}(2\mu+\lambda)}
+\frac{C(\gamma-1)^\frac{1}{8}E_0^\frac{7}{8}E_{10}^\frac{1}{2}}{\mu^\frac{9}{8}(2\mu+\lambda)^\frac{1}{2}}\nonumber\\[2mm]
&&\displaystyle+\frac{C\gamma^2(\gamma-1)^\frac{5}{8}E_0^\frac{3}{4}E_{9}^\frac{1}{2}E_{10}^\frac{1}{2}}{\mu^{\frac{5}{4}}(2\mu+\lambda)}
+\frac{C\gamma(\gamma-1)^\frac{5}{16}E_0^\frac{7}{8}E_{9}^\frac{1}{4}E_{10}^\frac{1}{4}}{\mu^{\frac{9}{8}}(2\mu+\lambda)^\frac{1}{2}}\nonumber\\[2mm]
&&\displaystyle+\frac{C(2\mu+\lambda)(\gamma-1)E_0}{\mu^2}\nonumber\\[2mm]
&\le&\displaystyle C\left(\frac{
(\gamma-1)^{\frac{1}{6}}E_0^{\frac{1}{2}}}{\mu^{\frac{1}{2}}}\right)^{\frac{3}{2}}E_{11},
\ee where \be\label{3d-dt-E-13} \displaystyle E_{11}&=&\displaystyle
\frac{(\gamma-1)^\frac{2}{3}}{\mu} +\frac{E_{10}}{\mu}
+\frac{E_{10}^{\frac{1}{2}}}{\mu^{\frac{3}{4}}}
+\frac{\gamma^2(\gamma-1)^\frac{3}{8}E_{9}^\frac{1}{2}E_{10}^\frac{1}{2}}{\mu^{\frac{3}{2}}}
+\frac{\gamma(\gamma-1)^{\frac{1}{48}}E_{9}^\frac{1}{4}E_{10}^\frac{1}{4}}{\mu^{\frac{3}{4}}}\nonumber\\[2mm]
&&\displaystyle+\left(2+\frac{\lambda}{\mu}\right)(\gamma-1)^\frac{2}{3}.
\ee Here we have used the facts that $\displaystyle
\frac{1}{\gamma-1} < 2$,
$\displaystyle\frac{(\gamma-1)^\frac{1}{3}E_0}{\mu}\le
1$ and $\mu+\lambda>0$. It thus follows from
(\ref{3d-dt-A-1+A-2-f1}) that \be\label{3d-dt-A-1+A-2-c1-ee}
A_1(T)+A_2(T) \le
C\left(\frac{(\gamma-1)^{\frac{1}{6}}E_0^{\frac{1}{2}}}{\mu^{\frac{1}{2}}}\right)^{\frac{3}{2}}E_{11}
\le\frac{(\gamma-1)^{\frac{1}{6}}E_0^{\frac{1}{2}}}{\mu^{\frac{1}{2}}},
\ee provided $\displaystyle
\frac{(\gamma-1)^{\frac{1}{6}}E_0^{\frac{1}{2}}}{\mu^{\frac{1}{2}}}\le
(CE_{11})^{-2}$.

By (\ref{3d-dt-A-1+A-2-dd}) and (\ref{3d-dt-A-1+A-2-c1-ee}), for
\textbf{Case 1} and \textbf{Case 2}, we conclude that if
$$ \frac{(\gamma-1)^\frac{1}{6}E_0^{\frac{1}{2}}}{\mu^{\frac{1}{2}}}\le
\min\left\{(CE_7)^{-3}\Big|_{(1< \gamma\leq \frac{3}{2})},
(CE_{11})^{-2}\Big|_{( \gamma>\frac{3}{2})}, \varepsilon_1,
\varepsilon_2 \right\}\triangleq \varepsilon_3,
$$
then \be\label{3d-dt-A-1+A-2-ff} A_1(T)+A_2(T)
\le\frac{(\gamma-1)^{\frac{1}{6}}E_0^{\frac{1}{2}}}{\mu^{\frac{1}{2}}}.
\ee The proof of
Lemma \ref{3d-le:3.8} is completed.

\endpf

  Next, we will derive the time-independent upper bound for the density $\rho$.
  The approach is motivated by Huang-Li-Xin in \cite{Huang-Li-Xin} and Li-Xin in \cite{lz}.
\begin{lemma}\label{3d-le:rho}Under the conditions of Proposition \ref{prop 3.1}, it holds that
\be\label{3d-upper bound of rho}
\sup_{0\leq t\leq T}\|\rho\|_{L^\infty}\leq\frac{7\bar{\rho}}{4}
\ee for any $(x,t)\in\mathbb{R}^3\times[0,T]$,
provided
\bex \begin{aligned}\frac{(\gamma-1)^\frac{1}{6}E_0^{\frac{1}{2}}}{\mu^{\frac{1}{2}}} \leq \min\left\{\varepsilon_3,(2C(\bar{\rho},M))^{-\frac{16}{3}}\mu^{4},(4C(\bar{\rho}))^{-2}
\right\}
\triangleq &  \ \varepsilon.\end{aligned} \eex
\end{lemma}
\pf  Denoting $D_t\rho=\rho_t+u\cdot\nabla\rho$ and expressing the equation of the mass conservation (\ref{3d-full N-S})$_1$ as
$$
D_t\rho=g(\rho)+b'(t),
$$
where
$$
g(\rho)\triangleq -\frac{\rho P}{2\mu+\lambda},\ \ \
b(t)\triangleq-\frac{1}{2\mu+\lambda}\int_0^t\rho Gd\tau.
$$
In accordance with Lemma \ref{GN}, \eqref{hp4}, \eqref{hp1}, \eqref{A2 1}
and Lemma \ref{3d-le:3.8}, for all $0 \leq t_1< t_2 \leq\sigma(T)$, we get
\be\label{b1}
 \displaystyle |b(t_2)-b(t_1)|&\leq&\displaystyle C\int_0^{\sigma(T)}\|(\rho G)(\cdot,t)\|_{L^\infty} dt\nonumber\\[2mm]
    &\leq&\displaystyle \frac{C(\bar{\rho})}{2\mu+\lambda}\int_0^{\sigma(T)}\|G(\cdot,t)\|_{L^6}^{\frac{1}{2}}\|\nabla G(\cdot,t)\|_{L^6}^{\frac{1}{2}}dt\nonumber\\[2mm]
    &\leq&\displaystyle \frac{C(\bar{\rho})}{2\mu+\lambda}\int_0^{\sigma(T)}\|\rho \dot{u}\|_{L^2}^\frac{1}{2}\|\nabla \dot{u}\|_{L^2}^\frac{1}{2}dt\nonumber\\[2mm]
    &\leq&\displaystyle \frac{C(\bar{\rho})}{2\mu+\lambda}\left(\int_0^{\sigma(T)}\|\rho \dot{u}\|_{L^2}^\frac{2}{3}t^{-\frac{1}{3}}dt\right)^{\frac{3}{4}}
                         \left(\int_0^{\sigma(T)}t\|\nabla \dot{u}\|_{L^2}^2dt\right)^{\frac{1}{4}}\nonumber\\[2mm]
    &\leq&\displaystyle \frac{C(\bar{\rho},M)}{2\mu+\lambda}\left(\int_0^{\sigma(T)}(\|\rho \dot{u}\|_{L^2}^2t)^{\frac{1}{3}}t^{-\frac{2}{3}}dt\right)^{\frac{3}{4}}\nonumber\\[2mm]
    &\leq&\displaystyle \frac{C(\bar{\rho},M)}{2\mu+\lambda}\sup_{0\leq t \leq \sigma(T)}\left(\|\rho \dot{u}\|_{L^2}^2t\right)^{\frac{1}{16}}\left(\int_0^{\sigma(T)}(\|\rho \dot{u}\|_{L^2}^2t)^\frac{1}{4}t^{-\frac{2}{3}}dt\right)^{\frac{3}{4}}\nonumber\\[2mm]
    &\leq&\displaystyle \frac{C(\bar{\rho},M)\mu^\frac{1}{16}}{2\mu+\lambda}\left(\int_0^{\sigma(T)}\|\rho \dot{u}\|_{L^2}^2t\right)^{\frac{3}{16}}\left(\int_0^{\sigma(T)}t^{-\frac{8}{9}}dt\right)^{\frac{9}{16}}\nonumber\\[2mm]
    &\leq&\displaystyle\frac{ C(\bar{\rho},M)}{\mu^\frac{3}{4}}A_{1}(\sigma (T))^{\frac{3}{16}}\nonumber\\[2mm]
    &\leq&\displaystyle\frac{ C(\bar{\rho},M)}{\mu^\frac{3}{4}}\Big(\frac{(\gamma-1)^\frac{1}{6}E_0^{\frac{1}{2}}}{\mu^{\frac{1}{2}}}\Big)^{\frac{3}{16}},
\ee
provided $\displaystyle \frac{(\gamma-1)^\frac{1}{6}E_0^{\frac{1}{2}}}{\mu^{\frac{1}{2}}}\leq \varepsilon_3$.
Therefore, for $t\in [0,\sigma(T)]$, one can choose $N_0$ and $N_{1}$ in Lemma \ref{zlo} as follows:
$$
N_{1}=0,\ \ \ \ \ \ N_0=\frac{ C(\bar{\rho},M)}{\mu^\frac{3}{4}}\left(\frac{(\gamma-1)^\frac{1}{6}E_0^{\frac{1}{2}}}{\mu^{\frac{1}{2}}}\right)^{\frac{3}{16}},
$$
and $\bar{\zeta}=0$. Then
$$
g(\zeta)=-\frac{\zeta P(\zeta)}{2\mu+\lambda}\leq -N_1=0  \ \  \mathrm{for\ \ all} \ \ \zeta\geq\bar{\zeta} =0.
$$
Thus
\be\label{md1}
\sup_{0\leq t \leq \sigma(T)}\|\rho\|_{L^\infty}\leq \max\Big\{\bar{\rho},0\Big\}+N_0\leq\bar{\rho}
+\frac{C(\bar{\rho},M)}{\mu^\frac{3}{4}}\left(\frac{(\gamma-1)^\frac{1}{6}E_0^{\frac{1}{2}}}{\mu^{\frac{1}{2}}}\right)^{\frac{3}{16}}\leq\frac{3\bar{\rho}}{2},
\ee
provided
\be\label{a1}
\frac{(\gamma-1)^\frac{1}{6}E_0^{\frac{1}{2}}}{\mu^{\frac{1}{2}}}\leq \min\left\{\varepsilon_3, (2C(\bar{\rho},M))^{-\frac{16}{3}}\mu^{4}\right\}.
\ee

On the other hand, by virtues of Lemma \ref{GN}, \eqref{hp4}, \eqref{A2 1}
and Lemma \ref{3d-le:3.8}, for $t\in [\sigma(T), T]$, one deduces that
\be\label{b2}
\displaystyle|b(t_2)-b(t_1)|&\leq&\displaystyle \frac{C(\bar{\rho})}{2\mu+\lambda}\int_{t_1}^{t_2}\|G(\cdot,t)\|_{L^\infty}dt\nonumber\\[2mm]
&\leq&\displaystyle \frac{1}{2\mu+\lambda}(t_2-t_1)+\frac{C(\bar{\rho})}{2\mu+\lambda}\int_{\sigma(T)}^T\|G\|_{L^\infty}^{4}\nonumber\\[2mm]
&\leq&\displaystyle \frac{1}{2\mu+\lambda}(t_2-t_1)+\frac{C(\bar{\rho})}{2\mu+\lambda}\int_{\sigma(T)}^T\|G\|_{L^6}^{2}\|\nabla G\|_{L^6}^{2}\nonumber\\[2mm]
&\leq&\displaystyle \frac{1}{2\mu+\lambda}(t_2-t_1)+\frac{C(\bar{\rho})}{2\mu+\lambda}\int_{\sigma(T)}^T\|\nabla G\|_{L^2}^{2}\|\nabla \dot{u}\|_{L^2}^{2}\nonumber\\[2mm]
&\leq&\displaystyle \frac{1}{2\mu+\lambda}(t_2-t_1)+\frac{C(\bar{\rho})}{2\mu+\lambda}\int_{\sigma(T)}^T\|\rho\dot{u}\|_{L^2}^{2}\|\nabla \dot{u}\|_{L^2}^{2}\nonumber\\[2mm]
&\leq&\displaystyle \frac{1}{2\mu+\lambda}(t_2-t_1)+C(\bar{\rho})A_2(T)\int_{\sigma(T)}^T\|\nabla \dot{u}\|_{L^2}^{2}\nonumber\\[2mm]
&\leq&\displaystyle \frac{1}{2\mu+\lambda}(t_2-t_1)+C(\bar{\rho})A^2_2(T)\nonumber\\[2mm]
&\leq&\displaystyle \frac{1}{2\mu+\lambda}(t_2-t_1)+C(\bar{\rho})\Big(\frac{(\gamma-1)^\frac{1}{6}E_0^{\frac{1}{2}}}{\mu^{\frac{1}{2}}}\Big)^2,
    \ee
provided $\displaystyle \frac{(\gamma-1)^\frac{1}{6}E_0^{\frac{1}{2}}}{\mu^{\frac{1}{2}}}\leq \varepsilon_3$.

Therefore, one can choose $N_1$ and $N_0$ in Lemma \ref{zlo} as
$$
N_1=\frac{1}{2\mu+\lambda},\ \ \ \ \ \ N_0=C(\bar{\rho})\left(\frac{(\gamma-1)^\frac{1}{6}E_0^{\frac{1}{2}}}{\mu^{\frac{1}{2}}}\right)^2.
$$
Note that
$$
g(\zeta)=-\frac{\zeta P(\zeta)}{2\mu+\lambda}\leq -N_1=-\frac{1}{2\mu+\lambda} \ \  \mathrm{for\ \ all} \ \ \zeta\geq 1,
$$
one can set $\bar{\zeta}=1$. Thus
\be\label{md2}
\begin{aligned}
\sup_{\sigma(T)\leq s \leq T}\|\rho\|_{L^\infty}\leq \max\left\{\frac{3}{2}\bar{\rho},1\right\}
+N_0 &\leq \frac{3}{2}\bar{\rho}+C(\bar{\rho})\left(\frac{(\gamma-1)^\frac{1}{6}E_0^{\frac{1}{2}}}{\mu^{\frac{1}{2}}}\right)^2\\[3mm]
&\le \frac{7\bar{\rho}}{4},
\end{aligned}
\ee
provided
\be\label{a2}
\frac{(\gamma-1)^\frac{1}{6}E_0^{\frac{1}{2}}}{\mu^{\frac{1}{2}}}\leq \min\Big\{\varepsilon_3, (2C(\bar{\rho},M))^{-\frac{16}{3}}\mu^{4}, (4C(\bar{\rho}))^{-2}\Big\}.
\ee
The combination of \eqref{md1} and \eqref{md2} completes the proof of Lemma \ref{3d-le:rho}.
\endpf

\bigbreak

\noindent \textbf{\textit{Proof of Theorem}} \ref{3d-th:1.1}: With Lemma \ref{3d-le:rho} and the higher norm estimates of the smooth solution $(\rho, u)$
in \cite{Huang-Li-Xin}, Theorem \ref{3d-th:1.1} follows.

\section{The proof of Theorem 1.2}

 \setcounter{equation}{0}\setcounter{theorem}{0}
\renewcommand{\theequation}{\thesection.\arabic{equation}}
\renewcommand{\thetheorem}{\thesection.\arabic{theorem}}
\setcounter{equation}{0}\setcounter{theorem}{0}
\renewcommand{\theequation}{\thesection.\arabic{equation}}
\renewcommand{\thetheorem}{\thesection.\arabic{theorem}}

 In this section, we will prove Theorem 1.2. Let
\be\label{A2}
\begin{cases}
\displaystyle A_1(T)=\sup\limits_{0\le t\le T}\sigma\int_{\mathbb{R}^3}|\nabla
u|^2+\int_0^T\int_{\mathbb{R}^3}\sigma\frac{\rho|\u|^2}{\mu},\\[3mm]
\displaystyle A_2(T)=\sup\limits_{0\le t\le
T}\sigma^3\int_{\mathbb{R}^3}\frac{\rho
|\u|^2}{\mu}+\int_0^T\int_{\mathbb{R}^3}\sigma^3|\nabla\u|^2,\\[3mm]
\displaystyle A_3(T)=\sup\limits_{0\le t\le T}\int_{\mathbb{R}^3}\frac{\rho|u|^3}{\mu^3}.
\end{cases}
\ee

 Throughout the
rest of the paper, we denote generic constant by $C$ depending on
$\bar{\rho}, M$ and some other known constants but independent of
$\mu,\lambda,\gamma-1,\wi{\rho}$ and $t$, and we write $C(\alpha)$ to
emphasize that $C$ may depend on $\alpha$. For simplicity of
presentation, we shall assume that \be\label{assumpution rb}
\displaystyle\frac{ (\gamma-1+\wi{\rho})E_0^\alpha}{\mu^\beta} \le
1, \ee where $\displaystyle \ \frac{1}{10}\leq \alpha\leq 100 ,\
\frac{2}{3}\leq \beta\leq 12$. It's worth noting that the assumption
\eqref{assumpution rb} is not essential for our paper. Without this
assumption, $E_{12}$-$E_{21}$ in the proof of Theorem
\ref{3d-th:1.2} should be more complex.
The following proposition plays a crucial role in this section.
\begin{proposition}\label{prop 4.1}
Assume that the initial data satisfies (\ref{initial data2}),
(\ref{3d-initial assumption2}) and (\ref{3d-compatibility2}),
$1<\gamma<2$. If the solution $(\rho, u)$ satisfies \be\label{a
priori assumption}
\begin{cases}
\displaystyle 0\le \rho\le 2\bar{\rho},\ \ \ \ \ \ \ \ \ \ \ \ \ \ \ \ \ \ \ \ \ \ \ \ \  A_1(T)\le 2\left\{\frac{\left((\gamma-1)^{\frac{1}{36}}+\wi{\rho}^{\frac{1}{6}}\right)E_0^\frac{1}{4}}{\mu^\frac{1}{3}}\right\}^\frac{3}{8}, \\[3mm]
\displaystyle A_2(T)\le\frac{2\left((\gamma-1)^\frac{1}{36}+\wi{\rho}^{\frac{1}{6}}\right)E_0^{\frac{1}{4}}}{\mu^{\frac{1}{3}}},\ A_3(\sigma(T))\le
2\left\{\frac{\left((\gamma-1)^{\frac{1}{36}}+\wi{\rho}^{\frac{1}{6}}\right)E_0^\frac{1}{4}}{\mu^\frac{1}{3}}\right\}^\frac{1}{2},
\end{cases}
\ee
then

 \be\label{4d-result from the a p}
 \begin{cases}
\displaystyle 0\le \rho\le \frac{7}{4}\bar{\rho},\ \ \ \ \ \ \ \ \ \ \ \ \ \ \ \ \ \ \ \ \ \ \ \ \ A_1(T)\le \left\{\frac{\left((\gamma-1)^{\frac{1}{36}}+\wi{\rho}^{\frac{1}{6}}\right)E_0^\frac{1}{4}}{\mu^\frac{1}{3}}\right\}^\frac{3}{8}, \\[3mm]
\displaystyle A_2(T)\le\frac{\left((\gamma-1)^\frac{1}{36}+\wi{\rho}^{\frac{1}{6}}\right)E_0^{\frac{1}{4}}}{\mu^{\frac{1}{3}}},\ A_3(\sigma(T))\le
\left\{\frac{\left((\gamma-1)^{\frac{1}{36}}+\wi{\rho}^{\frac{1}{6}}\right)E_0^\frac{1}{4}}{\mu^\frac{1}{3}}\right\}^\frac{1}{2},
\end{cases}
\ee
$(x,t)\in \mathbb{R}^3\times[0,T]$, provided $\displaystyle\frac{(\gamma-1)^\frac{1}{36}E_0^{\frac{1}{4}}}{\mu^{\frac{1}{3}}}\le \frac{\tilde{\rho}}{2 C}$ and $\displaystyle\frac{\left((\gamma-1)^\frac{1}{36}+\wi{\rho}^{\frac{1}{6}}\right)E_0^{\frac{1}{4}}}{\mu^{\frac{1}{3}}}\le\varepsilon$. Here
\bex
\begin{aligned}
\varepsilon= &\min\left\{\varepsilon_6,
(2C(\bar{\rho},M))^{-\frac{16}{3}}\mu^{4},
(4C(\bar{\rho}))^{-2}\right\},\end{aligned} \eex and
\bex\begin{aligned}
&\varepsilon_6=\min\left\{\Big(C(E_{18}+E_{19}+E_{20})\Big)^{-17},
\Big(C(E_{18}+E_{19}+E_{21})\Big)^{-8}, \varepsilon_5
\right\},\\[2mm]
&\varepsilon_5=\min \left\{
\Big(C(E_{15}E_{17}+E_{16})\Big)^{-4},\ \varepsilon_4\right\},\\[2mm]
&\varepsilon_4=\min\left\{\Big(4C(\bar{\rho})\Big)^{-6}, 1 \right\}.
\end{aligned}\eex
\end{proposition}

\pf Proposition \ref{prop 4.1} follows from Lemmas
\ref{3d-le-rho3}-\ref{4d:rho} below.

\begin{lemma}\label{3d-le-rho3}
Let $(\rho, u)$ be a smooth solution of (\ref{3d-full
N-S})-(\ref{3d-boundary-0}) with $0\le \rho \le 2\bar{\rho}$,
$\bar{\rho}\geq\wi{\rho}+1$ and $0< \gamma-1 < 1$. Then there exists
a positive constant $C$ such that \be\label{3d-le-rho3-1}
\displaystyle\int_{\mathbb{R}^3}|\rho-\wi{\rho}|^3\le
C(\bar{\rho})(\gamma-1)^{\frac{1}{4}}E_0.\ee
\end{lemma}
\pf A straightforward calculation implies that $$\displaystyle G(\rho)=\rho\int_{\widetilde{\rho}}^{\rho}\frac{P(s)-P(\widetilde{\rho})}{s^2}ds=\frac{1}{\gamma-1}[\rho^{\ga}-\wi{\rho}^{\ga}-\ga\wi{\rho}^{\ga-1}(\rho-\wi{\rho})].$$
Next, we claim that
\be\label{3d-G}
G(\rho)\geq\begin{cases}(\ga-1)^{-\frac{1}{4}}|\rho-\wi{\rho}|^{\ga-1}, \ \ \ \ \ \ |\rho-\wi{\rho}|>(\ga-1)^{\frac{1}{3}},\\[2mm]
(\ga-1)^{-\frac{1}{4}}|\rho-\wi{\rho}|^3, \ \ \ \ \ \ \ \ \ |\rho-\wi{\rho}|\le(\ga-1)^{\frac{1}{3}}.
\end{cases}
\ee

\noindent\underline{\textbf{Case 1:}} $|\rho-\wi{\rho}|>
(\ga-1)^{\frac{1}{3}}.$

Without loss of generality, we only consider the case of $\rho-\wi{\rho}> (\ga-1)^{\frac{1}{3}}$. Now we define $$f(\rho)=\rho^{\ga}-\wi{\rho}^{\ga}-\ga\wi{\rho}^{\ga-1}(\rho-\wi{\rho})-(\ga-1)^{\frac{3}{4}}(\rho-\wi{\rho})^{\ga-1},$$
and note that
$$ f(\wi{\rho})=0, \ \ \ \ f'(\rho)=\ga\rho^{\ga-1}-\ga\wi{\rho}^{\ga-1}-(\ga-1)^{\frac{7}{4}}(\rho-\wi{\rho})^{\ga-2}.$$
Thus
\be\label{3d-le-f} f(\rho)=f'(\xi)(\rho-\wi{\rho}),\ee
where $\xi=(1-\theta)\rho+\theta\wi{\rho}\ \  (0<\theta<1)$.
Then
\be
f'(\xi)&=&\ga\xi^{\ga-1}-\ga\wi{\rho}^{\ga-1}-(\ga-1)^{\frac{7}{4}}(\xi-\wi{\rho})^{\ga-2} \nonumber \\[2mm]
&=&\ga(\ga-1)(\xi-\wi{\rho})\left[(1-\theta_1)\xi+\theta_1\wi{\rho}\right]^{\ga-2}-(\ga-1)^{\frac{7}{4}}(\xi-\wi{\rho})^{\ga-2} \ \  (0<\theta_1<1)\nonumber\\[2mm]
&=&(\ga-1)(1-\theta)(\rho-\wi{\rho})\left[\ga((1-\theta_1)\xi+\theta_1\wi{\rho})^{\ga-2}-(\ga-1)^{\frac{3}{4}}(\xi-\wi{\rho})^{\ga-3}\right] \nonumber\\[2mm]
&>&(\ga-1)(1-\theta)(\rho-\wi{\rho})\left[C\bar{\rho}^{\ga-2}-C(\ga-1)^{\frac{3}{4}}(\ga-1)^{\frac{\ga-3}{3}}\right]\nonumber \\[2mm]
&>&(\ga-1)(1-\theta)(\rho-\wi{\rho})\left[C\bar{\rho}^{\ga-2}-C(\ga-1)^{\frac{1}{12}}\right]
>0, \ \ \ \ (1<\ga<2),\nonumber \ee
when $\displaystyle\ga-1\le
C\bar{\rho}^{12(\ga-2)}\le\frac{C}{\bar{\rho}^{12}}=\eta.$ This
together with (\ref{3d-le-f}) implies the first inequality of
(\ref{3d-G}).

\vspace{2ex}
\noindent\underline{\textbf{Case 2:}} $\displaystyle
|\rho-\wi{\rho}|\le (\ga-1)^{\frac{1}{3}}.$

Similar to \textbf{Case 1}, we only consider the case of $0\leq\rho-\wi{\rho}\le (\ga-1)^{\frac{1}{3}}$. Now let
$$g(\rho)=\rho^{\ga}-\wi{\rho}^{\ga}-\ga\wi{\rho}^{\ga-1}(\rho-\wi{\rho})-(\ga-1)^{\frac{3}{4}}(\rho-\wi{\rho})^{3},$$
and note that
\be
 &&g'(\rho)=\ga\rho^{\ga-1}-\ga\wi{\rho}^{\ga-1}-3(\ga-1)^{\frac{3}{4}}(\rho-\wi{\rho})^2,\nonumber \\[2mm]
 &&g''(\rho)=\ga(\ga-1)\rho^{\ga-2}-6(\ga-1)^{\frac{3}{4}}(\rho-\wi{\rho}).\nonumber
\ee
Thus $g(\wi{\rho})=g'(\wi{\rho})=0$, so that
\be\label{3d-le-g} g(\rho)=g''(\zeta)(\rho-\wi{\rho})^2,\ee
where $\zeta=(1-\theta_2)\rho+\theta_2\wi{\rho}\ \  (0<\theta_2<1)$. Then
\be
g''(\zeta)&=&\ga(\ga-1)\zeta^{\ga-2}-6(\ga-1)^{\frac{3}{4}}(\zeta-\wi{\rho}) \nonumber \\[2mm]
&=& \ga(\ga-1)\left[(1-\theta_2)\rho+\theta_2\wi{\rho}\right]^{\ga-2}-6(1-\theta_2)(\ga-1)^{\frac{3}{4}}(\rho-\wi{\rho}) \nonumber \\[2mm]
&>& (\ga-1)\left[C\bar{\rho}^{\ga-2}-C(\ga-1)^{\frac{1}{12}}\right]  \nonumber \\[2mm]
&>& 0, \nonumber \ee when $\displaystyle\ga-1\le \eta.$  This
together with (\ref{3d-le-g}) implies the second inequality of
(\ref{3d-G}). We thus obtain (\ref{3d-G}). Combining
(\ref{3d-basic-pr-1}) and (\ref{3d-G}), we get
\begin{eqnarray*}
&&\displaystyle \int_{\Sigma_1} |\rho-\wi{\rho}|^3\le C(\bar{\rho})\int_{\Sigma_1} |\rho-\wi{\rho}|^{\ga-1}\le C(\bar{\rho})(\ga-1)^{\frac{1}{4}}\int_{\mathbb{R}^3} G(\rho)\le  C(\bar{\rho})(\ga-1)^{\frac{1}{4}}E_0,\\
&&\displaystyle \int_{\Sigma_2} |\rho-\wi{\rho}|^3\le (\ga-1)^{\frac{1}{4}}\int_{\mathbb{R}^3} G(\rho)\le  (\ga-1)^{\frac{1}{4}}E_0,
\end{eqnarray*}
where \be \begin{cases}
\Sigma_1=\left\{x\in \mathbb{R}^3: \ |\rho(x,t)-\wi{\rho}|>(\ga-1)^{\frac{1}{3}}\right\}, \nonumber\\[2mm]
\Sigma_2=\left\{x\in \mathbb{R}^3: \ |\rho(x,t)-\wi{\rho}|\le(\ga-1)^{\frac{1}{3}}\right\}.\nonumber \end{cases}\ee
Thus
$$
\int_{\mathbb{R}^3} |\rho-\wi{\rho}|^3 = \int_{\Sigma_1} |\rho-\wi{\rho}|^3+\int_{\Sigma_2} |\rho-\wi{\rho}|^3\le  C(\bar{\rho})(\ga-1)^{\frac{1}{4}}E_0.
$$

On the other hand, for $\displaystyle 1>\gamma-1 > \eta$, it is clear that $\displaystyle G(\rho)\geq C(\bar{\rho})(\rho-\tilde{\rho})^2$. Then
$$
\int_{\mathbb{R}^3} |\rho-\wi{\rho}|^3 \le C(\bar{\rho})\int_{\mathbb{R}^3}|\rho-\tilde{\rho}|^2\le  C(\bar{\rho})\eta^{-\frac{1}{4}}(\ga-1)^{\frac{1}{4}}E_0
\le  C(\bar{\rho})(\ga-1)^{\frac{1}{4}}E_0.
$$
This completes the prove of the Lemma \ref{3d-le-rho3}.

\begin{lemma}\label{4.2}
Under the conditions of Proposition \ref{prop 4.1}, we have
\be\label{4d-u-1}\begin{aligned}
 \displaystyle\int_0^{\sigma(T)}\int_{\mathbb{R}^3}|\nabla u|^2\le
\frac{C(\ga-1)^{\frac{1}{13}}E_0^{\frac{25}{36}}E_{12}}{\mu^\frac{12}{13}}+\frac{C\wi{\rho}E_0}{\mu},
\end{aligned}\ee
provided
$\displaystyle\frac{(\gamma-1)^\frac{1}{36}E_0^{\frac{1}{4}}}{\mu^{\frac{1}{3}}}\le
\frac{\wi{\rho}}{2 C}$.

\end{lemma}
\pf Multiplying $(\ref{full N-S+1})_2$ by $u$ and then integrating the resulting equality over
$\mathbb{R}^3$, and using integration
by parts, we have
\be\label{3d-dt rho u2}\begin{aligned}
\frac{1}{2}\frac{d}{dt}\int_{\mathbb{R}^3}\rho|u|^2+\int_{\mathbb{R}^3}\left(\mu|\nabla u|^2+(\lambda+\mu)|\mathrm{div}u|^2\right)
= \int_{\mathbb{R}^3}(P-P(\wi{\rho}))\mathrm{div}u.
\end{aligned}\ee
Now, we turn to estimate the term on the right-hand side of (\ref{3d-dt rho u2}),
\be\label{3d-p-p}
\int_{\mathbb{R}^3}(P-P(\wi{\rho}))\mathrm{div}u &\le& \|P-P(\wi{\rho})\|_{L^3}\|\nabla u\|_{L^{\frac{3}{2}}}\nonumber \\[2mm]
&\le& \|P-P(\wi{\rho})\|_{L^3}\|u\|_{L^2}^{\frac{1}{2}}\|\nabla u\|_{L^2}^{\frac{1}{2}} \nonumber \\[2mm]
&\le&
\frac{C}{\mu}\|P-P(\wi{\rho})\|_{L^3}^{\frac{4}{3}}\|u\|_{L^2}^{\frac{2}{3}}+\frac{\mu}{16}\|\nabla
u\|_{L^2}^2 \nonumber \\[2mm]
&\le&
\frac{C}{\mu}(\gamma-1)^{\frac{1}{9}}E_0^{\frac{4}{9}}\|u\|_{L^2}^{\frac{2}{3}}+\frac{\mu}{16}\|\nabla
u\|_{L^2}^2. \ee Since \be\label{3d-p-p-p}
\wi{\rho}\int_{\mathbb{R}^3}|u|^2 &\le& \int_{\mathbb{R}^3}|\rho-\wi{\rho}||u|^2 +\int_{\mathbb{R}^3}\rho|u|^2 \nonumber \\[2mm]
&\le&\left(\int_{\mathbb{R}^3} |\rho-\wi{\rho}|^3\right)^{\frac{1}{3}}\left(\int_{\mathbb{R}^3}|u|^3\right)^{\frac{2}{3}}+E_0\nonumber \\[2mm]
&\le& C(\ga-1)^{\frac{1}{12}}E_0^{\frac{1}{3}}(\mu^{-\frac{1}{2}}\|u\|_{L^2}^2
+\mu^{\frac{1}{2}}\|\nabla u\|_{L^2}^2)+E_0\nonumber \\[2mm]
&\le& \frac{C(\ga-1)^{\frac{1}{12}}E_0^{\frac{1}{3}}}{\mu^{\frac{1}{2}}}\|u\|_{L^2}^2
+ C(\ga-1)^{\frac{1}{12}}E_0^{\frac{1}{3}}\mu^{\frac{1}{2}}\|\nabla u\|_{L^2}^2+E_0\nonumber\\[2mm]
&\le& \frac{\wi{\rho}}{2}\|u\|_{L^2}^2+
C(\ga-1)^{\frac{1}{12}}E_0^{\frac{1}{3}}\mu^{\frac{1}{2}}\|\nabla
u\|_{L^2}^2+E_0,\ee provided
$\displaystyle\frac{(\gamma-1)^\frac{1}{36}E_0^{\frac{1}{4}}}{\mu^{\frac{1}{3}}}\le
\frac{\tilde{\rho}}{2 C} ,\ \
\frac{(\gamma-1)^\frac{1}{18}E_0^{\frac{1}{12}}}{\mu^{\frac{1}{6}}}\le
1$. Then \be\label{3d-p-u} \int_{\mathbb{R}^3}|u|^2 &\le&
\frac{C(\ga-1)^{\frac{1}{12}}E_0^{\frac{1}{3}}\mu^{\frac{1}{2}}\|\nabla
u\|_{L^2}^2}{\wi{\rho}}
+\frac{C E_0}{\wi{\rho}}\nonumber\\[2mm]
&\le&
C(\ga-1)^{\frac{1}{18}}E_0^{\frac{1}{12}}\mu^{\frac{5}{6}}\|\nabla
u\|_{L^2}^2
+\frac{C\mu^\frac{1}{3}E_0^\frac{3}{4}}{(\gamma-1)^\frac{1}{36}}.
\ee Combining (\ref{3d-p-p}) and (\ref{3d-p-u}), then using Young
inequality, we get \be\label{3d-p1}
\int_{\mathbb{R}^3}\left(P-P(\wi{\rho})\right)\mathrm{div}u &\le&
\frac{C}{\mu}(\ga-1)^{\frac{1}{9}}E_0^{\frac{4}{9}}\left((\ga-1)^{\frac{1}{54}}E_0^{\frac{1}{36}}\mu^{\frac{5}{18}}\|\nabla
u\|_{L^2}^{\frac{2}{3}}+\frac{\mu^\frac{1}{9}E_0^{\frac{1}{4}}}{(\gamma-1)^\frac{1}{108}}\right)
+\frac{\mu}{16}\|\nabla u\|_{L^2}^2\nonumber \\[2mm]
&\le& \frac{C(\ga-1)^{\frac{7}{54}}E_0^{\frac{17}{36}}}{\mu^{\frac{13}{18}}}\|\nabla u\|_{L^2}^{\frac{2}{3}}
+\frac{(\ga-1)^{\frac{11}{108}}E_0^{\frac{25}{36}}}{\mu^\frac{8}{9}}
+\frac{\mu}{16}\|\nabla u\|_{L^2}^2\nonumber \\[2mm]
&\le& \frac{C(\ga-1)^{\frac{7}{36}}E_0^{\frac{17}{24}}}{\mu^{\frac{19}{12}}}
+\frac{C(\ga-1)^{\frac{11}{108}}E_0^{\frac{25}{36}}}{\mu^\frac{8}{9}}
+\frac{\mu}{8}\|\nabla u\|_{L^2}^2.
\ee
Substituting (\ref{3d-p1}) into (\ref{3d-dt rho u2}) and integrating the resulting inequality over $[0,\sigma(T)]$, we get
\be\label{3d-basic-pr-u2}
&&\displaystyle\sup\limits_{0\le t\le \sigma(T)}\frac{1}{2}\int_{\mathbb{R}^3}
\displaystyle\rho|u|^2+\int_0^{\sigma(T)}\int_{\mathbb{R}^3}\left(\frac{\mu}{2}|\nabla u|^2+(\lambda+\mu)|\mathrm{div}u|^2\right)\nonumber\\[2mm]
&\le &\displaystyle\frac{1}{2}\int_{\mathbb{R}^3}(\rho_0-\wi{\rho})|u_0|^2+\frac{1}{2}\int_{\mathbb{R}^3}\wi{\rho}|u_0|^2 +\frac{C(\ga-1)^{\frac{7}{36}}E_0^{\frac{17}{24}}}{\mu^{\frac{19}{12}}}
+\frac{C(\ga-1)^{\frac{11}{108}}E_0^{\frac{25}{36}}}{\mu^\frac{8}{9}}\nonumber\\[2mm]
&\le &\displaystyle C\|\rho_0-\wi{\rho}\|_{L^3}\|u_0\|_{L^3}^2+C\wi{\rho}\|u_0\|^2
+\frac{C(\ga-1)^{\frac{7}{36}}E_0^{\frac{17}{24}}}{\mu^{\frac{19}{12}}}
+\frac{C(\ga-1)^{\frac{11}{108}}E_0^{\frac{25}{36}}}{\mu^\frac{8}{9}}\nonumber\\[2mm]
&\le&\displaystyle C(\gamma-1)^\frac{1}{12}E_0^{\frac{5}{6}}+C\wi{\rho}E_0
+\frac{C(\ga-1)^{\frac{7}{36}}E_0^{\frac{17}{24}}}{\mu^{\frac{19}{12}}}
+\frac{C(\ga-1)^{\frac{11}{108}}E_0^{\frac{25}{36}}}{\mu^\frac{8}{9}}\nonumber\\[2mm]
&\le&\displaystyle
C(\ga-1)^{\frac{1}{13}}E_0^{\frac{25}{36}}\bigg((\ga-1)^{\frac{1}{156}}E_0^\frac{5}{36}
+\frac{(\ga-1)^{\frac{7}{36}-\frac{1}{13}}E_0^{\frac{1}{72}}}{\mu^{\frac{19}{12}}}
+\frac{(\ga-1)^{\frac{11}{108}-\frac{1}{13}}}{\mu^\frac{8}{9}}
\bigg)+C\wi{\rho}E_0. \ee
Claim: \be
 \displaystyle\int_0^{\sigma(T)}\int_{\mathbb{R}^3}|\nabla u|^2\le
\frac{C(\ga-1)^{\frac{1}{13}}E_0^{\frac{25}{36}}E_{12}}{\mu^\frac{12}{13}}
+\frac{C\wi{\rho}E_0}{\mu}, \ee where \be \displaystyle
E_{12}
&=&\displaystyle1+\frac{1}{\mu}
+\frac{(\ga-1)^\frac{35}{1404}}{\mu^\frac{113}{117}}. \ee
 In fact, \eqref{3d-basic-pr-u2} implies

\be
&&\displaystyle\int_0^{\sigma(T)}\int_{\mathbb{R}^3}|\nabla u|^2\nonumber\\[2mm]
&\le& \frac{C(\ga-1)^{\frac{1}{13}}E_0^{\frac{25}{36}}}{\mu}\left((\ga-1)^{\frac{1}{156}}E_0^\frac{5}{36}
+\frac{(\ga-1)^{\frac{7}{36}-\frac{1}{13}}E_0^{\frac{1}{72}}}{\mu^{\frac{19}{12}}}
+\frac{(\ga-1)^{\frac{11}{108}-\frac{1}{13}}}{\mu^\frac{8}{9}}
\right)
+\frac{C\wi{\rho}E_0}{\mu}\nonumber\\[2mm]
&=& \frac{C(\ga-1)^{\frac{1}{13}}E_0^{\frac{25}{36}}}{\mu^\frac{12}{13}}
\left(\frac{(\ga-1)^{\frac{1}{156}}E_0^\frac{5}{36}}{\mu^\frac{1}{13}}
+\frac{(\ga-1)^{\frac{7}{36}-\frac{1}{13}}E_0^{\frac{1}{72}}}{\mu^{\frac{19}{12}+\frac{1}{13}}}
+\frac{(\ga-1)^{\frac{11}{108}-\frac{1}{13}}}{\mu^{\frac{8}{9}+\frac{1}{13}}}
\right)+\frac{C\wi{\rho}E_0}{\mu}\nonumber\\[2mm]
&=& \frac{C(\ga-1)^{\frac{1}{13}}E_0^{\frac{25}{36}}}{\mu^\frac{12}{13}}
\left\{\left(\frac{(\ga-1)E_0^{\frac{5}{36}\times{156}}}{\mu^{12}}\right)^\frac{1}{156}
+\left(\frac{(\gamma-1)E_0^{\frac{1}{72}\times\frac{468}{55}}}{\mu^{\frac{103}{156}\times\frac{468}{55}}}\right)
^{\frac{468}{55}}\frac{1}{\mu}\right.\nonumber\\[2mm]
&&\left.+\frac{(\ga-1)^{\frac{11}{108}-\frac{1}{13}}}{\mu^{\frac{8}{9}+\frac{1}{13}}}
\right\}+\frac{C\wi{\rho}E_0}{\mu}\nonumber\\[2mm]
&\le&\frac{C(\ga-1)^{\frac{1}{13}}E_0^{\frac{25}{36}}E_{12}}{\mu^\frac{12}{13}}+\frac{C\wi{\rho}E_0}{\mu},
\ee where \eqref{assumpution rb} has been used.
\endpf

\begin{lemma}\label{4.3}
Under the conditions of Proposition \ref{prop 4.1}, we have
\be\label{4d-A1}
\displaystyle A_1(T)&\leq &\displaystyle
\sup_{0\le t\le T}\left\{\frac{4}{\mu}\int_{\mathbb{R}^3}\sigma\mathrm{div}u(P-P(\wi{\rho}))\right\}
+\frac{C\left(\gamma-1+P(\wi{\rho})\right)E_0}{\mu}\nonumber\\[3mm]
&&+\frac{C}{\mu}\int_0^{T}\int_{\mathbb{R}^3}\sigma |P-P(\wi{\rho})| |\nabla u|^2
+\frac{C(2\mu+\lambda)}{\mu}\int_0^T \sigma \|\nabla u\|_{L^3}^3
\ee
and
\be\label{4d-A2} \displaystyle A_2(T)&\leq& C A_1(\sigma(T))
+C\left(\frac{1}{\mu^2}+\frac{(\gamma-1)^2}{\mu^2}\right)\int_0^T\int_{\mathbb{R}^3}\sigma^3|P-P(\wi{\rho})|^4\nonumber\\[3mm]
&&+C\left(\frac{1}{\mu^2}+\frac{(\gamma-1)^2}{\mu^2}+\frac{(2\mu+\lambda)^2}{\mu^2}\right)\int_0^T\int_{\mathbb{R}^3}\sigma^3|\nabla u|^4
+\frac{CP(\wi{\rho})^2}{\mu^2}E_0.
\ee

\end{lemma}
\pf
The proof of Lemma \ref{4.3} is similar to  Lemma \ref{3d-le:3.2}, we just need to deal with the first terms  in
\eqref{3d-A1-eq} and \eqref{3d-A2-eq} again. Here $J_1$ and $J_2$ denote $I_1$ and $II_1$ in Lemma \ref{3d-le:3.2} respectively.
Integrating by parts gives
\be\label{4d-A1-J-1}
J_1&=&\displaystyle-\int_{\mathbb{R}^3}\sigma^m \dot{u}\cdot\nabla P\nonumber\\[3mm]
&=&\displaystyle\int_{\mathbb{R}^3}\sigma^m  \mathrm{div}u_t(P-P(\wi{\rho}))
+\int_{\mathbb{R}^3}\sigma^m  \mathrm{div}(u\cdot\nabla u)(P-P(\wi{\rho}))\nonumber\\[3mm]
&=&\displaystyle\frac{d}{dt}\left(\int_{\mathbb{R}^3}\sigma^m \mathrm{div}u(P-P(\wi{\rho}))\right)
-m\sigma^{m-1}\sigma'\int_{\mathbb{R}^3}\sigma^m \mathrm{div}u(P-P(\wi{\rho}))\nonumber\\[3mm]
&&+\int_{\mathbb{R}^3}\sigma^m\Big((\gamma-1)P(\mathrm{div}u)^2+(P-P(\wi{\rho}))\partial_iu^j\partial_ju^i
+P(\wi{\rho})(\mathrm{div}u)^2
\Big)
\ee
and
\begin{eqnarray}
\label{3d-A2-I-1}
J_2&=&\displaystyle-\int_0^T\int_{\mathbb{R}^3}\sigma^m\dot{u}^j\left[\partial_jP_t+\mathrm{div}(\partial_jP u)\right]\nonumber\\[3mm]
&=&\displaystyle\int_0^T\int_{\mathbb{R}^3}\sigma^m\partial_{j}\dot{u}^jP_t+\int_0^T\int_{\mathbb{R}^3} \sigma^m\partial_{k}\dot{u}^j(\partial_{j}Pu^k)\nonumber\\[3mm]
&=&\displaystyle(1-\gamma)\int_0^T\int_{\mathbb{R}^3}\sigma^m\mathrm{div}\dot{u}\mathrm{div}u(P-P(\wi{\rho}))
-\gamma P(\wi{\rho})\int_0^T\int_{\mathbb{R}^3}\sigma^m\mathrm{div}\dot{u}\mathrm{div}u\nonumber\\[3mm]
&&-\int_0^T\int_{\mathbb{R}^3}\sigma^m\partial_k\dot{u}^j\partial_ju^k(P-P(\wi{\rho})) \nonumber\\[3mm]
&\le& \displaystyle\frac{\mu}{4}\int_0^T\int_{\mathbb{R}^3}\sigma^m|\nabla\dot{u}|^2
+C\left(\frac{1}{\mu}+\frac{(\gamma-1)^2}{\mu}\right)\int_0^T\int_{\mathbb{R}^3}\sigma^m|P-P(\wi{\rho})|^4\nonumber\\[3mm]
&&+C\left(\frac{1}{\mu}+\frac{(\gamma-1)^2}{\mu}\right)\int_0^T\int_{\mathbb{R}^3}\sigma^m|\nabla u|^4
+\frac{CP(\wi{\rho})^2}{\mu}\int_0^T\int_{\mathbb{R}^3}\sigma^m |\nabla u|^2.
\end{eqnarray}

Then, from \eqref{3d-A1-I-5} and \eqref{3d-A2-I-4}, we have
\be\label{4d-A1-I-5}
&&\displaystyle\sup_{0\leq t\leq T}\left(\frac{\mu}{4}\sigma^m\|\nabla u\|_{L^2}^2+\frac{(\lambda+\mu)}{2}\sigma^m\|\mathrm{div}u\|_{L^2}^2\right)+\int_0^{T}\int_{\mathbb{R}^3}\sigma^m \rho|\dot{u}|^2\nonumber\\[3mm]
&\leq&\displaystyle \sup_{0\le t\le T}\left\{\int_{\mathbb{R}^3}\sigma^m\mathrm{div}u(P-P(\wi{\rho}))\right\}
-\int_0^{\sigma(T)}\int_{\mathbb{R}^3}m\sigma^{m-1}\mathrm{div}u(P-P(\wi{\rho}))\nonumber\\[3mm]
&&\displaystyle+C\int_0^T\int_{\mathbb{R}^3}\sigma^m\Big((\gamma-1)P(\mathrm{div}u)^2+(P-P(\wi{\rho}))|\nabla u|^2
+P(\wi{\rho})(\mathrm{div}u)^2
\Big)\nonumber\\[3mm]
&&\displaystyle+C(2\mu+\lambda)\int_0^{T}\sigma^m \int_{\mathbb{R}^3} |\nabla u|^3
\ee
 and

\be\label{4d-A2-I-4}
&&\displaystyle\sigma^m\int_{\mathbb{R}^3}\rho|\dot{u}|^2+\mu\int_0^T\int_{\mathbb{R}^3} \sigma^m |\nabla \dot{u}|^2+(\mu+\lambda)\int_0^T\int_{\mathbb{R}^3} \sigma^m |\mathrm{div}\dot{u}|^2\nonumber\\[3mm]
&\le &\displaystyle C\int_0^T\sigma^{m-1}\sigma'\int_{\mathbb{R}^3}\rho|\dot{u}|^2
+C\left(\frac{1}{\mu}+\frac{(\gamma-1)^2}{\mu}\right)\int_0^T\int_{\mathbb{R}^3}\sigma^m|P-P(\wi{\rho})|^4\nonumber\\[3mm]
&&+C\left(\frac{1}{\mu}+\frac{(\gamma-1)^2}{\mu}+\frac{(2\mu+\lambda)^2}{\mu}\right)\int_0^T\int_{\mathbb{R}^3}\sigma^m|\nabla u|^4\nonumber\\[3mm]
&&+\frac{CP(\wi{\rho})^2}{\mu}\int_0^T\int_{\mathbb{R}^3}\sigma^m |\nabla u|^2.
\ee
Choosing $m=1$ in \eqref{4d-A1-I-5} and $m=3$ in \eqref{4d-A2-I-4}, one gets \eqref{4d-A1} and \eqref{4d-A2}.
\endpf

\begin{lemma}\label{4d-le:3.3}Under the conditions of Proposition \ref{prop 4.1},
 it holds that
\be
&&\sup\limits_{0\le t\le
\sigma(T)}\int_{\mathbb{R}^3}|\nabla
u|^2+\int_0^{\sigma(T)}\int_{\mathbb{R}^3}\frac{\rho|\u|^2}{\mu}\le
E_{13} ,\label{4d-A-1}\\[2mm]
&&\displaystyle\sup_{0\leq t\leq
\sigma(T)}\sigma\int_{\mathbb{R}^3}\frac{\rho|\dot{u}|^2}{\mu}
+\int_0^{\sigma(T)}\int_{\mathbb{R}^3}\sigma|\nabla\dot{u}|^2 \le
E_{14}, \ee provided
$$
\displaystyle\frac{\left((\gamma-1)^{\frac{1}{36}}+\wi{\rho}^\frac{1}{6}\right)E_0^\frac{1}{4}}{\mu^\frac{1}{3}}
\leq \min\left\{\Big(4C(\bar{\rho})\Big)^{-6}, 1 \right\}\triangleq
\varepsilon_4,$$ where
\be\displaystyle E_{13}&=&C(\bar{\rho})+C(M+1)
+C\left(\frac{1}{\mu}+\gamma+P(\bar{\rho})\right)
\left(E_{12}+1\right),\\
E_{14}&=&E_{13}+C(\bar{\rho})\left(1+(\gamma-1)^2\right)
+C\left(\frac{1}{\mu}+\frac{(\gamma-1)^2}{\mu}+\frac{(2\mu+\lambda)^2}{\mu}\right)^2\frac{E^3_{13}}{\mu^4}\nonumber\\[3mm]
&&+C\left(\frac{1}{\mu^2}+\frac{(\gamma-1)^2}{\mu^2}
+\frac{(2\mu+\lambda)^2}{\mu^2}\right)\frac{E_{13}^\frac{1}{2}}{\mu^\frac{3}{2}}
+C{\wi{\rho}}^{2\gamma-1}. \ee

\end{lemma}

\pf
By an argument similar to the proof of Lemma \ref{3d-le:3.3}, one can derive from \eqref{3d-le:3.3-eq} that
\be
\label{4d-le:3.3-eq}
&&\displaystyle\frac{d}{dt}\left(\frac{\mu}{2}\|\nabla u\|_{L^2}^2+\frac{(\lambda+\mu)}{2}\|\mathrm{div}u\|_{L^2}^2-\int_{\mathbb{R}^3}\mathrm{div}u (P-P(\wi{\rho}))\right)+\int_{\mathbb{R}^3}\rho|\dot{u}|^2\nonumber\\[2mm]
&=&\displaystyle\int_{\mathbb{R}^3}\rho\dot{u}(u\cdot\nabla u)-\int_{\mathbb{R}^3}\mathrm{div}u P_t\nonumber\\[2mm]
&\le&\displaystyle C(\overline{\rho})\left(\int_{\mathbb{R}^3}\rho|\dot{u}|^2\right)^{\frac{1}{2}}\left(\int_{\mathbb{R}^3}\rho|u|^3\right)^{\frac{1}{3}}\|\nabla u\|_{L^6}
+\int_{\mathbb{R}^3}\mathrm{div}u\mathrm{div}((P-P(\wi{\rho}))u)\nonumber\\[2mm]
&&+(\gamma-1)\int_{\mathbb{R}^3}P|\mathrm{div}u|^2
+P(\wi{\rho})\int_{\mathbb{R}^3}|\mathrm{div}u|^2\nonumber\\[2mm]
&\le&\displaystyle C(\overline{\rho})\left(\int_{\mathbb{R}^3}\rho|\dot{u}|^2\right)A^{\frac{1}{3}}_3(\sigma(T))
+C(\bar{\rho})A_3^{\frac{1}{3}}(\sigma(T))(\gamma-1)^\frac{1}{12}E_0^\frac{1}{3}\nonumber\\[2mm]
&&+\frac{C}{2\mu+\lambda}\Big(\|\nabla u\|_{L^2}^2+\|P-P(\wi{\rho})\|_{L^4}^4\Big)
\displaystyle
+\frac{C(\overline{\rho})}{(2\mu+\lambda)^2}\|P-P(\wi{\rho})\|_{L^3}^2\|\nabla
u\|_{L^2}^2\nonumber\\[2mm]
&&+\frac{1}{4}\|\sqrt{\rho}
\dot{u}\|_{L^2}^2+C(\bar{\rho})(\gamma-1)\|\nabla u \|_{L^2}^2+P(\wi{\rho})\|\nabla u \|_{L^2}^2.
\ee
Integrating (\ref{4d-le:3.3-eq}) over $(0, \sigma(T))$ and using \eqref{3d-p1} give that
\be
\label{4A}
&&\displaystyle\frac{\mu}{2}\|\nabla u\|_{L^2}^2+\frac{(\lambda+\mu)}{2}\|\mathrm{div}u\|_{L^2}^2-\int_{\mathbb{R}^3}\mathrm{div}u (P-P(\wi{\rho}))+\frac{1}{2}\int_0^{\sigma(T)}\int_{\mathbb{R}^3}\rho|\dot{u}|^2\nonumber\\[2mm]
&\le&\displaystyle C(\bar{\rho})A_3^{\frac{1}{3}}(\sigma(T))(\gamma-1)^\frac{1}{12}E_0^\frac{1}{3}
+\frac{C}{(2\mu+\lambda)}\int_0^{\sigma(T)}\|P-P(\wi{\rho})\|_{L^4}^4\nonumber\\[2mm]
&&+C\left(\frac{1}{(2\mu+\lambda)}+\frac{(\gamma-1)^\frac{1}{6}E_0^{\frac{2}{3}}}{(2\mu+\lambda)^2}+(\gamma-1)+P(\wi{\rho})\right)\int_0^{\sigma(T)}\|\nabla u\|^2_{L^2}+C\mu (M+1).\nonumber
\ee
provided
$\displaystyle\frac{\left((\gamma-1)^{\frac{1}{36}}+\wi{\rho}^\frac{1}{6}\right)E_0^\frac{1}{4}}{\mu^\frac{1}{3}}
\leq \left(4C(\bar{\rho})\right)^{-6}.$

\vspace{3mm}
Then, using Lemma \ref{3d-le:3.1-1}, Lemma \ref{3d-le-rho3} and \eqref{a priori assumption}, we have
\be\label{4.23}
&&\sup\limits_{0\le t\le \sigma(T)}\int_{\mathbb{R}^3}|\nabla
u|^2+\int_0^{\sigma(T)}\int_{\mathbb{R}^3}\frac{\rho|\u|^2}{\mu}\nonumber\\[2mm]
&\le&\displaystyle{\frac{C(\bar{\rho})(\gamma-1)^\frac{1}{12}E_0^\frac{1}{3}}{\mu}}
\left(\frac{\left((\gamma-1)^{\frac{1}{36}}+\wi{\rho}^\frac{1}{6}\right)E_0^\frac{1}{4}}{\mu^\frac{1}{3}}\right)^\frac{1}{6}
+\frac{C(\gamma-1)^\frac{1}{4}E_0}{\mu^2}\nonumber\\[2mm]
&&+\frac{C}{\mu}\left(\frac{1}{(2\mu+\lambda)}+\frac{(\gamma-1)^\frac{1}{6}E_0^{\frac{2}{3}}}{(2\mu+\lambda)^2}+(\gamma-1)+P(\wi{\rho})\right)\int_0^{\sigma(T)}\|\nabla u\|^2_{L^2}+C(M+1)\nonumber\\[2mm]
&\le&\displaystyle{\frac{C(\bar{\rho})(\gamma-1)^\frac{1}{12}E_0^\frac{1}{3}}{\mu}}
+\frac{C(\gamma-1)^\frac{1}{4}E_0}{\mu^2}+C(M+1)\nonumber\\[2mm]
&&+\frac{C}{\mu}\bigg(\frac{1}{(2\mu+\lambda)}+\frac{(\gamma-1)^\frac{1}{6}E_0^{\frac{2}{3}}}{(2\mu+\lambda)^2}
+(\gamma-1)+P(\wi{\rho})\bigg)\bigg(\frac{(\ga-1)^{\frac{1}{13}}E_0^{\frac{25}{36}}E_{12}}{\mu^\frac{12}{13}}
+\frac{\wi{\rho}E_0}{\mu}\bigg). \ee Using \eqref{assumpution rb},
we get \be
&&\displaystyle{\frac{C(\bar{\rho})(\gamma-1)^\frac{1}{12}E_0^\frac{1}{3}}{\mu}}
+\frac{C(\gamma-1)^\frac{1}{4}E_0}{\mu^2}+C(M+1)\nonumber\\[2mm]
&&+\frac{C}{\mu}\left(\frac{1}{2\mu+\lambda}+\frac{(\gamma-1)^\frac{1}{6}E_0^{\frac{2}{3}}}{(2\mu+\lambda)^2}
+(\gamma-1)+P(\wi{\rho})\right)\left(\frac{(\ga-1)^{\frac{1}{13}}E_0^{\frac{25}{36}}E_{12}}{\mu^\frac{12}{13}}
+\frac{\wi{\rho}E_0}{\mu}\right)\nonumber\\[2mm]
&=& C(\bar{\rho})\left(\frac{(\gamma-1)E_0^4}{\mu^{12}}\right)^\frac{1}{12}
+C\left(\frac{(\gamma-1)E_0^4}{\mu^{8}}\right)^\frac{1}{4}
+C(M+1)\nonumber\\[2mm]
&&+C\left(\frac{1}{2\mu+\lambda}+\left(\frac{(\gamma-1)E_0^4}{(2\mu+\lambda)^{12}}\right)^\frac{1}{6}
+(\gamma-1)+P(\wi{\rho})\right)\left(\left(\frac{(\ga-1)E_0^{\frac{325}{36}}}{\mu^{12}}\right)^{\frac{1}{13}}E_{12}
+1\right)\nonumber\\[2mm]
&\le& C(\bar{\rho})+C(M+1)
+C\left(\frac{1}{\mu}+\gamma+P(\bar{\rho})\right)
\left(E_{12}+1\right)= E_{13}. \ee
\vspace{2mm}
Next taking $m=1$ in \eqref{4d-A2-I-4}, we have
\be\label{4.24}
&&\displaystyle\sigma\int_{\mathbb{R}^3}\rho|\dot{u}|^2+\mu\int_0^{\sigma(T)}\int_{\mathbb{R}^3} \sigma |\nabla \dot{u}|^2+(\mu+\lambda)\int_0^{\sigma(T)}\int_{\mathbb{R}^3} \sigma |\mathrm{div}\dot{u}|^2\nonumber\\[2mm]
&\le& \int_0^{\sigma(T)}\int_{\mathbb{R}^3}\rho|\dot{u}|^2+\left(\frac{1}{\mu}
+\frac{(\gamma-1)^2}{\mu}\right)\int_0^{\sigma(T)}\int_{\mathbb{R}^3}\sigma|P-P(\wi{\rho})|^4\nonumber\\[2mm]
&&+C\left(\frac{1}{\mu}+\frac{(\gamma-1)^2}{\mu}+\frac{(2\mu+\lambda)^2}{\mu}\right)\int_0^{\sigma(T)}\int_{\mathbb{R}^3}\sigma|\nabla u|^4+\frac{CP(\wi{\rho})^2}{\mu}\int_0^{\sigma(T)}\int_{\mathbb{R}^3}\sigma |\nabla u|^2\nonumber\\[2mm]
&\le& \int_0^{\sigma(T)}\int_{\mathbb{R}^3}\rho|\dot{u}|^2+C(\bar{\rho})\left(\frac{1}{\mu}
+\frac{(\gamma-1)^2}{\mu}\right)\|P-P(\wi{\rho})\|_{L^3}^3+\frac{CP(\wi{\rho})^2}{\mu}\int_0^{\sigma(T)}\int_{\mathbb{R}^3}\sigma|\nabla u|^2
\nonumber\\[2mm]
&&+C\left(\frac{1}{\mu}+\frac{(\gamma-1)^2}{\mu}+\frac{(2\mu+\lambda)^2}{\mu}\right)\int_0^{\sigma(T)}\sigma\|\nabla u\|_{L^2}\|\nabla u\|^3_{L^6}\nonumber\\[2mm]
&=&\sum_{i=1}^4K_i.
\ee
In fact, we just need to deal with $K_3$. Using \eqref{hp3} and  Cauchy inequality, we get
\be\label{4.25}
K_3&=&C\left(\frac{1}{\mu}+\frac{(\gamma-1)^2}{\mu}+\frac{(2\mu+\lambda)^2}{\mu}\right)\int_0^{\sigma(T)}\sigma\|\nabla u\|_{L^2}\|\nabla u\|^3_{L^6}\nonumber\\[2mm]
&\le& \left(\frac{1}{\mu}+\frac{(\gamma-1)^2}{\mu}+\frac{(2\mu+\lambda)^2}{\mu}\right)\frac{C}{(2\mu+\lambda)^3}\sup_{0\leq t\leq \sigma(T)}\|\nabla u\|_{L^2}
\int_0^{\sigma(T)}\sigma\|\sqrt{\rho}\dot{u}\|^3_{L^2}\nonumber\\[2mm]
&&+\left(\frac{1}{\mu}+\frac{(\gamma-1)^2}{\mu}+\frac{(2\mu+\lambda)^2}{\mu}\right)\frac{C}{(2\mu+\lambda)^3}\sup_{0\leq t\leq \sigma(T)}\|\nabla u\|_{L^2}
\int_0^{\sigma(T)}\sigma\|P-P(\wi{\rho})\|^3_{L^6}\nonumber\\[2mm]
&\le& \left(\frac{1}{\mu}+\frac{(\gamma-1)^2}{\mu}+\frac{(2\mu+\lambda)^2}{\mu}\right)\frac{CE^\frac{1}{2}_{9}}{(2\mu+\lambda)^3}
\sup_{0\leq t\leq \sigma(T)}\sigma^\frac{1}{2}\|\sqrt{\rho}\dot{u}\|_{L^2}\int_0^{\sigma(T)}\|\sqrt{\rho}\dot{u}\|^2_{L^2}\nonumber\\[2mm]
&&+\left(\frac{1}{\mu}+\frac{(\gamma-1)^2}{\mu}+\frac{(2\mu+\lambda)^2}{\mu}\right)\frac{CE^\frac{1}{2}_{9}}{(2\mu+\lambda)^3}
\int_0^{\sigma(T)}\sigma\|P-P(\wi{\rho})\|^3_{L^6}\nonumber\\[2mm]
&\le& \frac{1}{4}\sup_{0\leq t\leq
\sigma(T)}\sigma\|\sqrt{\rho}\dot{u}\|^2_{L^2}
+\left(\frac{1}{\mu}+\frac{(\gamma-1)^2}{\mu}+\frac{(2\mu+\lambda)^2}{\mu}\right)^2\frac{CE_{13}}{(2\mu+\lambda)^6}
\left(\int_0^{\sigma(T)}\|\sqrt{\rho}\dot{u}\|^2_{L^2}\right)^2\nonumber\\[2mm]
&&+\left(\frac{1}{\mu}+\frac{(\gamma-1)^2}{\mu}+\frac{(2\mu+\lambda)^2}{\mu}\right)\frac{CE^\frac{1}{2}_{9}(\gamma-1)^\frac{1}{8}
E_0^\frac{1}{2}}{(2\mu+\lambda)^3}
\nonumber\\[2mm]
&\le& \frac{1}{4}\sup_{0\leq t\leq \sigma(T)}\sigma\|\sqrt{\rho}\dot{u}\|^2_{L^2}
+\left(\frac{1}{\mu}+\frac{(\gamma-1)^2}{\mu}+\frac{(2\mu+\lambda)^2}{\mu}\right)^2\frac{CE^3_{9}\mu^2}{(2\mu+\lambda)^6}\nonumber\\[2mm]
&&+\left(\frac{1}{\mu}+\frac{(\gamma-1)^2}{\mu}+\frac{(2\mu+\lambda)^2}{\mu}\right)
\frac{CE^\frac{1}{2}_{9}(\gamma-1)^\frac{1}{8}E_0^\frac{1}{2}}{(2\mu+\lambda)^3}.
\ee
Substituting \eqref{4.25} into \eqref{4.24}, and using \eqref{4.23}, we have
\be
&&\displaystyle\sup_{0\leq t\leq \sigma(T)}\sigma\int_{\mathbb{R}^3}\frac{\rho|\dot{u}|^2}{\mu}
+\int_0^{\sigma(T)}\int_{\mathbb{R}^3}\sigma|\nabla\dot{u}|^2\nonumber\\[2mm]
&\le& E_{13}+C(\bar{\rho})\left(\frac{1}{\mu^2}
+\frac{(\gamma-1)^2}{\mu^2}\right)(\gamma-1)^\frac{1}{4}E_0+\left(\frac{1}{\mu}+\frac{(\gamma-1)^2}{\mu}
+\frac{(2\mu+\lambda)^2}{\mu}\right)^2\frac{CE^3_{9}\mu}{(2\mu+\lambda)^6}\nonumber\\[2mm]
&&+\left(\frac{1}{\mu^2}+\frac{(\gamma-1)^2}{\mu^2}+\frac{(2\mu+\lambda)^2}{\mu^2}\right)
\frac{CE^\frac{1}{2}_{9}(\gamma-1)^\frac{1}{8}E_0^\frac{1}{2}}{(2\mu+\lambda)^3}+\frac{CP(\wi{\rho})^2E_0}{\mu^3}\nonumber\\[2mm]
&\le&
E_{13}+C(\bar{\rho})\left(1+(\gamma-1)^2\right)\left(\frac{(\gamma-1)E_0^4}{\mu^8}\right)^\frac{1}{4}
+\left(\frac{1}{\mu}+\frac{(\gamma-1)^2}{\mu}+\frac{(2\mu+\lambda)^2}{\mu}\right)^2\frac{CE^3_{9}}{\mu^5}\nonumber\\[2mm]
&&+\left(\frac{1}{\mu^2}+\frac{(\gamma-1)^2}{\mu^2}
+\frac{(2\mu+\lambda)^2}{\mu^2}\right)\frac{CE_{13}^\frac{1}{2}}{\mu^\frac{3}{2}}\left(\frac{(\gamma-1)E_0^4}{\mu^{12}}\right)^\frac{1}{8}
+C{\wi{\rho}}^{2\gamma-1}\frac{\wi{\rho}E_0}{\mu^3}\nonumber\\[2mm]
&\le& E_{13}+C(\bar{\rho})\left(1+(\gamma-1)^2\right)
+C\left(\frac{1}{\mu}+\frac{(\gamma-1)^2}{\mu}+\frac{(2\mu+\lambda)^2}{\mu}\right)^2\frac{E^3_{9}}{\mu^5}\nonumber\\[2mm]
&&+C\left(\frac{1}{\mu^2}+\frac{(\gamma-1)^2}{\mu^2}
+\frac{(2\mu+\lambda)^2}{\mu^2}\right)\frac{E_{13}^\frac{1}{2}}{\mu^\frac{3}{2}}
+C{\wi{\rho}}^{2\gamma-1}=E_{14}, \ee where \eqref{assumpution rb}
has been used.
\endpf

Next, we will close the $a \ priori$ assumption on $A_3(\sigma(T))$.

\begin{lemma}\label{4d-A3}
Under the conditions of Proposition \ref{prop 4.1}, it holds that
\be A_3(\sigma(T))\leq
\left\{\frac{\left((\gamma-1)^{\frac{1}{36}}+\wi{\rho}^\frac{1}{6}\right)E_0^\frac{1}{4}}{\mu^\frac{1}{3}}\right\}^\frac{1}{2},
\ee provided \be
\frac{\left((\gamma-1)^{\frac{1}{36}}+\wi{\rho}^\frac{1}{6}\right)E_0^\frac{1}{4}}{\mu^\frac{1}{3}}
&\le &\min\Bigg\{\left(C(E_{15}E_{17}+E_{16})\right)^{-4},\
\varepsilon_4 \Bigg\} \triangleq  \varepsilon_5. \ee

\end{lemma}
\pf
Multiplying \eqref{full N-S+1}$_2$ by $3|u|u$ and integrating the resulting equation over $\mathbb{R}^3$,
using Lemma \ref{GN}, Lemma \ref{2.2} and H\"older inequality, we obtain
\be\label{4A3}
&&\displaystyle\frac{d}{dt}\int_{\mathbb{R}^3}\rho|u|^3\nonumber\\[2mm]
&\le& C\mu \int_{\mathbb{R}^3} |u||\nabla u|^2+C\int_{\mathbb{R}^3} |P-P(\wi{\rho})||u||\nabla u|\nonumber\\[2mm]
&\le& C\mu\|u\|_{L^6}\|\nabla u\|_{L^2}\|\nabla u\|_{L^3}+C\|P-P(\wi{\rho})\|_{L^3}\|u\|_{L^6}\|\nabla u\|_{L^2}\nonumber\\[2mm]
&\le& C\mu^\frac{1}{2}\|\nabla u\|^\frac{5}{2}_{L^2}\Big(\|\sqrt{\rho}\dot{u}\|_{L^2}^\frac{1}{2}
+\|P-P(\wi{\rho})\|_{L^6}^\frac{1}{2}\Big)
+C\|P-P(\wi{\rho})\|_{L^3}\|\nabla u\|^2_{L^2}.
\ee
Integrating \eqref{4A3} over $(0,\sigma(T))$ and using H\"older inequality, one gets
\be\label{4A3-1}
&&\sup_{0\le t \le \sigma(T)}\int_{\mathbb{R}^3}\rho|u|^3\nonumber\\[2mm]
&\le& C\mu^\frac{1}{2}\int_0^{\sigma(T)}\|\nabla u\|^\frac{5}{2}_{L^2}\|\sqrt{\rho}\dot{u}\|_{L^2}^\frac{1}{2}
+C\mu^\frac{1}{2}\int_0^{\sigma(T)}\|\nabla u\|^\frac{5}{2}_{L^2}\|P-P(\wi{\rho})\|_{L^6}^\frac{1}{2} \nonumber\\[2mm]
&&+C\int_0^{\sigma(T)}\|P-P(\wi{\rho})\|_{L^3}\|\nabla u\|^2_{L^2}+\sup_{0\le t \le \sigma(T)}\int_{\mathbb{R}^3}\rho_0|u_0|^3\nonumber\\[2mm]
&\le& C\mu^\frac{1}{2}\left(\int_0^{\sigma(T)}\|\nabla u\|_{L^2}^4\right)^{\frac{1}{2}}
\left(\int_0^{\sigma(T)}\|\sqrt{\rho}\dot{u}\|_{L^2}^2\right)^{\frac{1}{4}}
\left(\int_0^{\sigma(T)}\|\nabla u\|_{L^2}^2\right)^{\frac{1}{4}}\nonumber\\[2mm]
&&+C\mu^\frac{1}{2}\sup_{0\le t \le \sigma(T)}\|P-P(\wi{\rho})\|_{L^6}^\frac{1}{2}
\sup_{0\le t \le \sigma(T)}\|\nabla u\|_{L^2}^{\frac{1}{2}}
\int_0^{\sigma(T)}\|\nabla u\|_{L^2}^2
\nonumber\\[2mm]
&&+C\sup_{0\le t \le \sigma(T)}\|P-P(\wi{\rho})\|_{L^3}\int_0^{\sigma(T)}\|\nabla u\|_{L^2}^2
+\int_{\mathbb{R}^3}|\rho_0-\widetilde{\rho}||u_0|^3+\int_{\mathbb{R}^3}\widetilde{\rho}|u_0|^3\nonumber\\[2mm]
&\le& C\mu^\frac{1}{2}\sup_{0\leq t\leq \sigma(T)}\|\nabla u\|_{L^2}
\left(\int_0^{\sigma(T)}\|\nabla u\|_{L^2}^2\right)^{\frac{3}{4}}
\left(\int_0^{\sigma(T)}\|\sqrt{\rho}\dot{u}\|_{L^2}^2\right)^{\frac{1}{4}}
\nonumber\\[2mm]
&&+C\mu^\frac{1}{2}\sup_{0\le t \le \sigma(T)}\|P-P(\wi{\rho})\|_{L^6}^\frac{1}{2}
\sup_{0\le t \le \sigma(T)}\|\nabla u\|_{L^2}^{\frac{1}{2}}
\int_0^{\sigma(T)}\|\nabla u\|_{L^2}^2\nonumber\\[2mm]
&&+C\sup_{0\le t \le \sigma(T)}\|P-P(\wi{\rho})\|_{L^3}\int_0^{\sigma(T)}\|\nabla u\|_{L^2}^2
+C\|\rho_0-\widetilde{\rho}\|_{L^3}\|u_0\|_{L^6}\|u_0\|^2_{L^4}\nonumber\\[2mm]
&&+C\widetilde{\rho}\|u_0\|^3_{L^3}.
\ee
By Lemma \ref{3d-le:3.1-1}, Lemma \ref{3d-le-rho3}, Lemma \ref{4.2} and Lemma \ref{4d-le:3.3},
we obtain
\be\label{4A3-11}
&&\sup_{0\le t \le \sigma(T)}\int_{\mathbb{R}^3}\frac{\rho|u|^3}{\mu^3}\nonumber\\[3mm]
&\le& \frac{C}{\mu^2}\left(\int_0^{\sigma(T)}\|\nabla
u\|_{L^2}^2\right)^{\frac{3}{4}}
\Bigg(\frac{E_{13}^\frac{3}{4}}{\mu^\frac{1}{4}}
+\frac{E_{13}^\frac{1}{4}(\gamma-1)^\frac{1}{48}E_0^{\frac{1}{12}}}{\mu^{\frac{1}{2}}}
\frac{E_0^{\frac{1}{4}}}{\mu^{\frac{1}{4}}}+\frac{(\gamma-1)^\frac{1}{12}E_0^\frac{1}{3}}{\mu} \frac{E_0^\frac{1}{4}}{\mu^\frac{1}{4}}\Bigg)\nonumber\\[3mm]
&&+\frac{C}{\mu^3}\|\rho_0-\widetilde{\rho}\|_{L^3}\|\nabla u_0\|^\frac{5}{2}_{L^2}\|u_0\|^\frac{1}{2}_{L^2}
+\frac{C\widetilde{\rho}}{\mu^3}\|u_0\|^\frac{3}{2}_{L^2}\|\nabla u_0\|^\frac{3}{2}_{L^2}
\nonumber\\[3mm]
&\le&
C\left(\frac{(\gamma-1)^\frac{1}{13}E_0^\frac{25}{36}E_{12}}{\mu^\frac{12}{13}}
+\frac{\wi{\rho}E_0}{\mu}\right)^\frac{3}{4}E_{15}
+\frac{C(\gamma-1)^{\frac{1}{12}}E_0^\frac{7}{12}M^\frac{5}{4}}{\mu^3}
+\frac{C\widetilde{\rho}E_0^\frac{3}{4}M^\frac{3}{4}}{\mu^3}\nonumber\\[3mm]
&\le &
C\left(\frac{(\gamma-1)^\frac{1}{13}E_0^\frac{25}{36}E_{12}}{\mu^\frac{12}{13}}
+\frac{\wi{\rho}E_0}{\mu}\right)^\frac{3}{4}E_{15}
+C\left(\frac{(\gamma-1)^\frac{1}{36}E_0^\frac{1}{4}}{\mu^\frac{1}{3}}
+\frac{\widetilde{\rho}^\frac{1}{6}E_0^\frac{1}{4}}{\mu^\frac{1}{3}}\right)^\frac{3}{4}E_{16}\nonumber\\[3mm]
&\le& C\left(\frac{(\gamma-1)^\frac{1}{36}E_0^\frac{1}{4}}{\mu^\frac{1}{3}}
+\frac{\widetilde{\rho}^\frac{1}{6}E_0^\frac{1}{4}}{\mu^\frac{1}{3}}\right)^\frac{3}{4}
(E_{15}E_{17}+E_{16})\nonumber\\[3mm]
&\le&
C\left(\frac{(\gamma-1)^\frac{1}{36}E_0^\frac{1}{4}}{\mu^\frac{1}{3}}
+\frac{\wi{\rho}^\frac{1}{6}E_0^\frac{1}{4}}{\mu^\frac{1}{3}}\right)^\frac{1}{2},
\ee provided $\displaystyle
\frac{(\gamma-1)^\frac{1}{36}E_0^\frac{1}{4}}{\mu^\frac{1}{3}}
+\frac{\wi{\rho}^\frac{1}{6}E_0^\frac{1}{4}}{\mu^\frac{1}{3}}\le
\left(C(E_{15}E_{17}+E_{16})\right)^{-4}$, where
\be\label{e15-17}
&&\displaystyle E_{15}=\frac{E_{13}^{\frac{3}{4}}}{\mu^\frac{9}{4}}+\frac{E_{13}^{\frac{1}{4}}}{\mu^\frac{5}{2}}+\frac{1}{\mu^\frac{9}{4}},\ \ E_{16}=\frac{M^\frac{5}{4}}{\mu^2}+M^{\frac{3}{4}},\ \
E_{17}=\left(E_{12}+\frac{\widetilde{\rho}^\frac{1}{3}}{\mu}\right)^\frac{3}{4}.
\ee

In what follows, we give the details of the calculations of
$E_{15}$-$E_{17}$. \eqref{assumpution rb} gives
\be
E_{15}&=&\displaystyle\frac{1}{\mu^2}\left(\frac{E_{13}^\frac{3}{4}}{\mu^\frac{1}{4}}
+\frac{E_{13}^\frac{1}{4}(\gamma-1)^\frac{1}{48}E_0^{\frac{1}{12}}}{\mu^{\frac{1}{2}}}
\frac{E_0^{\frac{1}{4}}}{\mu^{\frac{1}{4}}}+\frac{(\gamma-1)^\frac{1}{12}E_0^\frac{1}{3}}{\mu} \frac{E_0^\frac{1}{4}}{\mu^\frac{1}{4}}\right)\nonumber\\[2mm]
&=&
\displaystyle\frac{1}{\mu^2}\left(\frac{E_{13}^\frac{3}{4}}{\mu^\frac{1}{4}}
+\frac{E_{13}^\frac{1}{4}}{\mu^\frac{1}{2}}\left(\frac{(\gamma-1)
E_0^{16}}{\mu^{12}}\right)^\frac{1}{48}
+\frac{1}{\mu^\frac{1}{4}}\left(\frac{(\gamma-1)E_0^7}{\mu^{12}}\right)^\frac{1}{12}\right) \nonumber\\[2mm]
&=&
\displaystyle\frac{E_{13}^{\frac{3}{4}}}{\mu^\frac{9}{4}}+\frac{E_{13}^{\frac{1}{4}}}{\mu^\frac{5}{2}}+\frac{1}{\mu^\frac{9}{4}},\label{s4 E11}\\[2mm]
E_{16}&=&\frac{(\gamma-1)E_0^{\frac{19}{48}}M^{\frac{5}{4}}}{\mu^{\frac{11}{4}}}
+\frac{\wi{\rho}^{\frac{7}{8}}E_0^{\frac{9}{16}}M^{\frac{3}{4}}}{\mu^{\frac{11}{4}}}\nonumber\\[2mm]
&=&\left(\frac{(\gamma-1)E_0^\frac{19}{3}}{\mu^{12}}\right)^\frac{1}{16}\frac{M^\frac{5}{4}}{\mu^2}
+M^{\frac{3}{4}}\left(\frac{\wi{\rho}E_0^{\frac{9}{14}}}{\mu^{\frac{22}{7}}}\right)^{\frac{7}{8}}\nonumber\\[2mm]
&=& \frac{M^\frac{5}{4}}{\mu^2}+M^{\frac{3}{4}} \ee and \be
&&C\left(\frac{(\gamma-1)^\frac{1}{13}E_0^\frac{25}{36}E_{12}}{\mu^\frac{12}{13}}
+\frac{\wi{\rho}E_0}{\mu}\right)^\frac{3}{4}\nonumber\\[2mm]
&=&C\left(\frac{(\gamma-1)^\frac{1}{36}E_0^\frac{1}{4}}{\mu^\frac{1}{3}}
\frac{(\gamma-1)^\frac{23}{468}E_0^\frac{4}{9}E_{12}}{\mu^\frac{23}{39}}
+\left(\frac{\widetilde{\rho}^\frac{1}{6}E_0^\frac{1}{4}}{\mu^\frac{1}{3}}\right)^4\frac{\widetilde{\rho}^\frac{1}{3}}{\mu}
\right)^\frac{3}{4}
\nonumber\\[2mm]
&=&C\left(\frac{(\gamma-1)^\frac{1}{36}E_0^\frac{1}{4}}{\mu^\frac{1}{3}}
\left(\frac{(\gamma-1)E_0^{\frac{4}{9}\times\frac{468}{23}}}{\mu^{12}}\right)^\frac{23}{468}E_{12}
+\left(\frac{\widetilde{\rho}^\frac{1}{6}E_0^\frac{1}{4}}{\mu^\frac{1}{3}}\right)^4\frac{\widetilde{\rho}^\frac{1}{3}}{\mu}
\right)^\frac{3}{4}\nonumber\\[2mm]
&\le&C\left(\frac{(\gamma-1)^\frac{1}{36}E_0^\frac{1}{4}}{\mu^\frac{1}{3}}
+\frac{\widetilde{\rho}^\frac{1}{6}E_0^\frac{1}{4}}{\mu^\frac{1}{3}}\right)^{\frac{3}{4}}
\left(E_{12}+\frac{\widetilde{\rho}^\frac{1}{3}}{\mu}\right)^\frac{3}{4}\nonumber\\[2mm]
&=&C\left(\frac{(\gamma-1)^\frac{1}{36}E_0^\frac{1}{4}}{\mu^\frac{1}{3}}
+\frac{\widetilde{\rho}^\frac{1}{6}E_0^\frac{1}{4}}{\mu^\frac{1}{3}}\right)^{\frac{3}{4}}E_{17}.\label{s 4 E12}
\ee
\endpf

\begin{lemma}\label{3d-le:4.7}Under the conditions of Proposition \ref{prop 4.1}, it holds that
\be
&\displaystyle A_1(T)&\le
\left(\frac{\left((\gamma-1)^{\frac{1}{36}}+\wi{\rho}^\frac{1}{6}\right)E_0^\frac{1}{4}}{\mu^\frac{1}{3}}\right)^\frac{3}{8},\label{3d-A1-le} \\[2mm]
&\displaystyle A_2(T)&\le
\frac{\left((\gamma-1)^{\frac{1}{36}}+\wi{\rho}^\frac{1}{6}\right)E_0^\frac{1}{4}}{\mu^\frac{1}{3}},\label{3d-A2-le}
\ee provided \be
\frac{\left((\gamma-1)^{\frac{1}{36}}+\wi{\rho}^\frac{1}{6}\right)E_0^\frac{1}{4}}{\mu^\frac{1}{3}}\le
\varepsilon_6,\nonumber \ee where \be \varepsilon_6= \min\left\{
\Big(C(E_{18}+E_{19}+E_{20})\Big)^{-17},
\Big(C(E_{18}+E_{19}+E_{21})\Big)^{-8},
\varepsilon_5\right\}\nonumber \ee and $ E_{18}$-$ E_{21}$ are given by
\eqref{s 4E14}, \eqref{s 4E15}, \eqref{s 4E16} and \eqref{s 4E17}
respectively.

\end{lemma}
\pf First, we will prove \eqref{3d-A2-le}. Recalling \eqref{4d-A2}, we have
\be\label{4d-A2-c} \displaystyle A_2(T)&\leq& C A_1(\sigma(T))
+C\left(\frac{1}{\mu^2}+\frac{(\gamma-1)^2}{\mu^2}\right)\int_0^T\int_{\mathbb{R}^3}\sigma^3|P-P(\wi{\rho})|^4\nonumber\\[2mm]
&&+C\left(\frac{1}{\mu^2}+\frac{(\gamma-1)^2}{\mu^2}+\frac{(2\mu+\lambda)^2}{\mu^2}\right)\int_0^T\int_{\mathbb{R}^3}\sigma^3|\nabla u|^4
+\frac{CP(\wi{\rho})^2}{\mu^2}E_0.
\ee
Now, we turn to estimate the term $\displaystyle \int_0^T\int_{\mathbb{R}^3}\sigma^3|\nabla u|^4$.  Due to \eqref{hp6},
\be\label{4d-A21}
&&\displaystyle\int_0^T\int_{\mathbb{R}^3}\sigma^3|\nabla u|^4\nonumber\\[2mm]
&\le&\displaystyle C\left(\frac{1}{(2\mu+\lambda)^2}+\frac{1}{\mu^2}\right)\int_0^T\sigma^3\|\nabla u\|_{L^3}^2\|\sqrt{\rho}\dot{u}\|_{L^2}^2
+\frac{C}{(2\mu+\lambda)^4}\int_0^T\sigma^3\|\sqrt{\rho}\dot{u}\||_{L^2}^2\|P-P(\wi{\rho})\|_{L^3}^2\nonumber\\[2mm]
&&+\frac{C}{(2\mu+\lambda)^4}\int_0^T\sigma^3\|P-P(\wi{\rho})\|_{L^4}^4
\displaystyle=\sum_{i=1}^3L_{i}.
\ee
By using H\"older inequality, Young inequality and \eqref{a priori assumption}, $L_1$ can be estimated as follows,
\be\label{4d-A21-1}
L_1&=&C\displaystyle\left(\frac{1}{(2\mu+\lambda)^2}
+\frac{1}{\mu^2}\right)\int_0^T\sigma^3\|\nabla u\|_{L^3}^2\|\sqrt{\rho}\dot{u}\||_{L^2}^2\nonumber\\[2mm]
&\le&C\displaystyle\left(\frac{1}{(2\mu+\lambda)^2}
+\frac{1}{\mu^2}\right)\int_0^T\sigma^3\|\nabla u\|_{L^2}^\frac{2}{3}\|\nabla u\|_{L^4}^\frac{4}{3}\|\sqrt{\rho}\dot{u}\||_{L^2}^2\nonumber\\[2mm]
&\le&\displaystyle
\frac{C}{\mu^2}\left(\int_0^T\sigma^3\|\nabla u\|_{L^4}^4\right)^\frac{1}{3}
\left(\int_0^T\sigma^3\|\nabla u\|_{L^2}\|\sqrt{\rho}\dot{u}\||_{L^2}^3\right)^\frac{2}{3}\nonumber\\[2mm]
&\le&\displaystyle
\frac{1}{4}\int_0^T\sigma^3\|\nabla u\|_{L^4}^4
+\frac{C}{\mu^3}\int_0^T\sigma^3\|\nabla u\|_{L^2}\|\sqrt{\rho}\dot{u}\||_{L^2}^3\nonumber\\[2mm]
&\le&\displaystyle
\frac{1}{4}\int_0^T\sigma^3\|\nabla u\|_{L^4}^4
+\frac{C}{\mu^3}\sup_{0 \le t \le T}\left(\sigma^\frac{1}{2}\|\nabla u\|_{L^2}\right)
\sup_{0 \le t \le T}\left(\sigma^\frac{3}{2}\|\sqrt{\rho} \dot{u}\|_{L^2}\right)
\int_0^T\sigma\|\sqrt{\rho} \dot{u}\|_{L^2}^2\nonumber\\[2mm]
&\le&\displaystyle
\frac{1}{4}\int_0^T\sigma^3\|\nabla u\|_{L^4}^4
+\frac{C}{\mu^3}A_1^{\frac{3}{2}}(T)A_2^\frac{1}{2}(T).
\ee
It follows from Lemma \ref{3d-le-rho3} and \eqref{a priori assumption} that
\be\label{4d-A21-2}
L_2&=&\displaystyle\frac{C}{(2\mu+\lambda)^4}\int_0^T\sigma^3\|\sqrt{\rho}\dot{u}\||_{L^2}^2\|P-P(\wi{\rho})\|_{L^3}^2\nonumber\\[2mm]
&\le& \displaystyle\frac{C}{(2\mu+\lambda)^4}\sup_{0 \le t \le T}\left(\|P-P(\wi{\rho})\|_{L^3}^2\right)\int_0^T\sigma^3\|\sqrt{\rho}\dot{u}\|_{L^2}^2\nonumber\\[2mm]
&\le & \frac{C(\gamma-1)^\frac{1}{6}E_0^{\frac{2}{3}}A_1(T)}{(2\mu+\lambda)^4}.
\ee

Following a process similar to \cite{Huang-Li-Xin}, we focus on estimating the term $\displaystyle\int_0^T\|P-P(\wi{\rho}\|^4_{L^4}$.
One deduces from $(\ref{full N-S+1})_1$ that $P-P(\wi{\rho})$ satisfies
\be\label{4P}
(P-P(\wi{\rho}))_t+u\cdot\nabla(P-P(\wi{\rho}))+\gamma(P-P(\wi{\rho}))\di u+\gamma P(\wi{\rho})\di u=0.
\ee
Multiplying \eqref{4P} by $3\sigma^3(P-P(\wi{\rho}))^2$ and integrating the resulting equality over $\mathbb{R}^3\times[0,T]$, using
$\displaystyle\di u = \displaystyle\frac{1}{2\mu+\lambda}\left(G+P-P(\wi{\rho})\right)$, we get
\be\label{4dP}
&&\frac{3\gamma-1}{2\mu+\lambda}\int_0^T\int_{\mathbb{R}^3}\sigma^3|P-P(\wi{\rho})|^4\nonumber\\[2mm]
&=&\int_{\mathbb{R}^3}\sigma^3(P-P(\wi{\rho}))^3
-\frac{(3\gamma-1)}{2\mu+\lambda}\int_0^T\int_{\mathbb{R}^3}\sigma^3(P-P(\wi{\rho}))^3G
+3\int_0^{\sigma(T)}\int_{\mathbb{R}^3}\sigma^2(P-P(\wi{\rho}))^3\nonumber\\[2mm]
&&-3\gamma P(\wi{\rho})\int_0^T\int_{\mathbb{R}^3}\sigma^3(P-P(\wi{\rho}))^2\di u
\nonumber\\[2mm]
&\le& C \sup_{0 \le t \le T}(\|P-P(\wi{\rho})\|_{L^3}^3)
+C\int_0^{\sigma(T)}\int_{\mathbb{R}^3}\sigma^2(P-P(\wi{\rho}))^3
+\frac{C(3\gamma-1)}{2\mu+\lambda}\int_0^T\int_{\mathbb{R}^3}\sigma^3|G|^4\nonumber\\[2mm]
&&+\frac{3\gamma-1}{2(2\mu+\lambda)}\int_0^T\int_{\mathbb{R}^3}\sigma^3|P-P(\wi{\rho})|^4
+\frac{CP(\wi{\rho})^2(2\mu+\lambda)}{3\gamma-1}\int_0^T\int_{\mathbb{R}^3}\sigma^3|\nabla u|^2
\nonumber\\[2mm]
&\le&C(\gamma-1)^\frac{1}{4}E_0+\frac{C}{2\mu+\lambda}\int_0^T\int_{\mathbb{R}^3}\sigma^3|G|^4
+\frac{3\gamma-1}{2(2\mu+\lambda)}\int_0^T\int_{\mathbb{R}^3}\sigma^3|P-P(\wi{\rho})|^4\nonumber\\[2mm]
&&
+\frac{C(2\mu+\lambda)P^2(\wi{\rho})E_0}{\mu}.
\ee
It follows from \eqref{hp5} that
\be\label{4d-A21-3}
L_3&=&\displaystyle\frac{C(\gamma-1)^\frac{1}{4}E_0}{(2\mu+\lambda)^3}
+\frac{C}{(2\mu+\lambda)^2}\int_0^T\sigma^3\|\sqrt{\rho}\dot{u}\|_{L^2}^2\|\nabla u\|_{L^3}^2\nonumber\\[2mm]
&&+\frac{C}{(2\mu+\lambda)^4}\int_0^T\sigma^3\|\sqrt{\rho}\dot{u}\|_{L^2}^2\|P-P(\wi{\rho})\|_{L^3}^2
+\frac{CP(\wi{\rho})^2E_0}{(2\mu+\lambda)^2\mu}\nonumber\\[2mm]
&\le &\displaystyle\frac{C(\gamma-1)^\frac{1}{4}E_0}{(2\mu+\lambda)^3}+ \frac{1}{4}\int_0^T\sigma^3\|\nabla u\|_{L^4}^4
+\frac{1}{(2\mu+\lambda)^3}\int_0^T\sigma^3\|\sqrt{\rho}\dot{u}\|_{L^2}^3\|\nabla u\|_{L^2}\nonumber\\[2mm]
&&+\frac{1}{(2\mu+\lambda)^4}\sup_{0 \le t \le T}\left(\|P-P(\wi{\rho})\|_{L^3}^2\right)\int_0^T\sigma^3\|\sqrt{\rho}\dot{u}\|_{L^2}^2
+\frac{CP(\wi{\rho})^2E_0}{(2\mu+\lambda)^2\mu}\nonumber\\[2mm]
&\le & \displaystyle\frac{C(\gamma-1)^\frac{1}{4}E_0}{(2\mu+\lambda)^3}+\frac{1}{4}\int_0^T\sigma^3\|\nabla u\|_{L^4}^4
+\frac{1}{(2\mu+\lambda)^4}\sup_{0 \le t \le T}\left(\|P-P(\wi{\rho})\|_{L^3}^2\right)\int_0^T\sigma^3\|\sqrt{\rho}\dot{u}\|_{L^2}^2\nonumber\\[2mm]
&&+\frac{1}{(2\mu+\lambda)^3}\sup_{0 \le t \le T}\left(\sigma^\frac{1}{2}\|\nabla u\|_{L^2}\right)
\sup_{0 \le t \le T}\left(\sigma^\frac{3}{2}\|\sqrt{\rho}\dot{u}\|_{L^2}\right)\int_0^T\sigma\|\sqrt{\rho}\dot{u}\|_{L^2}^2+\frac{CP(\wi{\rho})^2E_0}{(2\mu+\lambda)^2\mu}\nonumber\\[2mm]
&\le &\displaystyle\frac{C(\gamma-1)^\frac{1}{4}E_0}{(2\mu+\lambda)^3}+\frac{1}{4}\int_0^T\sigma^3\|\nabla u\|_{L^4}^4
+\frac{(\gamma-1)^\frac{1}{6}E_0^\frac{2}{3}A_1(T)}{(2\mu+\lambda)^4}+\frac{A^\frac{3}{2}_1(T)A^\frac{1}{2}_2(T)}{(2\mu+\lambda)^3}\nonumber\\[2mm]
&&+\frac{CP(\wi{\rho})^2E_0}{(2\mu+\lambda)^2\mu}.
\ee
Substituting $\eqref{4d-A21-1}$-$\eqref{4d-A21-3}$ into \eqref{4d-A21} shows that
\be\label{4d-A212}
&&\displaystyle\int_0^T\int_{\mathbb{R}^3}\sigma^3|\nabla u|^4\nonumber\\[2mm]
&\le&\displaystyle\frac{C(\gamma-1)^\frac{1}{4}E_0}{\mu^3}
+\frac{CA^\frac{3}{2}_1(T)A^\frac{1}{2}_2(T)}{\mu^3}
+\frac{C(\gamma-1)^\frac{1}{6}E_0^\frac{2}{3}A_1(T)}{\mu^4}+\frac{CP(\wi{\rho})^2E_0}{\mu^3}.
\ee And also we get \be
\int_0^T\int_{\mathbb{R}^3}\sigma^3|P-P(\widetilde{\rho})|^4 &\le&
C\mu(\ga-1)^{\frac{1}{4}}E_0+C\mu
A^\frac{3}{2}_1(T)A^\frac{1}{2}_2(T)\nonumber\\
&&+C(\ga-1)^{\frac{1}{6}}E^\frac{2}{3}_0A_1(T)+C\frac{(2\mu+\lambda)^2}{\mu}
P^2(\wi{\rho})E_0. \ee
\vspace{2mm}
Next, we turn to estimate $A_1(\sigma(T))$. \eqref{4d-A1} shows that
\be\label{4d-A1t-c}
\displaystyle A_1(\sigma(T))
&\leq &\displaystyle
\frac{4}{\mu}\int_{\mathbb{R}^3}\sigma\mathrm{div}u(P-P(\wi{\rho}))
+\frac{C(\gamma-1)E_0}{\mu}+\frac{CP(\wi{\rho})E_0}{\mu}\nonumber\\[3mm]
&&+\frac{C(2\mu+\lambda)}{\mu}\int_0^{\sigma(T)} \sigma \|\nabla u\|_{L^3}^3
+\frac{C}{\mu}\int_0^{\sigma(T)}\int_{\mathbb{R}^3}\sigma |P-P(\wi{\rho})| |\nabla u|^2.
\ee
Based on Lemma \ref{2.2}, Lemma \ref{4d-le:3.3} and \eqref{a priori assumption}, the last two terms in the right hand side of  \eqref{4d-A1t-c} can be estimated as follows:
\be\label{4A1-11}
&&\displaystyle\frac{C(2\mu+\lambda)}{\mu}\int_0^{\sigma(T)} \sigma \|\nabla u\|_{L^3}^3\nonumber\\[3mm]
&\le& \displaystyle C \int_0^{\sigma(T)} \sigma \|\nabla u\|_{L^2}^\frac{3}{2}\|\nabla u\|_{L^6}^\frac{3}{2}\nonumber\\[3mm]
&\le& \displaystyle \frac{C}{\mu^\frac{3}{2}} \int_0^{\sigma(T)} \sigma \|\nabla u\|_{L^2}^\frac{3}{2}
\left(\|\sqrt{\rho}\dot{u}\|_{L^2}+\|P-P(\wi{\rho})_{L^6}\|\right)^{\frac{3}{2}}\nonumber\\[3mm]
&\le& \displaystyle \frac{C}{\mu^\frac{3}{2}} \sup_{0 \le t \le \sigma(T)}\left(\sigma^\frac{1}{4}\|\nabla u\|_{L^2}^\frac{1}{2}\right)
\sup_{0 \le t \le \sigma(T)}\left(\|\nabla u\|_{L^2}^\frac{1}{2}\right)\left(\int_0^{\sigma(T)}\sigma\|\sqrt{\rho}\dot{u}\|_{L^2}^2\right)^\frac{3}{4}
\left(\int_0^{\sigma(T)}\|\nabla u\|_{L^2}^2\right)^\frac{1}{4}\nonumber\\[3mm]
&&+\frac{C}{\mu^\frac{3}{2}}\sup_{0 \le t \le \sigma(T)}\Big(\|P-P(\wi{\rho})\|_{L^6}\Big)^{\frac{3}{2}}
\left(\int_0^{\sigma(T)}\|\nabla u\|_{L^2}^2\right)^\frac{3}{4}
\nonumber\\[3mm]
&\le&\frac{CA_1(T)E_{13}^\frac{1}{4}}{\mu^\frac{3}{2}}\left(\int_0^{\sigma(T)}\|\nabla
u\|_{L^2}^2\right)^\frac{1}{4}
+\frac{C(\gamma-1)^\frac{1}{16}E_0^\frac{1}{4}}{\mu^\frac{3}{2}}
\left(\int_0^{\sigma(T)}\|\nabla u\|_{L^2}^2\right)^\frac{3}{4} \ee
and \be\label{4A1t-112}
\displaystyle\frac{C}{\mu}\int_0^{\sigma(T)}\int_{\mathbb{R}^3}\sigma
|P-P(\wi{\rho})| |\nabla u|^2
&\le&\displaystyle\frac{C(\bar{\rho})}{\mu}\int_0^{\sigma(T)}\int_{\mathbb{R}^3}|\nabla
u|^2. \ee Substituting \eqref{4A1-11}-\eqref{4A1t-112} into
\eqref{4d-A1t-c}, one has \be\label{4d-A1t-c1} \displaystyle
A_1(\sigma(T))&\leq &\displaystyle
\sup_{0 \le t \le \sigma(T)}\left\{\frac{4}{\mu}\int_{\mathbb{R}^3}\sigma\mathrm{div}u(P-P(\wi{\rho}))\right\}
+\frac{C(\gamma-1)E_0}{\mu}+\frac{CP(\wi{\rho})E_0}{\mu}\nonumber\\[3mm]
&&+\frac{CA_1(T)E_{13}^\frac{1}{4}}{\mu^\frac{3}{2}}\left(\int_0^{\sigma(T)}\|\nabla
u\|_{L^2}^2\right)^\frac{1}{4}
+\frac{C(\gamma-1)^\frac{1}{16}E_0^\frac{1}{4}}{\mu^\frac{3}{2}}
\left(\int_0^{\sigma(T)}\|\nabla u\|_{L^2}^2\right)^\frac{3}{4}\nonumber\\[3mm]
&&+\frac{C(\bar{\rho})}{\mu}\int_0^{\sigma(T)}\int_{\mathbb{R}^3}|\nabla u|^2.
\ee

It follows from \eqref{3d-p1} that
\be\label{4d-A1t-c12}
\displaystyle A_1(\sigma(T))&\leq &
\frac{C(\ga-1)^{\frac{7}{36}}E_0^{\frac{17}{24}}}{\mu^{\frac{31}{12}}}
+\frac{C(\ga-1)^{\frac{11}{108}}E_0^{\frac{25}{36}}}{\mu^\frac{17}{9}}+\frac{C(\gamma-1)E_0}{\mu}+\frac{CP(\wi{\rho})E_0}{\mu}\nonumber\\[3mm]
&&+\frac{CA_1(T)E_{13}^\frac{1}{4}}{\mu^\frac{3}{2}}\left(\int_0^{\sigma(T)}\|\nabla
u\|_{L^2}^2\right)^\frac{1}{4}
+\frac{C(\gamma-1)^\frac{1}{16}E_0^\frac{1}{4}}{\mu^\frac{3}{2}}
\left(\int_0^{\sigma(T)}\|\nabla u\|_{L^2}^2\right)^\frac{3}{4}\nonumber\\[3mm]
&&+\frac{C(\bar{\rho})}{\mu}\int_0^{\sigma(T)}\int_{\mathbb{R}^3}|\nabla u|^2.
\ee
Collecting \eqref{4d-u-1}, \eqref{4d-A2-c}, \eqref{4d-A21-3}, \eqref{4d-A212} and \eqref{4d-A1t-c12} implies that
\be\label{4a2}
\displaystyle A_2(T)&\leq&
\underbrace{\frac{C(\ga-1)^{\frac{7}{36}}E_0^{\frac{17}{24}}}{\mu^{\frac{31}{12}}}
+\frac{C(\ga-1)^{\frac{11}{108}}E_0^{\frac{25}{36}}}{\mu^\frac{17}{9}}
+\frac{C(\gamma-1)E_0}{\mu}+\frac{CP(\wi{\rho})E_0}{\mu}}_{N_1}\nonumber\\[2mm]
&&\displaystyle+\underbrace{\frac{CA_1(T)E_{13}^\frac{1}{4}}{\mu^\frac{3}{2}}\left(
\frac{(\ga-1)^{\frac{1}{13}}E_0^{\frac{25}{36}}E_{12}}{\mu^\frac{12}{13}}
+\frac{\wi{\rho}E_0}{\mu}
\right)^\frac{1}{4}}_{N_2}\nonumber\\[2mm]
&&\displaystyle+\underbrace{\frac{C(\gamma-1)^\frac{1}{16}E_0^\frac{1}{4}}{\mu^\frac{3}{2}}
\left(
\frac{(\ga-1)^{\frac{1}{13}}E_0^{\frac{25}{36}}E_{12}}{\mu^\frac{12}{13}}
+\frac{\wi{\rho}E_0}{\mu}
\right)^\frac{3}{4}}_{N_3}\nonumber\\[2mm]
&&+\underbrace{\frac{C(\bar{\rho})(\ga-1)^{\frac{1}{13}}E_0^{\frac{25}{36}}E_{12}}{\mu^\frac{25}{13}}
+\frac{C(\bar{\rho})\wi{\rho}E_0}{\mu^2}}_{N_4}
+\underbrace{\displaystyle\Bigg(1+\frac{1}{\mu^4}+\frac{(\gamma-1)^2}{\mu^4}+\frac{(2\mu+\lambda)^2}{\mu^4}\Bigg)
}_{N_5}\nonumber\\[2mm]
&&\times\underbrace{\Bigg(\frac{C(\gamma-1)^\frac{1}{4}E_0}{\mu}\frac{CA_{1}^{\frac{3}{2}}(T)A_{2}^{\frac{1}{2}}(T)}{\mu}
+\frac{C(\gamma-1)^\frac{1}{6}E_0^\frac{2}{3}A_1(T)}{\mu^2}
+\frac{CP(\wi{\rho})^2E_0}{\mu}\Bigg)
}_{N_6}\nonumber\\[2mm]
&&+\underbrace{\frac{CP(\wi{\rho})E_0}{\mu^2}}_{N_7}.
\ee
Next we focus on dealing with $N_1$-$N_6$. In fact, \eqref{assumpution rb} leads to
\be\label{N12}
N_1+N_2 &\le& C\left(\frac{(\gamma-1)^\frac{1}{36}E_0^\frac{1}{4}}{\mu^\frac{1}{3}}
+\frac{\wi{\rho}^\frac{1}{6}E_0^\frac{1}{4}}{\mu^\frac{1}{3}}\right)^\frac{18}{17}
\Bigg\{\frac{(\gamma-1)^\frac{101}{612}E_0^\frac{181}{408}}{\mu^\frac{455}{204}}
+\frac{(\gamma-1)^\frac{133}{1836}E_0^\frac{263}{612}}{\mu^\frac{235}{153}}\nonumber\\[2mm]
&&+\frac{(\gamma-1)^\frac{33}{34}E_0^\frac{25}{34}}{\mu^\frac{11}{17}}
+\frac{(\gamma-1)^\frac{5}{21216}E_0^\frac{13}{4896}E_{12}^\frac{1}{4}E_{13}^\frac{1}{4}}{\mu^\frac{2657}{1768}}
\Bigg\}
+C\left(\frac{\widetilde{\rho}^\frac{1}{6}E_0^\frac{1}{4}}{\mu^\frac{1}{3}}\right)^2
\left(\frac{\widetilde{\rho}E_0^\frac{3}{2}}{\mu^\frac{3}{2}}\right)^\frac{2}{9}\widetilde{\rho}^{\gamma-\frac{5}{9}}\nonumber\\[2mm]
&&+C\left(\frac{\widetilde{\rho}^\frac{1}{6}E_0^\frac{1}{4}}{\mu^\frac{1}{3}}\right)^\frac{9}{8}\frac{\widetilde{\rho}^\frac{1}{12}}{\mu^\frac{17}{12}}
\nonumber\\[2mm]
&=& C\left(\frac{(\gamma-1)^\frac{1}{36}E_0^\frac{1}{4}}{\mu^\frac{1}{3}}
+\frac{\wi{\rho}^\frac{1}{6}E_0^\frac{1}{4}}{\mu^\frac{1}{3}}\right)^\frac{18}{17}
\Bigg\{\left(\frac{(\gamma-1)E_0^{\frac{543}{202}}}{\mu^{12}}\right)^\frac{101}{612}\frac{1}{\mu^\frac{1}{4}}+\left(\frac{(\gamma-1)E_0^{\frac{789}{133}}}{\mu^{12}}\right)^\frac{133}{1836}\frac{1}{\mu^\frac{2}{3}}
\nonumber\\[2mm]
&&+\left(\frac{(\gamma-1)E_0^\frac{25}{33}}{\mu^\frac{2}{3}}\right)^\frac{33}{34}+\left(\frac{(\gamma-1)E_0^{\frac{169}{15}}}{\mu^{12}}\right)^\frac{5}{21216}
\frac{E_{12}^\frac{1}{4}E_{13}^\frac{1}{4}}{\mu^\frac{3}{2}}
+\left(\frac{\widetilde{\rho}^\frac{1}{6}E_0^\frac{1}{4}}{\mu^\frac{1}{3}}\right)^\frac{16}{17}
\left(\frac{\widetilde{\rho}E_0^\frac{3}{2}}{\mu^\frac{3}{2}}\right)^\frac{2}{9}\widetilde{\rho}^{\gamma-\frac{5}{9}}\nonumber\\[2mm]
&&+\left(\frac{\widetilde{\rho}^\frac{1}{6}E_0^\frac{1}{4}}{\mu^\frac{1}{3}}\right)^\frac{1}{136}\frac{\widetilde{\rho}^\frac{1}{12}}{\mu^\frac{17}{12}}
\Bigg\}\nonumber\\[2mm]
&\le&
C\left(\frac{(\gamma-1)^\frac{1}{36}E_0^\frac{1}{4}}{\mu^\frac{1}{3}}
+\frac{\wi{\rho}^\frac{1}{6}E_0^\frac{1}{4}}{\mu^\frac{1}{3}}\right)^\frac{18}{17}E_{18},
\ee where \be\label{s 4E14} \displaystyle
E_{18}=\frac{1}{\mu^\frac{1}{4}}+\frac{1}{\mu^\frac{2}{3}}+1
+\frac{E_{12}^\frac{1}{4}E_{13}^\frac{1}{4}}{\mu^\frac{3}{2}}
+\widetilde{\rho}^{\gamma-\frac{5}{9}}+\frac{\widetilde{\rho}^\frac{1}{12}}{\mu^\frac{17}{12}}.\ee

\be\label{N34}
N_3+N_4&\le& C\left(\frac{(\gamma-1)^\frac{1}{36}E_0^\frac{1}{4}}{\mu^\frac{1}{3}}
+\frac{\wi{\rho}^\frac{1}{6}E_0^\frac{1}{4}}{\mu^\frac{1}{3}}\right)^\frac{18}{17}
\Bigg\{\frac{(\gamma-1)^\frac{321}{3536}E_0^\frac{413}{816}}{\mu^\frac{813}{442}}
+\frac{(\gamma-1)^\frac{21}{442}E_0^\frac{263}{612}}{\mu^\frac{347}{221}}\nonumber\\[2mm]
&&+\left(\frac{(\gamma-1)E_0^4}{\mu^{12}}\right)^\frac{1}{16}
\left(\frac{\wi{\rho}^\frac{1}{6}E_0^\frac{1}{4}}{\mu^\frac{1}{3}}\right)^\frac{33}{17}
\frac{\widetilde{\rho}^\frac{1}{4}}{\mu^\frac{1}{2}}
+\left(\frac{\wi{\rho}^\frac{1}{6}E_0^\frac{1}{4}}{\mu^\frac{1}{3}}\right)^\frac{50}{17}
\frac{\widetilde{\rho}^\frac{1}{3}}{\mu^\frac{2}{3}}
\Bigg\}\nonumber\\[2mm]
&=&C\left(\frac{(\gamma-1)^\frac{1}{36}E_0^\frac{1}{4}}{\mu^\frac{1}{3}}
+\frac{\wi{\rho}^\frac{1}{6}E_0^\frac{1}{4}}{\mu^\frac{1}{3}}\right)^\frac{18}{17}
\Bigg\{
\left(\frac{(\gamma-1)E_0^{\frac{5369}{963}}}{\mu^{12}}\right)^\frac{321}{3536}\frac{1}{\mu^\frac{3}{4}}\nonumber\\[2mm]
&&+\left(\frac{(\gamma-1)E_0^{\frac{263}{612}\times\frac{442}{21}}}{\mu^{12}}\right)^\frac{21}{442}\frac{1}{\mu}
+\frac{\widetilde{\rho}^\frac{1}{4}}{\mu^\frac{1}{2}}
+\frac{\widetilde{\rho}^\frac{1}{3}}{\mu^\frac{2}{3}}
\Bigg\}\nonumber\\[2mm]
&\le&
C\left(\frac{(\gamma-1)^\frac{1}{36}E_0^\frac{1}{4}}{\mu^\frac{1}{3}}
+\frac{\wi{\rho}^\frac{1}{6}E_0^\frac{1}{4}}{\mu^\frac{1}{3}}\right)^\frac{18}{17}E_{19},
\ee where \be\label{s 4E15} \displaystyle
E_{19}=\frac{1}{\mu^\frac{3}{4}}+\frac{1}{\mu}+\frac{\widetilde{\rho}^\frac{1}{4}}{\mu^\frac{1}{2}}
+\frac{\widetilde{\rho}^\frac{1}{3}}{\mu^\frac{2}{3}}.\ee

\be\label{N45}
N_5\times N_6+N_7&\le& C\Bigg(1+\frac{1}{\mu^4}+\frac{(\gamma-1)^2}{\mu^4}+\frac{(2\mu+\lambda)^2}{\mu^4}\Bigg)
\left(\frac{(\gamma-1)^\frac{1}{36}E_0^\frac{1}{4}}{\mu^\frac{1}{3}}
+\frac{\wi{\rho}^\frac{1}{6}E_0^\frac{1}{4}}{\mu^\frac{1}{3}}\right)^\frac{18}{17}\nonumber\\[2mm]
&&\times\Bigg\{\frac{(\gamma-1)^\frac{15}{68}E_0^\frac{25}{34}}{\mu^\frac{11}{17}}
+\frac{(\gamma-1)^{\frac{1}{9792}}E_0^{\frac{1}{1088}}}{\mu^\frac{817}{816}}
+\frac{(\gamma-1)^\frac{7}{51}E_0^\frac{41}{102}}{\mu^\frac{28}{17}}\nonumber\\[2mm]
&&+\left(\frac{\wi{\rho}^\frac{1}{6}E_0^\frac{1}{4}}{\mu^\frac{1}{3}}\right)^\frac{16}{17}
\left(\frac{\widetilde{\rho}E_0^\frac{3}{2}}{\mu}\right)^\frac{1}{3}\widetilde{\rho}^{2\gamma-\frac{2}{3}}
\Bigg\}
+C\left(\frac{\wi{\rho}^\frac{1}{6}E_0^\frac{1}{4}}{\mu^\frac{1}{3}}\right)^4\frac{\widetilde{\rho}^{\gamma-\frac{2}{3}}}{\mu^\frac{2}{3}}\nonumber\\[2mm]
&=&C\Bigg(1+\frac{1}{\mu^4}+\frac{(\gamma-1)^2}{\mu^4}+\frac{(2\mu+\lambda)^2}{\mu^4}\Bigg)
\left(\frac{(\gamma-1)^\frac{1}{36}E_0^\frac{1}{4}}{\mu^\frac{1}{3}}
+\frac{\wi{\rho}^\frac{1}{6}E_0^\frac{1}{4}}{\mu^\frac{1}{3}}\right)^\frac{18}{17}\nonumber\\[2mm]
&&\times\Bigg\{
\left(\frac{(\gamma-1)E_0^\frac{10}{3}}{\mu^\frac{44}{15}}\right)^\frac{15}{68}
+\left(\frac{(\gamma-1)E_0^9}{\mu^{12}}\right)^{\frac{1}{36}\times\frac{1}{272}}\frac{1}{\mu}
+\left(\frac{(\gamma-1)E_0^\frac{41}{14}}{\mu^{12}}\right)^{\frac{7}{51}}
\Bigg\}\nonumber\\[2mm]
&&+C\left(\frac{(\gamma-1)^\frac{1}{36}E_0^\frac{1}{4}}{\mu^\frac{1}{3}}
+\frac{\wi{\rho}^\frac{1}{6}E_0^\frac{1}{4}}{\mu^\frac{1}{3}}\right)^\frac{18}{17}
\left(\frac{\wi{\rho}^\frac{1}{6}E_0^\frac{1}{4}}{\mu^\frac{1}{3}}\right)^\frac{50}{17}\frac{\widetilde{\rho}^{\gamma-\frac{2}{3}}}{\mu^\frac{2}{3}}\nonumber\\[2mm]
&\le&
C\left(\frac{(\gamma-1)^\frac{1}{36}E_0^\frac{1}{4}}{\mu^\frac{1}{3}}
+\frac{\wi{\rho}^\frac{1}{6}E_0^\frac{1}{4}}{\mu^\frac{1}{3}}\right)^\frac{18}{17}E_{20},
\ee \be\label{s 4E16} \displaystyle
E_{20}=\Bigg(1+\frac{1}{\mu^4}+\frac{(\gamma-1)^2}{\mu^4}+\frac{(2\mu+\lambda)^2}{\mu^4}\Bigg)
\left(1+\frac{1}{\mu}\right)+\frac{\widetilde{\rho}^{\gamma-\frac{2}{3}}}{\mu^\frac{2}{3}}.\ee
\vspace{3mm}
It thus follows \eqref{N12}, \eqref{N34} and \eqref{N45} that
\be
\displaystyle A_2(T)&\leq&C\left(\frac{(\gamma-1)^\frac{1}{36}E_0^\frac{1}{4}}{\mu^\frac{1}{3}}
+\frac{\wi{\rho}^\frac{1}{6}E_0^\frac{1}{4}}{\mu^\frac{1}{3}}\right)^\frac{18}{17}
(E_{18}+E_{19}+E_{20})\nonumber\\[3mm]
&\le&\frac{(\gamma-1)^\frac{1}{36}E_0^\frac{1}{4}}{\mu^\frac{1}{3}}
+\frac{\wi{\rho}^\frac{1}{6}E_0^\frac{1}{4}}{\mu^\frac{1}{3}}, \ee
provided $\displaystyle
\left(\frac{(\ga-1)^{\frac{1}{36}}E_0^{\frac{1}{4}}}{\mu^\frac{1}{3}}
+\frac{\wi{\rho}^\frac{1}{6}E_0^\frac{1}{4}}{\mu^\frac{1}{3}}\right)
\le \Big(C(E_{18}+E_{19}+E_{20})\Big)^{-17}$.

\vspace{3mm}

Finally, to finish the proof of Lemma \ref{3d-le:4.7}, it remains to prove \eqref{3d-A1-le}. With \eqref{4d-A1} and \eqref {4d-A1t-c1} at hand, we just have to estimate the terms $\displaystyle\frac{2\mu+\lambda}{\mu}\int_{\sigma(T)}^T \sigma \|\nabla u\|_{L^3}^3$ and $\displaystyle\frac{C}{\mu}\int_{\sigma(T)}^T\int_{\mathbb{R}^3}\sigma |P-P(\wi{\rho})| |\nabla u|^2$. By H\"older inequality, we have
\be\label{4A1t-12}
\displaystyle\frac{2\mu+\lambda}{\mu}\int_{\sigma(T)}^T \sigma \|\nabla u\|_{L^3}^3
&\le& \displaystyle C\left(\int_{\sigma(T)}^T\|\nabla u\|_{L^2}^2\right)^\frac{1}{2}
\left(\int_{\sigma(T)}^T\|\nabla u\|_{L^4}^4\right)^\frac{1}{2}\nonumber\\[3mm]
&\le& \displaystyle \left(\frac{CE_0}{\mu}\right)^\frac{1}{2} \left(\int_{\sigma(T)}^T\|\nabla u\|_{L^4}^4\right)^\frac{1}{2}
\ee
and
\be\label{4A1-112}
&&\frac{C}{\mu}\int_{\sigma(T)}^T\int_{\mathbb{R}^3}\sigma |P-P(\wi{\rho})| |\nabla u|^2\nonumber\\[3mm]
&\le&
\frac{C}{\mu}\left(\int_{\sigma(T)}^T\int_{\mathbb{R}^3} |P-P(\wi{\rho})|^4\right)^\frac{1}{4}
\left(\int_{\sigma(T)}^T\int_{\mathbb{R}^3} |\nabla u|^4\right)^\frac{1}{4}
\left(\int_{\sigma(T)}^T\int_{\mathbb{R}^3} |\nabla u|^2\right)^\frac{1}{2}.
\ee
It follows from Lemma \ref{4.2}, \eqref{4dP}, \eqref{4d-A212}, Lemma \ref{4.3}, \eqref{4A1t-12} and \eqref{4A1-112} that
\be\label{four A one}
A_1(T)&\le& A_1(\sigma(T))+\frac{C(2\mu+\lambda)}{\mu}\int_{\sigma(T)}^T \sigma \|\nabla u\|_{L^3}^3
+\frac{C}{\mu}\int_{\sigma(T)}^T\int_{\mathbb{R}^3}\sigma |P-P(\wi{\rho})| |\nabla u|^2\nonumber\\[2mm]
&\le & \displaystyle
\underbrace{\frac{C(\ga-1)^{\frac{7}{36}}E_0^{\frac{17}{24}}}{\mu^{\frac{31}{12}}}
+\frac{C(\ga-1)^{\frac{11}{108}}E_0^{\frac{25}{36}}}{\mu^\frac{17}{9}}
+\frac{C(\gamma-1)E_0}{\mu}+\frac{CP(\wi{\rho})E_0}{\mu}}_{N_8}\nonumber\\[2mm]
&&\displaystyle\underbrace{+\frac{CA_1(T)E_{13}^\frac{1}{4}}{\mu^\frac{3}{2}}\left(
\frac{(\ga-1)^{\frac{1}{13}}E_0^{\frac{25}{36}}E_{12}}{\mu^\frac{12}{13}}
+\frac{\wi{\rho}E_0}{\mu}
\right)^\frac{1}{4}}_{N_9}\nonumber\\[2mm]
&&\displaystyle\underbrace{+\frac{C(\gamma-1)^\frac{1}{16}E_0^\frac{1}{4}}{\mu^\frac{3}{2}}
\left(
\frac{(\ga-1)^{\frac{1}{13}}E_0^{\frac{25}{36}}E_{12}}{\mu^\frac{12}{13}}
+\frac{\wi{\rho}E_0}{\mu}
\right)^\frac{3}{4}}_{N_{10}}\nonumber\\[2mm]
&&
\displaystyle\underbrace{+\displaystyle\frac{C(\bar{\rho},\wi{\rho})}{\mu}
\left(
\frac{(\ga-1)^{\frac{1}{13}}E_0^{\frac{25}{36}}E_{12}}{\mu^\frac{12}{13}}
+\frac{\wi{\rho}E_0}{\mu}
\right)}_{N_{11}}\nonumber\\[2mm]
&&\displaystyle\underbrace{+\frac{CE_0^\frac{1}{2}}{\mu^\frac{1}{2}}
\left(\frac{(\gamma-1)^\frac{1}{4}E_0}{\mu^3}
+\frac{A^\frac{3}{2}_1(T)A^\frac{1}{2}_2(T)}{\mu^3}
+\frac{(\gamma-1)^\frac{1}{6}E_0^\frac{2}{3}A_1(T)}{\mu^4}+\frac{P(\wi{\rho})^2E_0}{\mu^3}\right)^\frac{1}{2}}
_{N_{12}}.
\ee
Following a process similar to $N_1$-$N_7$, one gets $N_8$-$N_{12}$ as follows:
\be\label{N78}
N_8+N_9=N_1+N_2&\le& C\left(\frac{(\gamma-1)^\frac{1}{36}E_0^\frac{1}{4}}{\mu^\frac{1}{3}}
+\frac{\wi{\rho}^\frac{1}{6}E_0^\frac{1}{4}}{\mu^\frac{1}{3}}\right)^\frac{18}{17}E_{18}\nonumber\\[2mm]
&\le&
C\left(\frac{(\gamma-1)^\frac{1}{36}E_0^\frac{1}{4}}{\mu^\frac{1}{3}}
+\frac{\wi{\rho}^\frac{1}{6}E_0^\frac{1}{4}}{\mu^\frac{1}{3}}\right)^\frac{1}{2}E_{18},
\ee
\be\label{N910} N_{10}+N_{11}=N_3+N_4&\le&
C\left(\frac{(\gamma-1)^\frac{1}{36}E_0^\frac{1}{4}}{\mu^\frac{1}{3}}
+\frac{\wi{\rho}^\frac{1}{6}E_0^\frac{1}{4}}{\mu^\frac{1}{3}}\right)^\frac{18}{17}E_{19}\nonumber\\[2mm]
&\le&
C\left(\frac{(\gamma-1)^\frac{1}{36}E_0^\frac{1}{4}}{\mu^\frac{1}{3}}
+\frac{\wi{\rho}^\frac{1}{6}E_0^\frac{1}{4}}{\mu^\frac{1}{3}}\right)^\frac{1}{2}E_{19}
\ee
and
\be\label{N11} N_{12}
&\le&
C\left(\frac{(\gamma-1)^\frac{1}{36}E_0^\frac{1}{4}}{\mu^\frac{1}{3}}
+\frac{\wi{\rho}^\frac{1}{6}E_0^\frac{1}{4}}{\mu^\frac{1}{3}}\right)^\frac{1}{2}
\left\{\frac{(\gamma-1)^\frac{1}{9}E_0^\frac{7}{8}}{\mu^\frac{11}{6}}\right.\nonumber\\[2mm]
&&\left.+\left(\frac{(\gamma-1)^\frac{1}{36}E_0^\frac{1}{4}}{\mu^\frac{1}{3}}
+\frac{\wi{\rho}^\frac{1}{6}E_0^\frac{1}{4}}{\mu^\frac{1}{3}}\right)^\frac{1}{32}\frac{E_0^\frac{1}{2}}{\mu^2}
+\frac{(\gamma-1)^\frac{5}{72}E_0^\frac{17}{24}A_1(T)}{\mu^\frac{7}{3}}
+\frac{\widetilde{\rho}^{\gamma-\frac{1}{12}}E_0^\frac{7}{8}}{\mu^\frac{11}{6}}
\right\}
\nonumber\\[2mm]
&\le&C\left(\frac{(\gamma-1)^\frac{1}{36}E_0^\frac{1}{4}}{\mu^\frac{1}{3}}
+\frac{\wi{\rho}^\frac{1}{6}E_0^\frac{1}{4}}{\mu^\frac{1}{3}}\right)^\frac{1}{2}
\left\{\left(\frac{(\gamma-1)E_0^{\frac{63}{8}}}{\mu^{12}}\right)^\frac{1}{9}\frac{1}{\mu^\frac{1}{2}}
+\left(\frac{\widetilde{\rho}}{2C}\right)^\frac{1}{32}\frac{E_0^\frac{1}{2}}{\mu^2}\right.\nonumber\\[2mm]
&&\left.+\left(\frac{\wi{\rho}^\frac{1}{6}E_0^\frac{1}{4}}{\mu^\frac{1}{3}}\right)^\frac{1}{32}\frac{E_0^\frac{1}{2}}{\mu^2}
+\left(\frac{(\gamma-1)E_0^{\frac{51}{5}}}{\mu^{12}}\right)^\frac{5}{72}
\frac{1}{\mu^\frac{1}{2}}
+\left(\frac{\widetilde{\rho}E_0^\frac{3}{2}}{\mu^\frac{22}{7}}\right)^\frac{7}{12}\widetilde{\rho}^{\gamma-\frac{2}{3}}\right.
\Bigg\}\nonumber\\[2mm]
&\le& C\left(\frac{(\gamma-1)^\frac{1}{36}E_0^\frac{1}{4}}{\mu^\frac{1}{3}}
+\frac{\wi{\rho}^\frac{1}{6}E_0^\frac{1}{4}}{\mu^\frac{1}{3}}\right)^\frac{1}{2}
\left\{\frac{1}{\mu^\frac{1}{2}}
+\left(\frac{\widetilde{\rho}E_0^{16}}{\mu^{12}}\right)^\frac{1}{32}\frac{1}{\mu^\frac{13}{8}}
+\left(\frac{\widetilde{\rho}E_0^{\frac{195}{2}}}{\mu^2}\right)^\frac{1}{192}\frac{1}{\mu^2}
+\widetilde{\rho}^{\gamma-\frac{2}{3}}\right\}\nonumber\\[2mm]
&\le&
C\left(\frac{(\gamma-1)^\frac{1}{36}E_0^\frac{1}{4}}{\mu^\frac{1}{3}}
+\frac{\wi{\rho}^\frac{1}{6}E_0^\frac{1}{4}}{\mu^\frac{1}{3}}\right)^\frac{1}{2}
E_{21}, \ee where \be\label{s 4E17} \displaystyle
E_{21}=\frac{1}{\mu^\frac{1}{2}} +\frac{1}{\mu^\frac{13}{8}}
+\frac{1}{\mu^2} +\widetilde{\rho}^{\gamma-\frac{2}{3}}, \ee and we
have used
$\displaystyle\frac{(\gamma-1)^\frac{1}{36}E_0^\frac{1}{4}}{\mu^\frac{1}{3}}
\leq \frac{\widetilde{\rho}}{2C}$.

Substituting \eqref{N78}, \eqref{N910} and \eqref{N11} into \eqref{four A one},
we obtain
\be
A_1(T)&\le &\displaystyle C\left(\frac{(\ga-1)^{\frac{1}{36}}E_0^{\frac{1}{4}}}{\mu^\frac{1}{3}}
+\frac{\wi{\rho}^\frac{1}{6}E_0^\frac{1}{4}}{\mu^\frac{1}{3}}\right)^\frac{1}{2}(E_{18}+E_{19}+E_{21})\nonumber\\[2mm]
&\le
&\displaystyle\left(\frac{(\ga-1)^{\frac{1}{36}}E_0^{\frac{1}{4}}}{\mu^\frac{1}{3}}
+\frac{\wi{\rho}^\frac{1}{6}E_0^\frac{1}{4}}{\mu^\frac{1}{3}}\right)^\frac{3}{8},
\ee provided $\displaystyle
\left(\frac{(\ga-1)^{\frac{1}{36}}E_0^{\frac{1}{4}}}{\mu^\frac{1}{3}}
+\frac{\wi{\rho}^\frac{1}{6}E_0^\frac{1}{4}}{\mu^\frac{1}{3}}\right)
\le \Big(C(E_{18}+E_{19}+E_{21})\Big)^{-8}$.
\endpf

Now we are in a position to close the $a \ priori$ assumption on $\rho$.
\begin{lemma}\label{4d:rho}
Under the conditions of Proposition \ref{prop 4.1}, it holds that
\be\label{4d-upper bound of rho} \sup_{0\leq t\leq
T}\|\rho\|_{L^\infty}\leq\frac{7\bar{\rho}}{4} \ee for any
$(x,t)\in\mathbb{R}^3\times[0,T]$, provided \bex
\begin{aligned}\frac{\left((\gamma-1)^{\frac{1}{36}}+\wi{\rho}^\frac{1}{6}\right)E_0^\frac{1}{4}}{\mu^\frac{1}{3}}
\leq\varepsilon\triangleq
\min\left\{\varepsilon_6,(2C(\bar{\rho},M))^{-\frac{16}{3}}\mu^{4},(4C(\bar{\rho}))^{-2}
\right\}.\end{aligned} \eex
\end{lemma}
\pf
In fact, the proof is similar to the one in Lemma \ref{3d-le:rho}, then we just list some
differences.
Here we rewrite \eqref{b1} as follows:
\be
\displaystyle |b(t_2)-b(t_1)|&\leq&\displaystyle \frac{ C(\bar{\rho})}{\mu^\frac{3}{4}}\left(\frac{(\ga-1)^{\frac{1}{36}}E_0^{\frac{1}{4}}}{\mu^\frac{1}{3}}
+\frac{\wi{\rho}^\frac{1}{6}E_0^\frac{1}{4}}{\mu^\frac{1}{3}}\right)^{\frac{3}{16}}.
\ee
Therefore, for $t\in [0,\sigma(T)]$, one can choose $N_0$ and $N_{1}$ in Lemma \ref{zlo} as follows:
$$
N_{1}=0,\ \ \ \ \
\ N_0=\frac{ C(\bar{\rho})}{\mu^\frac{3}{4}}\left(\frac{(\ga-1)^{\frac{1}{36}}E_0^{\frac{1}{4}}}{\mu^\frac{1}{3}}
+\frac{\wi{\rho}^\frac{1}{6}E_0^\frac{1}{4}}{\mu^\frac{1}{3}}\right)^{\frac{3}{16}},
$$
and $\bar{\zeta}=0$. Then
$$
g(\zeta)=-\frac{\zeta P(\zeta)}{2\mu+\lambda}\leq -N_1=0  \ \  \mathrm{for\ \ all} \ \ \zeta\geq\bar{\zeta} =0.
$$
Thus \be\label{md3} \sup_{0\leq t \leq
\sigma(T)}\|\rho\|_{L^\infty}\leq
\max\{\bar{\rho},0\}+N_0\leq\bar{\rho}
+\frac{C(\bar{\rho})}{\mu^\frac{3}{4}}\Big(\frac{(\ga-1)^{\frac{1}{36}}E_0^{\frac{1}{4}}}{\mu^\frac{1}{3}}
+\frac{\wi{\rho}^\frac{1}{6}E_0^\frac{1}{4}}{\mu^\frac{1}{3}}\Big)^{\frac{3}{16}}\leq\frac{3\bar{\rho}}{2},
\ee provided \be\label{a3}
\frac{(\ga-1)^{\frac{1}{36}}E_0^{\frac{1}{4}}}{\mu^\frac{1}{3}}
+\frac{\wi{\rho}^\frac{1}{6}E_0^\frac{1}{4}}{\mu^\frac{1}{3}} \leq
\min\left\{\varepsilon_6,
(2C(\bar{\rho},M))^{-\frac{16}{3}}\mu^{4}\right\}. \ee
\vspace{2mm}
On the other hand, for $t\in [\sigma(T), T]$, we can rewrite
\eqref{b2} as follows:
 \be\label{b4} \displaystyle|b(t_2)-b(t_1)|
&\leq&\displaystyle \frac{1}{2\mu+\lambda}(t_2-t_1)
+C(\bar{\rho})\left(\frac{(\ga-1)^{\frac{1}{36}}E_0^{\frac{1}{4}}}{\mu^\frac{1}{3}}
+\frac{\wi{\rho}^\frac{1}{6}E_0^\frac{1}{4}}{\mu^\frac{1}{3}}\right)^2,
\ee provided $\displaystyle
\frac{(\ga-1)^{\frac{1}{36}}E_0^{\frac{1}{4}}}{\mu^\frac{1}{3}}
+\frac{\wi{\rho}^\frac{1}{6}E_0^\frac{1}{4}}{\mu^\frac{1}{3}}\leq
\varepsilon_6$. Therefore, one can choose $N_1$ and $N_0$ in Lemma
\ref{zlo} as
$$
N_1=\frac{1}{2\mu+\lambda},\ \ \ \ \ \ N_0=C(\bar{\rho})\left(\frac{(\ga-1)^{\frac{1}{36}}E_0^{\frac{1}{4}}}{\mu^\frac{1}{3}}
+\frac{\wi{\rho}^\frac{1}{6}E_0^\frac{1}{4}}{\mu^\frac{1}{3}}\right)^2.
$$
Note that
$$
g(\zeta)=-\frac{\zeta P(\zeta)}{2\mu+\lambda}\leq -N_1=-\frac{1}{2\mu+\lambda} \ \  \mathrm{for\ \ all} \ \ \zeta\geq 1,
$$
one can set $\bar{\zeta}=1$. Thus
\be\label{md4}
\begin{aligned}[b]
\sup_{\sigma(T)\leq s \leq T}\|\rho\|_{L^\infty}\leq \max\left\{\frac{3}{2}\bar{\rho},1\right\}
+N_0 &\leq \frac{3}{2}\bar{\rho}+C(\bar{\rho})\left(\frac{(\ga-1)^{\frac{1}{36}}E_0^{\frac{1}{4}}}{\mu^\frac{1}{3}}
+\frac{\wi{\rho}^\frac{1}{6}E_0^\frac{1}{4}}{\mu^\frac{1}{3}}\right)^2\\[3mm]
&\le \frac{7\bar{\rho}}{4},
\end{aligned}
\ee provided \be\label{a22}
\frac{(\ga-1)^{\frac{1}{36}}E_0^{\frac{1}{4}}}{\mu^\frac{1}{3}}
+\frac{\wi{\rho}^\frac{1}{6}E_0^\frac{1}{4}}{\mu^\frac{1}{3}} \leq
\min\Big\{\varepsilon_6, (2C(\bar{\rho},M))^{-\frac{16}{3}}\mu^{4},
(4C(\bar{\rho}))^{-2}\Big\}. \ee The combination of \eqref{md3} and
\eqref{md4} completes the proof of Lemma \ref{4d:rho}.
\endpf

\bigbreak Now, the proof of Proposition \ref{prop 4.1} is completed.
Next, following a process similar to that in the proof of Theorem
\ref{3d-th:1.1}, we can prove that the results obtained in
Proposition \ref{prop 4.1} still hold in the case of $\gamma \geq
2$.  At last, we will derive the time-dependent higher norm
estimates of the smooth solution $(\rho, u)$. In fact, the proofs are
the same as the ones in \cite{Huang-Li-Xin}. For the
convenience, we omit them here.

\section*{Acknowledgements}
This work was
supported  by the National Natural Science Foundation of China
$\#$11331005, the Program for Changjiang Scholars and Innovative
Research Team in University $\#$IRT13066.

\end{document}